\documentclass[12pt,a4paper]{article}
\usepackage[english]{babel}
\usepackage[utf8]{inputenc}
\usepackage{amsmath}
\usepackage{amsfonts}
\usepackage{amssymb}
\usepackage{graphicx}
\usepackage{listings}
\usepackage{listingsutf8}
\usepackage{tikz}
\usepackage{cite}
\usepackage{algorithm2e}
\usepackage{array}
\begin{document}

\setcounter{tocdepth}{2}

\nocite{*}
\newcounter{currentcounter}

\newcounter{example}
\newcommand\exc{\theexample}

\newenvironment{general}[1]{\refstepcounter{currentcounter} \bigskip \noindent \textbf{#1 \thecurrentcounter :}\itshape}{\normalfont \medskip}

\newenvironment{generalbis}[2]{\refstepcounter{currentcounter} \bigskip \noindent \textbf{#1 \thesubsection.\thecurrentcounter ~#2 :}\itshape}{\normalfont \medskip}

\newenvironment{proof}{\noindent \textbf{Proof :} \newline \noindent \hspace*{0.2cm}}{\hspace*{\fill}$\square$ \bigskip}

\newenvironment{proofbis}[1]{\noindent \textbf{#1:} \newline \noindent \hspace*{0.2cm}}{\newline \hspace*{\fill}$\square$ \newline}

\newenvironment{famous}[1]{\medskip \noindent \textbf{#1:} \newline \noindent  \itshape}{\normalfont \medskip}

\newenvironment{court}[1]{\refstepcounter{currentcounter} \smallskip \noindent \textbf{#1 \thesubsection.\thecurrentcounter :} \newline \noindent \itshape}{\normalfont \smallskip}

\newcolumntype{C}{>{$}c<{$}}
\newcommand\toline{\smallskip \newline}
\newcommand\refp[1]{(\ref{#1})}
\renewcommand\mod{\mathrm{~mod~}}
\renewcommand\gcd[2]{\mathrm{gcd}(#1,#2)}
\newcommand\lcm[2]{\mathrm{lcm}(#1,#2)}
\newcommand\naturaliso{\cong}
\newcommand\iso{\simeq}
\newcommand\cross[1]{#1^{\times}}
\newcommand\crosslong[1]{\cross{(#1)}}
\newcommand\poldegree[1]{\mathrm{deg}(#1)}
\newcommand\cardinal[1]{\# \left(#1\right)}
\newcommand\cardinalshort[1]{\# #1}
\renewcommand\det{\mathrm{det}}

\newcommand\signature{\mathrm{sgn}}

\newcommand\bb[1]{\mathbb{#1}}
\newcommand\kk{\bb{K}}
\newcommand\ok{\mathcal{O}_{\kk}}
\newcommand\zz{\bb{Z}}
\newcommand\qq{\bb{Q}}
\newcommand\rr{\bb{R}}
\newcommand\cc{\bb{C}}
\newcommand\goth[1]{\mathfrak{#1}}
\newcommand\zsz[1]{\zz/#1\zz}
\newcommand\zszcross[1]{\zsz{#1}^{\times}}
\newcommand\dualsimple[1]{#1^{\vee}}
\newcommand\dual[1]{(#1)^{\vee}}

\newcommand\dx[1]{\mathrm{d}#1}
\newcommand\ddx[1]{\frac{\mathrm{d}}{\dx{#1}}}
\newcommand\dkdxk[2]{\frac{\mathrm{d}^{#2}}{\dx{#1}^{#2}}}
\newcommand\partialx[1]{\partial #1}
\newcommand\partialdx[1]{\frac{\partial}{\partialx{#1}}}
\newcommand\partialkdxk[2]{\frac{\partial^{#2}}{\partialx{#1}^{#2}}}
\newcommand\integral[3]{\int_{#1}{#2}\dx{#3}}
\newcommand\integralsimple[1]{\integrale{G}{#1}{\lambda}}
\newcommand\sprod[2]{\langle #1, #2 \rangle}

\newcommand\sign[1]{sg\left(#1 \right)}
\newcommand\limi[1]{\lim_{#1 \to + \infty}}

\newcommand\sln[2]{\mathrm{SL}_{#1}(#2)}
\newcommand\slnz[1]{\sln{#1}{\zz}}
\newcommand\gln[2]{\mathrm{GL}_{#1}(#2)}
\newcommand\glnz[1]{\gln{#1}{\zz}}

\newcommand\bars[1]{\underline{#1}}
\newcommand\taubar{\bars{\tau}}
\newcommand\qbar{\bars{q}}
\newcommand\sigmabar{\bars{\sigma}}
\newcommand\rbar{\bars{r}}
\newcommand\nbar{\bars{n}}
\newcommand\mbar{\bars{m}}
\newcommand\mubar{\bars{\mu}}
\newcommand\xbar{\bars{x}}
\newcommand\Xbar{\bars{X}}
\newcommand\abar{\bars{a}}
\newcommand\bbar{\bars{b}}
\newcommand\alphabar{\bars{\alpha}}
\newcommand\omegabar{\bars{\omega}}
\newcommand\taubarsj[1]{\taubar^{-}(#1)}
\newcommand\onebar{\bars{1}}

\newcommand\hh{\bb{H}}
\newenvironment{psmallmatrix}
  {\left(\begin{smallmatrix}}
  {\end{smallmatrix}\right)}
  
\newcommand\declarefunction[5]{#1 := \begin{cases}
\hfill #2 \hfill & \to ~ #3 \\
\hfill #4 \hfill & \to ~ #5 \\
\end{cases}}

\newcommand\sume{\sideset{}{_e}\sum}

\newcommand\tendstowhen[1]{\xrightarrow[#1]{}}

\newcommand\pending{$\bigskip \newline \blacktriangle \blacktriangle \blacktriangle \blacktriangle \blacktriangle \blacktriangle \blacktriangle \blacktriangle \bigskip \newline \hfill$}

\newcommand\pendingref{\textbf{[?]}}

\newcommand\opc[1]{\mathcal{O}^{+,\times}_{#1}}
\newcommand\opcf{\opc{\goth{f}}}
\newcommand\eps{\varepsilon}
\newcommand\units{\ok^{\times}}
\newcommand\norm[1]{\mathcal{N}(#1)}

\newcommand\ta{\tilde{a}}
\newcommand\tabar{\tilde{\abar}}
\newcommand\talpha{\tilde{\alpha}}
\newcommand\talphabar{\tilde{\alphabar}}
\newcommand\ts{\tilde{s}}
\newcommand\ttt{\tilde{t}}
\newcommand\tlambda{\tilde{\lambda}}
\newcommand\tc{\tilde{c}}
\newcommand\tcone{\tc_1}
\newcommand\tctwo{\tc_2}
\newcommand\te{\tilde{e}}
\newcommand\teone{\te_1}
\newcommand\tetwo{\te_2}
\newcommand\tmu{\tilde{\mu}}
\newcommand\tA{\tilde{A}}

\newcommand\spvgamma{\Gamma_{\goth{f}, \goth{b}, \goth{a}}(\eps; h)}
\newcommand\spvgammasign{\Gamma_{\goth{f}, \goth{b}, \goth{a}}^{\pm}(\eps; h)}
\newcommand\spvgr{G_{r, \goth{f}, \goth{b}, \goth{a}}^{\pm}(u_1, \dots, u_r; h)}
\newcommand\spvgrrho{G_{r, \goth{f}, \goth{b}, \goth{a}}^{\pm}([\eps_{\rho(1)}|\dots|\eps_{\rho(r)}]; h_\rho)}
\newcommand\spvgrcomplete{I_{r,\goth{f}, \goth{b}, \goth{a}}(\eps_1, \dots, \eps_r ; \bars{h}, \bars{\mu}, \bars{\nu})}

\newcommand\spvgammac{\Gamma_{\goth{f}, \goth{b}, \goth{a}}(\eps; h, \sigma_{\cc})}
\newcommand\spvgrc{G_{r, \goth{f}, \goth{b}, \goth{a}}^{\pm}([\eps_{\rho(1)}|\dots|\eps_{\rho(r)}]; h_\rho, \sigma_{\cc})}
\newcommand\spvgrrhoc{G_{r, \goth{f}, \goth{b}, \goth{a}}^{\pm}([\eps_{\rho(1)}|\dots|\eps_{\rho(r)}]; h_\rho)}
\newcommand\spvgrcompletec{I_{r,\goth{f}, \goth{b}, \goth{a}}(\eps_1, \dots, \eps_r ; \bars{h}, \bars{\mu}, \bars{\nu}, \sigma_{\cc})}

\newcommand\ula{u_{L, \goth{a}}}
\newcommand\ulah{u_{L, \goth{a}, \bars{h}}}
    
\newcommand\homlz{\mathrm{Hom}_{\zz}(L, \zz)}
\newcommand\homlc{\mathrm{Hom}_{\zz}(L, \cc)}
\newcommand\homlambdaz{\mathrm{Hom}_{\zz}(\Lambda, \zz)}
\newcommand\homlambdac{\mathrm{Hom}_{\zz}(\Lambda, \cc)}
\newcommand\zexc{z}
  
\newcommand\zfone{\mathcal{Z}_{\goth{f}}^1}

\begin{center}
\end{center}
\begin{center}
\large \noindent Elliptic units above fields with exactly one complex place
\end{center}
\medskip
\begin{center}
Pierre L. L. Morain\footnotemark[1]\footnotetext[1]{Sorbonne Université and Université Paris Cité, CNRS, IMJ-PRG, F-75005 Paris, France. This PhD work is funded by the École polytechnique,  Palaiseau, France.}
\end{center}
\medskip
\begin{center}
\noindent \textbf{Abstract:}
\end{center}
In this work we explore the construction of abelian extensions of number fields with exactly one complex place using multivariate analytic functions in the spirit of Hilbert's 12th problem. To this end we study the special values of the multiple elliptic Gamma functions introduced in the early 2000s by Nishizawa following the work of Felder and Varchenko on Ruijsenaars' elliptic Gamma function. We construct geometric variants of these functions enjoying transformation properties under an action of $\mathrm{SL}_{d}(\mathbb{Z})$ for $d \geq 2$. The evaluation of these functions at points of a degree $d$ field $\mathbb{K}$ with exactly one complex place following the scheme of a recent article by Bergeron, Charollois and Garc\'ia \cite{BCG} seems to produce algebraic numbers. More precisely, we conjecture that such infinite products yield algebraic units in abelian extensions of $\mathbb{K}$ related to conjectural Stark units and we provide numerical evidence to support this conjecture for cubic, quartic and quintic fields.

\normalsize \bigskip \bigskip

\noindent \textbf{Acknowledgments:} 
This work is part of an on-going PhD work and the author would like to thank his advisors Pierre Charollois and Antonin Guilloux for their guidance and their helpful comments on this work.
\newpage

\tableofcontents

\section{Introduction}

In explicit class field theory the main task is to provide a way to compute the abelian extensions of a number field $\kk$. The celebrated theorem of Kronecker and Weber shows that the exponential function plays a key role in the description of the abelian extensions of $\qq$. Since then, the theory of Complex Multiplication has provided explicit constructions for abelian extensions of imaginary quadratic fields. There are many candidates for the function which plays the role of the exponential function here, including Klein's $j$-invariant, Weber's functions or Jacobi's $\theta$ functions, which are all related to modular forms for $\slnz{2}$. Elliptic units are obtained by evaluating ratios of $\theta$ functions at points of an imaginary quadratic field. For instance, define 
$$\theta(z, \tau) = \prod_{m \geq 0} (1 - \exp(2i\pi((m+1)\tau-z)))(1-\exp(2i\pi(m\tau + z)))$$
for $z \in \cc$ and $\tau \in \hh = \{w \in \cc, \Im(w) > 0\}$. Then we may compute 
$$\theta\left(\frac{7}{91}, \frac{10+e^{2i\pi/3}}{91}\right)^7\theta\left(\frac{7}{13}, \frac{10+e^{2i\pi/3}}{13}\right)^{-1} \approx 0.0196475... - i\cdot0.6399892...$$
which is a unit outside an ideal above $13$ in the field $\qq[x]/(x^4+3x^3+32x^2+13)$, an abelian extension of $\qq(e^{2i\pi/3})$ ramified only above a prime ideal of norm $13$. For a more complete presentation on the subject, we refer to \cite{Robert}.

Finding general analytic functions which play the role of the exponential function for other number fields constitutes Hilbert's 12th problem. The situation becomes more complicated outside of the two cases above due to the presence of an infinite set of units in $\kk$. Other efforts have been made to compute elements in class fields using analytic functions, especially in the case of real quadratic fields. In \cite{Stark}, Stark analyses regulators in abelian extensions and studies how they can be factored. This led him to refine Dirichlet's class number formula and to conjecture the existence of special units in class fields which should yield the exponential of the values of $L$-functions at $s=0$ as a combination of their valuations at infinite places. To relate this to Hilbert's 12th problem, we are interested in finding a way to compute these units directly, instead of computing absolute values. It was later shown that there are two simpler cases in the Stark conjectures (see for instance \cite{rankoneStark}). The first case deals with totally real number fields and it was solved using $p$-adic methods by Dasgupta and Kakde in \cite{Dasgupta}. The second case deals with number fields with exactly one complex place, the so-called ``almost totally real fields'' (ATR for short).

In this article we are concerned with the case of ATR number fields of degree $d \geq 3$ for which Hilbert's 12th problem and Stark's conjectures are still open (the case $d = 2$ of imaginary quadratic fields is dealt with using CM theory). Previous work on the abelian Stark conjecture for complex cubic fields was carried out by Ren and Sczech in \cite{RenSczech} where a conjectural explicit formula is given for the Stark unit. Recently, Bergeron, Charollois and Garc\'ia \cite{BCG} used Ruijsenaars' elliptic Gamma function coming from mathematical physics (see \cite{Ruijsenaars}), to build conjectural elliptic units above complex cubic fields and related their construction to the conjectural Stark units. For instance, if $z = e^{2i\pi/3}2^{1/3}$ is the root of the polynomial $x^3-2$ in the upper half-plane $\hh$ we may compute with high precision (more than 1000 digits) the quotient
$$ \frac{\Gamma\left(\frac{-1}{3}, \frac{z^2+z+3}{15}, \frac{z-8}{15}\right)^{-5}}{\Gamma\left(\frac{-5}{3}, \frac{z^2+z+3}{3}, \frac{z-8}{3}\right)^{-1}} \approx -1.2937005... + i\cdot1.4743341...$$
(see section \ref{defgamma} for the definition of the elliptic Gamma function). It is close to a root of the polynomial $x^6 + 3x^5 + 6x^4 + 5x^3 + 6x^2 + 3x + 1$ which is the minimal polynomial over $\qq$ of a unit in an abelian extension of $\qq(z)$. We treat this example in section \ref{cubicsimple} below. The elliptic Gamma function had been studied by Felder and Varchenko \cite{FV} and shown to satisfy modular identities for the group $\slnz{3}$ involving $\theta$ functions. The construction of conjectural elliptic units in \cite{BCG} is inspired by the construction of elliptic units above imaginary quadratic fields using $\theta$ functions, and the constructions carried out in this article will follow the same principles for abritrary ATR number fields. Previous work was also carried out in the case of even degree $2d'$ ATR number fields containing a maximal real subfield $\mathbb{F}$ of degree $d'$ in \cite{CharolloisDarmon} using Hilbert modular forms for the modular group $\sln{2}{\mathcal{O}_{\mathbb{F}}}$. It is interesting to note that the existence of such a real subfield will often be an obstacle to the approach developed in this article. 

This work proposes a generalisation beyond complex cubic fields of the construction carried out in \cite{BCG} using multiple elliptic Gamma functions (denoted $G_r$ for $r \geq 0$) which were introduced by Nishizawa in \cite{Nishizawa} and studied by Narukawa in \cite{Narukawa} (see section \ref{defgr}). In order to achieve this, we first construct geometric variants of the $G_r$ functions of Nishizawa in section \ref{sectiongammaell}, upgrading the construction of the geometric elliptic Gamma function in \cite{FDuke} and we show that these functions enjoy nice properties under the action of $\slnz{r+2}$. In section \ref{sectionarithmetic}, we show how to evaluate these functions at points in an ATR number field $\kk$ of degree $d=r+2$ to obtain conjectural algebraic units in abelian extensions of $\kk$. This is expressed in Conjecture \ref{conjecture} which generalises the conjecture in \cite{BCG} to higher degree ATR number fields. In section \ref{sectioncomputing}, under some assumptions, we show how to choose the evaluation points for the geometric $G_r$ functions to efficiently obtain conjectural algebraic units. We also describe the algorithmic approach to obtaining numerical evidence to support Conjecture \ref{conjecture}, and in section \ref{sectionnumerical}, we provide many examples for ATR cubic, quartic and quintic fields.

\section{Multiple elliptic Gamma functions and their geometric variants}\label{sectiongammaell}
In this section, we recall the properties of the elliptic Gamma function and its geometric variants. Then we will present the multiple elliptic Gamma functions of Nishizawa and construct geometric variants of these functions, upgrading the construction of the geometric elliptic Gamma functions in \cite{FDuke}. 

\subsection{The geometric elliptic Gamma functions}

\subsubsection{The elliptic Gamma function}\label{defgamma}

The elliptic Gamma function was first introduced by Ruijsenaars \cite{Ruijsenaars} and was studied later at length by Felder and Varchenko \cite{FV} who made the connexion with the above $\theta$ function. For $z \in \cc$ and $\tau, \sigma \in \hh$ put
$$\Gamma(z,\tau,\sigma) = \prod_{m,n \geq 0}\left(\frac{1-\exp(2i\pi((m+1)\tau +(n+1)\sigma-z))}{1-\exp(2i\pi(m\tau+n\sigma+z))} \right)$$
This function converges absolutely for $z \in \cc, z \not\in \zz +\zz_{\leq 0}\tau + \zz_{\leq 0}\sigma$ and it is meromorphic in this domain. It is worth noting that this function has a nice expression as the exponential of an infinite sum involving sinuses when $|\Im(2z - \tau - \sigma)| < |\Im(\tau)| + |\Im(\sigma)|$:
\begin{equation}
\Gamma(z, \tau, \sigma) = \exp\left(\sum_{j \geq 1} \frac{1}{2ij} \frac{\sin(\pi j(2z - \tau - \sigma))}{\sin(\pi j \tau) \sin(\pi j \sigma)}\right)
\end{equation}
The range of parameters $\tau, \sigma$ can thus be extended to $\tau, \sigma \in \cc - \rr$ by putting:
$$\Gamma(z, -\tau, \sigma) = \frac{1}{\Gamma(z + \tau, \tau, \sigma)},~~~~ \Gamma(z, \tau, -\sigma) = \frac{1}{\Gamma(z + \sigma, \tau, \sigma)}$$
In [\hspace{1sp}\cite{FV}, Theorems 3.1 and 4.1], Felder and Varchenko proved the following properties for the elliptic Gamma function related to the $\theta$ function of Complex Multiplication. If $\tau, \sigma \in \cc-\rr$ then
\begin{align}
\Gamma(z + 1, \tau, \sigma) = \Gamma(z, \tau + 1, \sigma) & = \Gamma(z, \tau, \sigma+1) = \Gamma(z, \tau, \sigma) \nonumber \\
\Gamma(z + \tau, \tau, \sigma) & = \theta(z, \sigma)\Gamma(z, \tau, \sigma) \nonumber \\
\Gamma(z + \tau + \sigma, \tau, \sigma) & = \Gamma(-z, \tau, \sigma)^{-1} \nonumber
\end{align}
Finally, if $\sigma/\tau \not\in \rr$ then
\begin{equation} \label{modgamma}
\Gamma(z, \tau, \sigma)^{-1}\Gamma\left(\frac{z}{\tau}, \frac{-1}{\tau}, \frac{\sigma}{\tau}\right) \Gamma\left(\frac{z-\tau}{\sigma}, -\frac{\tau}{\sigma}, -\frac{1}{\sigma}\right)^{-1} = \exp(i\pi P_1(z,\tau, \sigma))
\end{equation}
where $$P_1(z, \tau, \sigma) = \frac{z^3}{3\tau\sigma} - \frac{\tau + \sigma - 1}{2\tau\sigma}z^2 + \frac{\tau^2 + \sigma^2 + 3\tau\sigma - 3\tau - 3\sigma + 1}{6\tau\sigma}z $$
$$+ \frac{1}{12}(\tau + \sigma - 1)\left(\frac{1}{\tau} + \frac{1}{\sigma} - 1\right)$$ 

These properties are transformation properties for the elliptic Gamma function under the action of $\slnz{3}$. This will be more apparent for the geometric elliptic Gamma functions which we introduce now.

\subsubsection{Geometric variants of the elliptic Gamma function}\label{geomgamma}

In \cite{FDuke}, Felder, Henriques, Rossi and Zhu have built elliptic Gamma functions associated to pairs of primitive vectors in a free $\zz$-module of rank $3$. We briefly review the construction which we will adapt to Nishizawa's $G_r$ functions. Consider two linearly independent primitive vectors $a,b$ in $\Lambda$ a free $\zz$-module of rank $3$. Fix a basis $B$ (and therefore an orientation) of $\Lambda$ and write $\det_{B}(a,b,\cdot) = s\gamma$ for a unique integer $s>0$ and a unique linear form $\gamma \in L = \homlambdaz$. The basis $B$ identifies $\Lambda \iso \zz^3$ which fixes an action of $\slnz{3}$ on $\Lambda$. Now write
$$C(a,b) = \{\delta \in L, \delta(a) > 0, \delta(b) \leq 0\}$$
The geometric elliptic $\Gamma_{a,b}$ function associated to the pair $(a,b)$ is defined as follows:
$$\Gamma_{a,b}(w, x, \Lambda) = \frac{\prod_{\delta \in C(a,b)/\zz\gamma}(1-\exp(-2i\pi(\delta(x) - w)/\gamma(x)))}{\prod_{\delta \in C(b,a)/\zz\gamma}(1-\exp(2i\pi(\delta(x) - w)/\gamma(x)))}$$
where $w \in \cc$ and $x \in \homlc$ and where we understand $\delta(x) = x(\delta)$ by biduality. The study of the convergence properties of this function was carried out in \cite{FDuke} and is better understood with an alternative definition of this function. Consider $\alpha, \beta \in L$ such that $\beta(a) = \alpha(b) = 0, \alpha(a) > 0$, and $  \beta(b) > 0$. We call such a pair a positive dual family to $(a,b)$. We define 
$$ F(\alpha, \beta) = \{\delta \in L, 0 \leq \delta(a) < \alpha(a), 0 \leq \delta(b) < \beta(b)\}$$
Then in \cite{FDuke} the authors prove that the following alternative definition holds:
\begin{equation}\label{productformula}\Gamma_{a,b}(w, x, \Lambda) =  \prod_{\delta \in F(\alpha, \beta)/\zz\gamma} \Gamma\left(\frac{w+\delta(x)}{\gamma(x)}, \frac{\alpha(x)}{\gamma(x)}, \frac{\beta(x)}{\gamma(x)}\right)
\end{equation}
The right-hand side is a finite product of elliptic Gamma functions and is independent of the choice for $\alpha, \beta$. As a function of $x$ it converges outside the hyperplane $\gamma(x) = 0$ whenever $\alpha(x)/\gamma(x)$ and $\beta(x)/\gamma(x)$ are not real. As a function of $w$, it converges outside the discrete set of poles of every elliptic Gamma function involved. Note that the conditions $\alpha(x)/\gamma(x), \beta(x)/\gamma(x) \in \rr$ describe a $\rr$-hyperplane of the $\rr$-vector space $\homlc$, so that the geometric elliptic Gamma function is well-defined on a dense open set of $\cc \times \homlc$ endowed with the finite dimensional $\cc$-vector space topology.

\textbf{Remark : } Formula \refp{productformula} could have been taken as the definition of the geometric elliptic Gamma function associated to the pair $(a,b)$. This will be done below (see Proposition \ref{propdefgeom}) when considering $G_r$ functions associated to $r+1$ vectors in a free $\zz$-module of rank $r+2$. It will remain to check that this definition is independent of the choice for $\alpha, \beta$, and this can be done by direct computation in reverting the proof of formula \refp{productformula} carried out in \cite{FDuke}.

On the computational side, the geometric function $\Gamma_{a,b}$ is here computed as a product of $\#|F(\alpha, \beta)/\zz\gamma|$ ordinary elliptic Gamma functions. A standard choice gives $\alpha(a) = \beta(b) = s = \#|F/\zz\gamma|$. This choice is given by considering a third vector $c$ such that $\det_B(a,b,c) = s$ and writing the transpose of the comatrix of $(a,b,c)$ as $(\alpha, \beta, s\gamma)^{T}$. The properties of the elliptic Gamma function translate into nice properties for the geometric $\Gamma_{a,b}$ functions under the action of $\slnz{3}$ as shown in \cite{FDuke}. Namely, in the domains where every term converges (which are dense open sets), we have the following inversion relation:
$$\Gamma_{a,b}(w, x, \Lambda)\Gamma_{b,a}(w,x, \Lambda) = 1$$
We also have for $a,b,c$ linearly independent:
\begin{equation}\label{cocyclegamma}\Gamma_{a,b}(w, x, \Lambda)\Gamma_{b,c}(w,x, \Lambda)\Gamma_{c,a}(w,x, \Lambda) = \exp\left(i\pi P_{a,b,c}(w,x, \Lambda)\right)
\end{equation}
where $P_{a,b,c}$ is a polynomial in $w$ of degree 3. Furthermore, we have the equivariance relations: for all $g \in \slnz{3}$,
$$\Gamma_{ga, gb}(w, gx, \Lambda) = \Gamma_{a, b}(w,x, \Lambda)$$
$$P_{ga, gb, gc}(w,gx, \Lambda) = P_{a,b,c}(w,x, \Lambda)$$
The polynomial $P_{a,b,c}$ will be presented below alongside other Bernoulli polynomials serving the same purpose for higher degree (see definition \ref{defgeombern} in section \ref{grgeom}). It is worth noting that the 3-term relation \refp{cocyclegamma} can be understood as a $1$-cocycle relation for $\slnz{3}$. Consider the cocycle $\psi_a(g) = \Gamma_{a, ga}$ for a primitive vector $a \in \Lambda$ and $g \in \slnz{3}$. Then we may rewrite formula $\refp{cocyclegamma}$ for fixed $w, x$ as
$$\Gamma_{a, g_1a} \Gamma_{g_1a, g_1g_2a} \Gamma_{g_1g_2a, a} = \exp(i\pi P_{a, g_1a, g_1g_2a})$$
which yields the coboundary relation:
\begin{equation}\label{cocycleone}
d\psi_a(g_1, g_2) := \psi_a(g_1) (g_1\cdot\psi_a(g_2)) \psi_a(g_1g_2)^{-1} = \exp(i\pi P_{a, g_1a, g_1g_2a})
\end{equation}

\subsection{Multiple elliptic Gamma functions}\label{defgr}

Here we recall the definition and properties of Nishizawa's $G_r$ functions which generalise the $\theta$ and elliptic Gamma functions above, and we construct their geometric variants $G_{r,\abar}$ upgrading the construction in \cite{FDuke} to higher degrees.

\subsubsection{Nishizawa's $G_r$ functions}

In \cite{Nishizawa}, Nishizawa introduced the hierarchy of $G_r$ functions that generalise both the $\theta$ and the elliptic Gamma functions. Later, Narukawa proved the modularity property for the $G_r$ functions by relating them to multiple sine functions.
To simplify the notations, we will write as in \cite{Nishizawa}, with a clear separation for base and exponential variables:
\begin{align*}
\taubar = (\tau_j)_{0 \leq j \leq r} \in \cc^{r+1}~&;~\qbar = (q_j)_{0 \leq j \leq r} = (\exp(2i\pi \tau_j))_{0 \leq j \leq r} \\
|\taubar| = \sum_{j = 0}^r \tau_j~&;~|\qbar| = \prod_{j = 0}^r q_j \\
\taubar^{-}(j) = (\tau_n)_{0 \leq n \neq j \leq r}~&;~\qbar^{-}(j) = (q_n)_{0 \leq n \neq j \leq r} \\
\taubar[j] = (\tau_0, \dots, \tau_{j-1}, -\tau_j, \tau_{j+1}, \dots, \tau_r)~&;~\qbar[j] = (q_0, \dots, q_{j-1}, q_j^{-1}, q_{j+1}, \dots, q_r) 
\end{align*}
Moreover, for a $(r+1)$-uple $\mbar = (m_0, \dots, m_r) \in \zz^{r+1}$, we write:
$$ \mbar + \bars{1} = (m_j + 1)_{0 \leq j \leq r},~\mbar\cdot\taubar = (m_j\tau_j)_{0 \leq j \leq r},~\qbar^{\mbar} = (q_j^{m_j})_{0 \leq j \leq r}$$
$$ \prod_{\mbar \,\geq 0} = \prod_{m_0 \geq 0} \dots \prod_{m_r \geq 0}$$

\noindent Let $r \geq 0$ and $\taubar = (\tau_0, \dots, \tau_r) \in \hh^{r+1}$. For $z \in \cc$, ($z \not\in \zz + \sum_{j=0}^r \zz_{\leq 0}\tau_j$ if $r$ is odd), put $x = \exp(2i\pi z)$ and define as in \cite{Nishizawa}:
$$ G_r(z, \taubar) = \prod_{\mbar \,\geq 0} \left(1 - \left|\qbar^{\mbar +\onebar}\right|x^{-1}\right)\left(1-\left|\qbar^{\mbar}\right|x\right)^{(-1)^r}$$

\noindent Note that $G_0$ and $G_1$ are respectively the $\theta$ function and the elliptic Gamma function. It is once again worth noting that the $G_r$ functions have nice expressions as the exponential of infinite sums involving sinuses (see [\hspace{1sp}\cite{Nishizawa}, Proposition 3.6]), namely:
\begin{equation}\label{exponentialformula}
G_r(z, \taubar) = \begin{cases} \exp\left(\sum_{j \geq 1} \frac{1}{(2i)^rj} \frac{\sin(\pi j(2z - |\taubar|))}{\prod_{k = 0}^{r} \sin(\pi j \tau_k)}\right) \mathrm{~if~}r~\mathrm{is~odd} \\
\exp\left(\sum_{j \geq 1} \frac{2}{(2i)^{r+1}j} \frac{\cos(\pi j(2z - |\taubar|))}{\prod_{k = 0}^{r} \sin(\pi j \tau_k)}\right) \mathrm{~if~}r~\mathrm{is~even}
\end{cases}
\end{equation}
This expression is valid provided that $\taubar \in (\cc - \rr)^{r+1}$ and $|\Im(2z - |\taubar|)| < \sum_{j=0}^r |\Im(\tau_j)|$ and allows us to extend the range of parameters to $\taubar \in (\cc - \rr)^{r+1}$ by putting:
\begin{equation}\label{inversionGr}
G_r(z, \taubar[j]) = G_r(z + \tau_j, \taubar)^{-1}
\end{equation}
The $G_r$ functions satisfy modular properties similar to that of the elliptic Gamma function, as proven by Nishizawa in \cite{Nishizawa} and later by Narukawa in \cite{Narukawa}. First, the $G_r$ functions are $1$-periodic in each of their arguments and they satisfy:
\begin{align}
G_r(z + |\taubar|, \taubar) & = G_r(-z, \taubar)^{(-1)^r} \nonumber \\
G_r(-z, -\taubar) & = G_r(z, \taubar)^{-1} \label{inversiontotale}
\end{align}
Furthermore, if $r \geq 1$, the $G_r$ functions are almost periodic in $z$ with periods $\tau_j$ for $0 \leq j \leq r$ with a correction factor involving a lower degree function:
$$G_r(z + \tau_j, \taubar) = G_{r-1}(z, \taubarsj{j})G_{r}(z, \taubar)$$
A consequence of these properties which we will use is the following formula for complex conjugation:
\begin{equation}\label{grconjugate}
G_r(\overline{z}, \overline{\tau_0}, \dots, \overline{\tau_r}) = \overline{G_r(z, \tau_0, \dots, \tau_r)}^{(-1)^r}
\end{equation}

\subsubsection{Bernoulli polynomials and Narukawa's theorem}\label{Bernoulli}

The last missing property is crucial for arithmetic purposes and was later proved by Narukawa. To state Narukawa's theorem, we first need to introduce the multiple Bernoulli polynomials. Consider an integer $d \geq 1$. Let $\omegabar = (\omega_1, \dots, \omega_d) \in (\cc-\{0\})^{d}$ with $\omega_j \neq 0$ for all $1 \leq j \leq d$. We define the multiple Bernoulli polynomials $B_{d,n}^{*}(z, \omegabar)$ with the following generating function:
$$e^{zt}\prod_{j = 1}^d\frac{\omega_j t}{e^{\omega_jt}-1} = \sum_{n \geq 0} B_{d,n}^{*}(z, \omegabar)\frac{t^n}{n!} = \left(\sum_{m \geq 0} \frac{z^mt^m}{m!}\right)\prod_{1 \leq j \leq d} \left(\sum_{k_j \geq 0} B_{k_j} \omega_j^{k_j}\frac{t^{k_j}}{k_j!}\right)$$
where the $B_{k_j}$ are the usual Bernoulli numbers. The classic Bernoulli polynomials are obtained with the choice $d = 1$, $\omega_1 = 1$, and the Bernoulli numbers are obtained for $z = 0$. These polynomials obey many relations which we sum up below and which can easily be obtained from the properties of the generating function.
For $\omegabar \in (\cc-\{0\})^{d},~ B_{d,n}^{*}(z, \omegabar)$ is a degree $n$ homogeneous polynomial in $d+1$ variables, which is symmetric in the $d$ variables of $\omegabar$. Moreover, the coefficients of $B_{d,n}^{*}(z, \omegabar)$ are rational numbers. In \cite{Narukawa}, Narukawa used the rescaled version $B_{d,n}(z, \omegabar) = (\prod_{j = 1}^d \omega_j^{-1}) B_{d,n}^{*}(z, \omegabar)$ and proved the following properties:
\begin{align}
B_{d,n}(\alpha z, \alpha \omegabar) & = \alpha^{n-d} B_{d,n}(z, \omegabar),~\forall \alpha \in \cc - \{0\} \nonumber \\
B_{d,n}(z + |\omegabar|, \omegabar) & = (-1)^n B_{d,n}(-z, \omegabar) \nonumber \\
B_{d,n}(z + \omega_j, \omegabar) & = B_{d,n}(z, \omegabar) + nB_{d-1,n-1}(z, \omegabar^{-}(j)) \nonumber \\
B_{d,n}(z + \omega_j, \omegabar) & = -B_{d,n}(z, \omegabar[j]) \nonumber \\
\partialdx{z}B_{d,n}(z, \omegabar) & = nB_{d,n-1}(z, \omegabar) \nonumber
\end{align}
Narukawa's modularity theorem (see [\hspace{1sp}\cite{Narukawa}, Theorem 7]) involves the polynomials $B_{d,d}^{*}$ or rather the rescaled versions $B_{d,d}$. For instance, we have
$$(\omega_1\omega_2)B_{2,2}(z, \omega_1, \omega_2) = B_{2,2}^{*}(z, \omega_1, \omega_2) = z^2 - z(\omega_1 + \omega_2) + \frac{\omega_1^2+\omega_2^2+3\omega_1\omega_2}{6}$$
Then, Narukawa's theorem can be stated as follows. Let $r \geq 0$. Take $\omegabar \in (\cc - \{0\})^{r+2}$ and suppose that $\omega_j/\omega_n \in \cc - \rr$ for all $j \neq n$. Then for $z \in \cc$ outside the discrete set of poles of the left-hand side
\begin{equation} \label{modGr}
\prod_{j = 0}^{r+1} G_r\left(\frac{z}{\omega_j}, \left(\frac{\omega_n}{\omega_j}\right)_{n \neq j}\right) = \exp\left(\frac{-2i\pi}{(r+2)!}B_{r+2,r+2}(z, \omegabar)\right)
\end{equation}

\subsubsection{Geometric $G_{r, \abar}$ functions}\label{grgeom}

In this section we upgrade the construction of section \ref{geomgamma} and define geometric variants of the $G_r$ functions. We also prove (Theorem \ref{theoremslnz}) the modular property of the geometric variants and their equivariance property under the action of $\slnz{r+2}$. In the literature one may find other variants of these functions, such as the $G_r$ functions associated to cones by Winding in \cite{winding}. Although our construction inspired by \cite{FDuke} differs from Winding's construction, it would be interesting to understand the links between the two variants. As mentioned, we will use the alternative definition of \ref{geomgamma} to start and prove that this definition is indeed valid. Throughout this section, we consider a free $\zz$-module $\Lambda$ of rank $r+2$ with an orientation form given by the determinant in a fixed basis. Let us start by defining positive dual families which will be used in the rest of this work.

\begin{general}{Definition}
Let $(a_0, \dots, a_m)$ be a family of $m+1$ independent primitive vectors in a lattice $\Lambda$. We call $(\alpha_0, \dots, \alpha_m) \in L = \homlambdaz$ a positive dual family to $\abar = (a_0, \dots, a_m)$ if for all $0 \leq j \leq m$ the following holds:
$$ \alpha_j(a_j) > 0, ~~~~\alpha_j(a_k) = 0, \, \forall\, k \neq j$$
\end{general}

\noindent The important cases in this work will be those where $m = r$ and $m = r+1$ in the lattice $\Lambda$ of rank $r+2$. The following lemmas show that in these cases two positive dual families to the same family $\abar$ are closely related.

\begin{general}{Lemma}\label{lemmachoice}
Let $\abar = (a_0, \dots, a_r)$ be a family of $r+1$ linearly independent primitive vectors in the oriented lattice $\Lambda$ of rank $r+2$. Take two positive dual families $(\alpha_0, \dots, \alpha_r), (\alpha'_0, \dots, \alpha'_r)$ to $\abar$ in $L = \homlambdaz$. For $0 \leq j \leq r$ write $s_j = \alpha_j(a_j)$ and $s'_j = \alpha_j'(a_j)$. Then, there are rational numbers $r_j$ such that $\alpha_j = (s_j/s'_j)\alpha'_j + r_j \det(a_0, \dots, a_r, \cdot)$ for all $0 \leq j \leq r$. 
\end{general}

\begin{proof}
Consider $\gamma_j = s'_j\alpha_j-s_j\alpha'_j$. Then for all $0 \leq k \leq r$, we get $\gamma_j(a_k) = 0$. This means that either $\gamma_j = 0$ and then $r_j = 0$ or the $\qq$-linear forms $\gamma_j$ and $\det(a_0, \dots, a_r, \cdot)$ defined on the $\qq$-vector space $\Lambda \otimes \qq$ share the same kernel, and therefore, they must be linearly dependent. 
\end{proof}

\begin{general}{Lemma}\label{lemmaunique}
A family $\abar = (a_0, \dots, a_{r+1})$ of $r+2$ linearly independent primitive vectors in a lattice $\Lambda$ of rank $r+2$ has exactly one positive dual family $\alphabar$ containing only primitive vectors. We call this family the primitive positive dual family to $\abar$.
\end{general}

\begin{proof}
Consider two such families $\alphabar, \,\alphabar'$ and write again $s_j = \alpha_j(a_j), s'_j = \alpha'_j(a_j)$. Consider once again $\gamma_j =  s'_j\alpha_j-s_j\alpha'_j$ such that for all $0 \leq k \leq r+1$, $\gamma_j(a_k) = 0$. The family $(a_0, \dots, a_{r+1})$ is a basis of the $\qq$-vector space $\Lambda \otimes \qq$, so $\gamma_j = 0$ for all $0 \leq j \leq r+1$. This means that $\alpha_j = s_j\alpha'_j/s'_j$ and we write $s_jn_j + s'_jn'_j = d_j$ where $d_j = gcd(s_j, s'_j)$. Then $n_j\alpha_j +n'_j\alpha'_j = d_j\alpha'_j/s'_j \in L$. Because $\alpha'_j$ is a primitive vector in $L$ by assumption, we must have $d_j = s'_j$, and then because $\alpha_j = s_j\alpha'_j/d_j$ is primitive in $L$ we must also have $d_j = s_j$. The conditions $\alpha_j(a_j) = s_j$, $\alpha_j(a_k) = 0$ for $k \neq j$ describe exactly one linear form on the $\qq$-vector space $\Lambda \otimes \qq$. To compute $\alphabar$ we compute the comatrix of $\abar$ in a suitable basis of $\Lambda$ and we rescale the linear forms to obtain primitive elements in $L$.
\end{proof}

\noindent Now we may define the geometric variants of the $G_r$ functions using lemma \refp{lemmachoice} and adapting formula \ref{productformula} to higher rank lattices:

\begin{general}{Proposition}\label{propdefgeom}
Let $\abar = (a_0, \dots, a_r)$ be a family of $r+1$ linearly independent primitive vectors in the oriented lattice $\Lambda$ of rank $r+2$. There is a unique $s \in \zz_{>0}$ and a unique primitive element $\gamma \in L = \homlambdaz$ such that $\forall c \in \Lambda,~ \det(\abar, c) = s\gamma(c)$. For any choice of positive dual family $(\alpha_0, \dots, \alpha_r)$ to $\abar$ in $L$ the function 
\begin{equation}\label{defgeom}
G_{r, \abar}^{\alphabar}(w, x, \Lambda) = \prod_{\delta \in F(\alphabar)/\zz\gamma} G_r\left(\frac{w+\delta(x)}{\gamma(x)}, \frac{1}{\gamma(x)}\alphabar(x) \right)
\end{equation}
$$F(\alphabar) = \{ \delta \in L, \forall\,0 \leq j \leq r,~0 \leq \delta(a_j) < \alpha_j(a_j) \}$$
is well defined for $(w, x)$ in a dense open set of the $\cc$-vector space $\cc \times \homlc \iso \cc \times \cc^{r+2}$ endowed with the finite dimensional $\cc$-vector space topology. Furthermore, it is independent of the choice for $\alphabar$.
\end{general}

\begin{proof} 
The set $\homlc$ is a $\rr$-vector space of dimension $2r+4$. The condition $\gamma(x) = 0$ describes a subspace of dimension $2r+2$ and the condition $\alpha_j(x)/\gamma(x) \in \rr$ describes a subspace of dimension $2r+3$. The function $G_{r,\abar}^{\alphabar}(w, x, \Lambda)$ is therefore well defined for $x$ outside a finite union of $\rr$-hyperplanes in $\homlc$, i.e. in a dense open set. The right-hand side in \refp{defgeom} is a finite product of meromorphic functions in $w$ and it has a discrete set of poles. 

Let us now prove that this definition is indeed independent of the choice for the positive dual family $\alphabar$. Take $x \in \homlc$ such that $\gamma(x) \neq 0$ and $\alpha_j(x)/\gamma(x) \not\in\rr$ for all $0 \leq j \leq r$. Put for $0 \leq j \leq r$, $d_j = \pm 1$ such that $d_j\alpha_j(x)/\gamma(x) \in \hh$. Put also $D = \sum_{j = 0}^{r} (d_j-1)/2$. Then using the inversion relation (\ref{inversionGr}) we get:
\begin{multline*}
G_{r, \abar}^{\alphabar}(w, x, \Lambda)^{(-1)^D}  = \\
\prod_{\delta \in F(\alphabar)/\zz\gamma} G_r\left(\frac{w+\delta(x)}{\gamma(x)} + \sum_{j = 0}^r \frac{d_j-1}{2}\frac{\alpha_j(x)}{\gamma(x)}, \left(\frac{1}{\gamma(x)}d_j\alpha_j(x)\right)_{0 \leq j\leq r}\right)
\end{multline*}
\begin{multline*}
G_{r, \abar}^{\alphabar}(w, x, \Lambda)^{(-1)^D} = \prod_{\delta \in F(\alphabar)/\zz\gamma} \prod_{\mbar\,\geq 0} \Big[\left(1-e^{2i\pi\left(\sum_{j = 0}^r \frac{(d_jm_j+(1+d_j)/2)\alpha_j(x)}{\gamma(x)} - \frac{w + \delta(x)}{\gamma(x)}\right)}\right) \\
\times \left(1-e^{2i\pi\left(\sum_{j = 0}^r \frac{(d_jm_j+(1-d_j)/2)\alpha_j(x)}{\gamma(x)} + \frac{w + \delta(x)}{\gamma(x)}\right)}\right)^{(-1)^r}\Big]
\end{multline*}
Now write $C^{\pm}(\alphabar, x)$ for the set of $\delta' \in L$ satisfying for all $0 \leq j \leq r$:
$$\begin{cases}
\delta'(a_j) > 0, \mathrm{~if~} d_j = \pm 1 \\
\delta'(a_j) \leq 0, \mathrm{~if~} d_j = \mp 1
\end{cases}$$
Consider $\delta' \in C^{+}(\alphabar, x)$. This means that if $d_j = 1$ then $\delta'(a_j) > 0$ and we can perform Euclidian division by $\alpha_j(a_j)$ so that there exists a unique integer $m_j > 0$ such that $0 \leq (m_j\alpha_j-\delta')(a_j) < \alpha_j(a_j)$. On the contrary, if $d_j = -1$ then $\delta'(a_j) \leq 0$ and there exists a unique integer $m_j \geq 0$ such that $0 \leq (-m_j\alpha_j-\delta')(a_j) < \alpha_j(a_j)$. Then for all $0 \leq k \leq r$: 
$$0 \leq \left(\sum_{j = 0}^r d_jm_j\alpha_j - \delta'\right)(a_k) < \alpha_k(a_k)$$
This shows that $C^{+}(\alphabar, x)$ can be written as a disjoint union
$$C^{+}(\alphabar, x) = \bigcup_{\delta \in F(\alphabar)/\zz\gamma} \bigcup_{\mbar\, \geq 0} \left\{-\delta + \sum_{j = 0}^r (d_jm_j + (1+d_j)/2)\alpha_j\right\} + \zz \gamma$$
and a similar statement holds for $C^{-}(\alphabar, x)$ so that $G_{r, \abar}^{\alphabar}(w, x, \Lambda)^{(-1)^D}$ is equal to
$$ \prod_{\delta' \in C^{+}(\alphabar, x)/\zz\gamma}\left(1-e^{2i\pi\left(\frac{\delta'(x)-w}{\gamma(x)}\right)}\right)\prod_{\delta' \in C^{-}(\alphabar,x )/\zz\gamma}\left(1-e^{-2i\pi\left(\frac{\delta'(x)-w}{\gamma(x)}\right)}\right)^{(-1)^r}$$
Then, we only need to show that the sets $C^{\pm}(\alphabar, x)$ are independent of the choice for $\alphabar$. Consider another positive dual family $\alphabar'$ to $\abar$. Write $\alpha_j(a_j) = s_j > 0$ and $\alpha_j'(a_j) = s'_j>0$. Then from lemma \ref{lemmachoice} there is a rational number $r_j$ such that
$$\alpha'_j = \frac{s'_j}{s_j}\alpha_j + r_j\gamma$$
For any $x \in \homlc$ we get:
$$d_j\frac{\alpha'_j(x)}{\gamma(x)} = d_jr_j + d_j\frac{s'_j}{s_j}\frac{\alpha_j(x)}{\gamma(x)} \in \hh$$

This shows that the signs $d_j$ (and therefore also $D$) are independent of the choice for $\alphabar$. By construction we get $C^{\pm}(\alphabar, x) = C^{\pm}(\alphabar', x)$ and the definition of the geometric $G_{r, \abar}^{\alphabar}$ function is independent of the choice for $\alphabar$.
\end{proof}

From now on we denote by $G_{r, \abar} := G_{r, \abar}^{\alpha}$ the geometric $G_r$ function associated to $\abar$ for any suitable choice of $\alphabar$. A choice that is interesting for us is defined by the comatrix or rather by the rescaled comatrix. Let $c$ be a vector such that $\gamma(c) = 1$. Then the elementary divisors of the concatenation $(\abar, c)$ can be written as $[A_r, A_{r-1}, \dots, A_1, 1, 1]$ with $A_i | A_{i+1}$ so that a possible choice for $\alphabar$ is given by:

\begin{equation}\label{comatrixchoice}
\begin{pmatrix} \abar & c \end{pmatrix} \begin{pmatrix} \alphabar \\ A_r\gamma\end{pmatrix} = A_rI_{r+2}
\end{equation}

To express the transformation properties of the function $G_{r, \abar}$ under the action of $\slnz{r+2}$ we need to introduce the family of geometric Bernoulli polynomials which encompass $P_{a,b,c}$. In what follows, we want to discuss an analogue of formula \refp{cocyclegamma} for $G_{r,\abar}$, and we must consider families of $r+2$ vectors. Recall that we have defined the polynomials $B_{r+2, r+2}^{*}$ in section \ref{Bernoulli}.

\begin{general}{Definition}\label{defgeombern}
Let $w \in \cc$ and $x \in \homlc \iso \cc^{r+2}$. Let $a_0, \dots, a_r, a_{r+1}$ be a family of $r+2$ linearly independent primitive vectors in the oriented lattice $\Lambda$. Let $\epsilon$ be the sign of $(-1)^{r+1}\det(a_0, \dots, a_{r+1})$ and put
$$B^{*}_{r+2, a_0, \dots, a_{r+1}}(w, x, \Lambda) = \frac{-2\epsilon}{(r+2)!} \sum_{\delta \in F} B^{*}_{r+2,r+2}(w + \delta(x), \alpha_0(x), \dots, \alpha_{r+1}(x))$$ 
where $\alpha_0, \dots, \alpha_{r+1}$ is the unique primitive positive dual family to $(a_0, \dots, a_{r+1})$ in $L = \homlambdaz$ (see lemma \ref{lemmaunique}) and
$$F = F(\alpha_0, \dots, \alpha_{r+1}) = \{\delta \in L, 0 \leq \delta(a_j) < \alpha_j(a_j), \forall\,0 \leq j \leq r+1 \} $$ 
is a finite set.
\end{general}

\noindent These geometric Bernoulli ``polynomials'' have the form
$$B^{*}_{r+2, a_0, \dots, a_{r+1}}(w, x, \Lambda) = \sum_{l + k_0 + \dots + k_{r+1} = r+2} c_{l, k_0, \dots, k_{r+1}} w^l \prod_{j = 0}^{r+1} \alpha_j(x)^{k_j}$$
where the sum is taken over non-negative integers and the coefficients $c_{l, k_0, \dots, k_{r+1}}$ are rational numbers which explicitly depend on the set $F = F(\alpha_0, \dots, \alpha_{r+1})$. The modularity property for the geometric $G_{r, \abar}$ functions will involve the rescaled versions of the geometric Bernoulli polynomials 
$$B_{r+2, a_0, \dots, a_{r+1}}(w, x, \Lambda) = \left(\prod_{j = 0}^{r+1} \alpha_j(x)\right)^{-1}B^{*}_{r+2,a_0, \dots, a_{r+1}}(w, x, \Lambda)$$
For $r=1$ the polynomial $B_{3, a, b, c}$ may be identified with the polynomial $P_{a,b,c}$ appearing in formula \refp{cocyclegamma}. We can now generalise formula (\ref{cocyclegamma}).

\begin{general}{Theorem}\label{theoremslnz}
\begin{enumerate}
\item $[\,$Modular property$\,]$ For a family $a_0, \dots, a_{r+1}$ of $r+2$ linearly independent primitive vectors in $\Lambda$:
\begin{equation}\label{modulargeom}
\prod_{j= 0}^{r+1} G_{r, (a_k)_{k \neq j}}(w, x, \Lambda)^{(-1)^j} = \exp(i\pi B_{r+2,a_0, \dots, a_{r+1}}(w,x, \Lambda))\end{equation}
\item $[\,$Equivariance relations$\,]$ For a family $a_0, \dots, a_r$ of $r+1$ linearly independent primitive vectors in $\Lambda$, for $w \in \cc$ and $x \in \homlc$, for all $g \in \slnz{r+2}$:
$$G_{r,g\abar}(w, gx, \Lambda) = G_{r,\abar}(w, x, \Lambda) $$
$$B_{r+2,ga_0, \dots, ga_{r+1}}(w, gx, \Lambda) = B_{r+2,a_0, \dots, a_{r+1}}(w, x, \Lambda) $$
\end{enumerate}
\end{general}

\begin{proof}
1. The modular relation for the geometric $G_{r, \abar}$ functions is a consequence of the modular property \refp{modGr} for the ordinary $G_r$ function. Let $\alphabar$ be the primitive positive dual family to $\abar$ in $L$ as given in lemma \ref{lemmaunique}. Write $\epsilon$ for the sign of $\det_B(a_0, \dots, a_{r+1}).(-1)^{r+1}$. Then for all $0 \leq j \leq r+1$ we have 
$$\det_B(a_0, \dots, a_{j-1}, a_{j+1}, \dots, a_{r+1}, \cdot) = s_j.\epsilon.(-1)^j\alpha_j$$
so that $\gamma_j = \epsilon.(-1)^j\alpha_j$ and
$$G_{r, (a_k)_{k \neq j}}(w, x, \Lambda) = \prod_{\delta \in F_j/\zz\alpha_j}G_r\left(\frac{w + \delta(x)}{\alpha_j(x)}, \left(\frac{\alpha_k(x)}{\alpha_j(x)}\right)_{k \neq j}\right)^{(-1)^j.\epsilon}$$
where $$F_j = \{\delta \in L, 0 \leq \delta(a_k) < \alpha_k(a_k), \forall\, 0 \leq k \neq j \leq r+1 \} $$ 
Then
$$\prod_{j= 0}^{r+1} G_{r, ((a_k)_{k \neq j}}(w, x, \Lambda)^{(-1)^j} = \prod_{j= 0}^{r+1}  \prod_{\delta \in F_j/\zz\alpha_j}G_r\left(\frac{w + \delta(x)}{\alpha_j(x)}, \left(\frac{\alpha_k(x)}{\alpha_j(x)}\right)_{k \neq j}\right)^{\epsilon}$$
Put $\mathcal{F} = F(\alpha_0, \dots, \alpha_{r+1}) = \{\delta \in L, 0 \leq \delta(a_j) < \alpha_j(a_j), \forall\,0 \leq j \leq r+1\}$ so that we can write uniformly $\mathcal{F} \iso F_j/\zz\alpha_j$ for all $0 \leq j \leq r+1$. Then, using Narukawa's theorem (\ref{modGr}) for each $\delta$ in the finite set $\mathcal{F}$ yields:
$$\prod_{j= 0}^{r+1} G_{r, ((a_k)_{k \neq j}}(w, x, \Lambda)^{(-1)^j} = \prod_{\delta \in \mathcal{F}}
\exp\left(\frac{-2i\pi\epsilon}{(r+2)!}B_{r+2, r+2}(w+\delta(x), \alphabar(x))\right)$$
The identification of the right-hand side gives the conclusion:
$$\prod_{j= 0}^{r+1} G_{r, ((a_k)_{k \neq j}}(w, x, \Lambda)^{(-1)^j} =\exp(i\pi B_{r+2, a_0, \dots, a_{r+1}}(w, x, \Lambda))$$

2. Consider $g \in \slnz{r+2}$. In the construction of $G_{r, \abar}$ or $B_{r, a_0, \dots, a_{r+1}}$ replacing $\abar$ with $g \abar$ replaces $\alphabar$ with $\alphabar g^{-1}$, $\gamma$ with $\gamma g^{-1}$ and $F(\alphabar)$ with $F(\alphabar g^{-1})$ but $\epsilon$ is left unchanged. Therefore we may write:
$$G_{r, g\abar}(w, gx, \Lambda) = \prod_{\delta \in F(\alphabar g^{-1})/\zz\gamma g^{-1}} G_r\left(\frac{w + \delta(gx)}{\gamma g^{-1}(gx)}, \frac{1}{\gamma g^{-1}(gx)} \alphabar g^{-1}(gx) \right) $$
Then, identifying $F(\alphabar g^{-1}) = F(\alphabar)g^{-1}$ and putting $\delta' = \delta g$ gives
$$G_{r, g\abar}(w, gx, \Lambda) = \prod_{\delta' \in F(\alphabar)/\zz\gamma} G_r\left(\frac{w + \delta'(x)}{\gamma (x)}, \frac{1}{\gamma (x)} \alphabar (x) \right) $$
$$G_{r, g\abar}(w, gx, \Lambda) = G_{r, \abar}(w, x, \Lambda)$$
As for $B_{r+2, ga_0, \dots, ga_{r+1}}$, write $\mathcal{F} = F(\alpha_0, \dots, \alpha_{r+1})$:
$$B_{r+2, ga_0, \dots, ga_{r+1}}(w, gx, \Lambda) = \frac{-2\epsilon}{(r+2)!} \sum_{\delta \in \mathcal{F}g^{-1}} B_{r+2,r+2}(w + \delta(gx), \alphabar(g^{-1}gx))$$
Put once again $\delta' = \delta g$, which gives:
$$B_{r+2, ga_0, \dots, ga_{r+1}}(w, gx, \Lambda) = \frac{-2\epsilon}{(r+2)!} \sum_{\delta' \in F(\alpha_0, \dots, \alpha_{r+1})} B_{r+2,r+2}(w + \delta(x), \alphabar(x))$$
Identify the right-hand side to conclude that
$$B_{r+2, ga_0, \dots, ga_{r+1}}(w, gx, \Lambda) = B_{r+2, a_0, \dots, a_{r+1}}(w, x, \Lambda)$$
\end{proof}

We also add that the geometric $G_{r, \abar}$ functions behave nicely under permutation of vectors, namely for all permutation $\sigma \in \goth{S}_{r+1}$, $G_{r, \sigma(\abar)} = G_{r, \abar}^{\eps(\sigma)}$ where $\eps(\sigma)$ is the sign of the permutation $\sigma$. As in the cubic case, the relation \refp{modulargeom} can be restated and best understood as a multiplicative $r$-cocycle relation for $\slnz{r+2}$. For $g_1, \dots, g_l \in \slnz{r+2}$ put
$$[g_1|\dots| g_l] = \left(g_1, \dots, \prod_{j = 1}^{k}g_j, \dots, \prod_{j=1}^{l}g_j\right)$$
Then we may write $\psi_a(g_1, \dots, g_r) = G_{r, a, g_1a, \dots, (\prod_{j = 1}^r g_j)a} = G_{r, [1| g_1| \dots| g_r]a}$ for a primitive vector $a \in \Lambda$ and $g_1, \dots, g_r \in \slnz{r+2}$. The coboundary of $\psi_a$ is given by the modular property \refp{modulargeom}:
\begin{equation} \label{cocycler}
d\psi_a(g_1, \dots,  g_{r+1}) = \exp(i\pi B_{r+2, [1|g_1|\dots|g_{r+1}]a})
\end{equation}
where $1$ stands for the identity matrix $I_{r+2}$. This property is essential to understand the type of evaluation of the $G_r$ functions we carry out in section \ref{arithmeticsection}. Indeed, in \cite{Sczech} Sczech uses Eisenstein $(d-1)$-cocycles for $\slnz{d}$ to express the values of $L$-functions of totally real fields of degree $d$ at negative integers. In the case of a totally real field $\kk$ of degree d, fixing a basis of $\ok$, the unit group $\ok^{\times}$ may be embedded inside $\slnz{d}$ as an abelian group of rank $d-1$. In the case of ATR fields which we are concerned with in this article, we have here obtained a $(d-2)$-cocycle $\psi_a$ (up to the $\exp(i\pi B_{r+2, [1|g_1|\dots|g_{r+1}]a})$ term) for $\slnz{d}$ which is compatible with the unit group of an ATR field of degree $d$ having rank $d-2$.

\section{An algebraicity conjecture on special values of $G_r$ functions}\label{sectionarithmetic}

\subsection{Arithmetic $G_r$ functions}\label{arithmeticsection}

In their recent article \cite{BCG}, Bergeron, Charollois and Garc\'ia describe special values of the elliptic Gamma function evaluated at specific points in complex cubic fields which they relate to partial zeta values. They also conjecture that these special values are units lying in specific abelian extensions. We shall present a generalisation of this construction for fields of degree $d \geq 2$ with exactly one pair of complex embeddings yielding special values of $G_{d-2}$ functions.

Consider a degree $d \geq 2$ ATR number field $\kk$. Fix a complex embedding $\sigma_{\cc}$ of $\kk$ and fix an ordering $\sigma_1, \dots, \sigma_{d-2}, \sigma_{d-1} = \sigma_{\cc}, \sigma_d = \overline{\sigma_{d-1}}$ on the embeddings of $\kk$. We fix a global orientation of $\kk$ as follows. For any basis $B = (e_1, \dots, e_d)$ of the $\qq$-vector space $\kk$, the determinant $i\cdot\det((\sigma_j(e_k))_{1 \leq j,k\leq d})$ is a non zero real number. We call a basis $B$ positive if the latter is positive. Consider an integral ideal $\goth{f} \neq \ok$. We write $q\zz = \zz \cap \goth{f}$. Fix an integral ideal $\goth{b}$ coprime to $\goth{f}$ and put $L = \goth{f}\goth{b}^{-1}$. Fix an integral ideal $\goth{a}$ coprime to $\goth{f}\goth{b}$ such that $\goth{a}^{-1}L/L$ is a cyclic group. We call such an ideal a \textit{smoothing} ideal and we write $N=\norm{\goth{a}}$ for the norm of $\goth{a}$. We consider specific vectors in $L$ which we will use to evaluate geometric $G_r$ functions.

\begin{general}{Definition}\label{defadmissible}
A vector $h \in L = \goth{f}\goth{b}^{-1}$ is called admissible for the data $\goth{f}, \goth{b}, \goth{a}$ if $h/q \equiv 1 \mathrm{~mod~} L$ and $h/N$ generates the cyclic group $\goth{a}^{-1}L/L$, where $q\zz = \zz \cap \goth{f}$ and $N= \norm{\goth{a}}$ is the norm of the ideal $\goth{a}$.
\end{general}

\noindent In this section we fix $h \in L$ an admissible vector and we fix a $\zz$-basis $B_L$ of the $\zz$-module $L$ which is a positive basis of $\kk$. The unit group
\begin{equation}\label{defopcf}
\opc{\goth{f}} = \{\eps \in \units, \eps \equiv 1 \mod \goth{f}, \sigma(\eps) > 0, \mathrm{for~all~real~embedding~} \sigma~\mathrm{of}~\kk\}
\end{equation}
is a free $\zz$-module of rank $r = d-2$ acting on $L$, and we fix for now $u_1, \dots, u_r$ a set of fundamental units for $\opc{\goth{f}}$ such that the linear form $\det_{B_L}(h, u_1 h, \dots, u_r h, \cdot)$ is non zero. Denote by $a$ the unique primitive element in $\homlz$ satisfying $na = \det_{B_L}(h, u_1 h, \dots, u_r h, \cdot)$ for a positive integer $n$. 

\begin{general}{Definition}\label{defgrarithmetic}
Let $u_1, \dots, u_r$ be a set of fundamental units for $\opc{\goth{f}}$ as above. Let $h$ be an admissible vector in $\goth{f}\goth{b}^{-1}$ and $a$ be the linear form described above. We define for the complex embedding $\sigma_{\cc}$ of $\kk$ and the choice of sign $\pm$ the arithmetic $G_r$ function:
$$G_{r, \goth{f}, \goth{b}, \goth{a}}^{+}(u_1, \dots, u_r ; h, \sigma_{\cc}) = \frac{G_{r, (a, au_1, \dots, au_r)}(\frac{h}{q}, \sigma_{\cc},\mathrm{Hom}_{\zz}(\goth{f}\goth{b}^{-1}, \zz))^N}{G_{r, (a, au_1, \dots, au_r)}(\frac{h}{q}, \sigma_{\cc}, \mathrm{Hom}_{\zz}(\goth{f}(\goth{b}\goth{a})^{-1},\zz))}$$
$$G_{r, \goth{f}, \goth{b}, \goth{a}}^{-}(u_1, \dots, u_r ; h, \sigma_{\cc}) = \frac{G_{r, (- a, - au_1, \dots, -au_r)}(\frac{h}{q}, \sigma_{\cc},\mathrm{Hom}_{\zz}(\goth{f}\goth{b}^{-1}, \zz))^N}{G_{r, (-a, -au_1, \dots, -au_r)}(\frac{h}{q}, \sigma_{\cc}, \mathrm{Hom}_{\zz}(\goth{f}(\goth{b}\goth{a})^{-1},\zz))}$$
where in both cases the right-hand side is a quotient of two geometric $G_r$ functions as defined in Proposition \ref{propdefgeom}.
\end{general}

\textbf{Remark} This is well-defined provided that a positive dual family $\alphabar$ to $a, au_1, \dots, au_r$ satisfies $\sigma_{\cc}(\alpha_j/h) \not\in \rr$ for all $0 \leq j \leq r$. This is always the case, except possibly when $\kk$ is of even degree $2d'$ and contains a real subfield $\mathbb{F}$ of degree $\leq d'$. We will briefly discuss this case in section \ref{quarticreal}. 

One may use formula \refp{grconjugate} to relate the two values $G_{r, \goth{f}, \goth{b}, \goth{a}}^{\pm}(u_1, \dots, u_r ; h, \sigma_{\cc})$ and $G_{r, \goth{f}, \goth{b}, \goth{a}}^{\pm}(u_1, \dots, u_r ; h, \overline{\sigma_{\cc}})$ for the two complex embeddings of $\kk$ as:
$$G_{r, \goth{f}, \goth{b}, \goth{a}}^{\pm}(u_1, \dots, u_r ; h, \overline{\sigma_{\cc}}) = \overline{G_{r, \goth{f}, \goth{b}, \goth{a}}^{\pm}(u_1, \dots, u_r ; h, \sigma_{\cc})}^{(-1)^r}$$
This means that we only need to treat the case of one of the complex embeddings. Throughout the rest of this work, we will consider fields with a specific defining polynomial equation and we fix the complex embedding $\sigma_{\cc}$ corresponding to the root of the polynomial lying in the upper half-plane. We shall then drop the $\sigma_{\cc}$ in the writing. In the special cases $r = 0, \,r = 1$ we will identify
$$\theta_{\goth{f}, \goth{b}, \goth{a}}^{\pm}(h) = G_{0,\goth{f}, \goth{b}, \goth{a}}^{\pm}(\emptyset; h)$$
$$\Gamma_{\goth{f}, \goth{b}, \goth{a}}^{\pm}(\eps; h) = G_{1, \goth{f}, \goth{b}, \goth{a}}^{\pm}(\eps; h)$$

In \cite{BCG} the authors use modular symbols (see the cocycle relation \refp{cocycleone}) for $\slnz{3}$ built from smoothed Eisenstein series to establish a link between this special value and the partial zeta values associated to the conductor $\goth{f}$ and the ideal class $[\goth{b}]$ of the integral ideal $\goth{b}$ in the narrow ideal class group $Cl^{+}(\goth{f})$. Namely, the partial $\zeta$ function associated to this class is
$$\zeta_{\goth{f}}([\goth{b}], s) = \sum_{I \in [\goth{b}]} \norm{I}^{-s}$$ 
when $\Re(s) > 1$ and where the sum ranges over all integral ideals in the class $[\goth{b}]$. This function can be extended into a meromorphic function over $\cc$ with a pole at $s = 1$. In the specific case where $\kk$ is a complex cubic field, this partial $\zeta$ function vanishes at the origin. In \cite{BCG} the authors have proven that
when $\opcf = \eps^{\zz}$ with $\sigma_{\rr}(\eps) < 1$ then the modulus of $\Gamma_{\goth{f}, \goth{b}, \goth{a}}^{-}(\eps, h)^{-1} = \Gamma_{\goth{f}, \goth{b}, \goth{a}}^{+}(\eps^{-1}, h)$ is independent of the choice for the admissible vector $h$ and the following Kronecker type limit formula holds:
$$\norm{\goth{a}}\zeta'_{\goth{f}}([\goth{b}], 0) - \zeta'_{\goth{f}}([\goth{a}\goth{b}], 0) = \log\left|\Gamma_{\goth{f}, \goth{b}, \goth{a}}^{-}(\eps, h)^{-1}\right|^2 = \log\left|\Gamma_{\goth{f}, \goth{b}, \goth{a}}^{+}(\eps^{-1}, h)\right|^2 $$

They also conjectured that $\Gamma_{\goth{f}, \goth{b}, \goth{a}}^{-}(\eps, h)^{-1}$ is a unit in the associated class field $\kk^{+}(\goth{f})$ related to the conjectural Stark unit. We believe that a combination of special values of arithmetic $G_r$ functions for $r \geq 2$ will be linked to the partial $\zeta$ functions of an ATR field of degree $d = r+2$ using modular symbols (see the cocycle relation \refp{cocycler}) for $\slnz{d}$ related to some smoothed Eisenstein series. In view of this cohomological interpretation and following the work of Sczech in \cite{Sczech} on the Eisenstein cocycles for totally real fields, we consider the following combination of arithmetic $G_r$ functions:

\begin{general}{Definition}\label{defIrfba}
Let $(\eps_1, \dots, \eps_r)$ be a set of fundamental units for $\opc{\goth{f}}$. Let $\bars{h} = (h_{\rho})$ be a set of admissible vectors for $\goth{f}, \goth{b}, \goth{a}$ where $\rho$ ranges over the permutations in $\goth{S}_r$, the symmetric group on $r$ elements. For any choice of orientations for the $G_r$ functions $\bars{\mu} = (\mu_{\rho}), \bars{\nu} = (\nu_{\rho})$ we define:

$$\spvgrcomplete = \prod_{\rho \in \goth{S}_{r}} G_{r, \goth{f}, \goth{b}, \goth{a}}^{\mu_{\rho}}([\eps_{\rho(1)} | \dots | \eps_{\rho(r)}] ; h_\rho)^{\nu_{\rho}}$$
where again
$$[\eps_1 | \dots | \eps_r] = \left(\eps_1, \dots, \prod_{j = 1}^k\eps_j, \dots, \prod_{j=1}^r\eps_j\right)$$
and the orientations $\mu_{\rho}, \nu_{\rho} \in \{-1, +1\}$. This is well-defined if for all permutation $\rho$, the family $(1, \eps_{\rho(1)}, \dots, \prod_{j = 1}^r \eps_{\rho(j)})$ is a family of linearly independent vectors in $\kk$.
\end{general}

We will now present a conjecture on special values of these arithmetic $G_r$ functions for specific choices of admissible vectors $(h_{\rho})_{\rho \in \goth{S}_r}$ and orientations $(\mu_{\rho})_{\rho \in \goth{S}_r}, (\nu_{\rho})_{\rho \in \goth{S}_r}$ which we will describe in the next section.

\subsection{Formulation of the conjecture}\label{formulation}

One aim of this article is to generalise the following conjecture of Bergeron, Charollois and Garc\'ia \cite{BCG} on the special values of the elliptic Gamma function evaluated at points in complex cubic fields.

\begin{famous}{Conjecture (Bergeron, Charollois and Garc\'ia)}
Suppose that $\kk$ is a complex cubic field and $\goth{f} \neq \ok$ is such that $\kk^{+}(\goth{f})$ is totally complex. Take $\goth{b}, \goth{a}$ as in section \ref{arithmeticsection} and further ask that $(\mathcal{N}(\goth{a}), 6) = 1$. Define $\eps$ to be the generator of $\opcf$ satisfying $\sigma_{\rr}(\eps) < 1$. Then
for any admissible vector $h \in L = \goth{f}\goth{b}^{-1}$:
\begin{itemize} 
\item the number $\Gamma_{\goth{f}, \goth{b}, \goth{a}}^{-}(\eps, h)^{-1} = \Gamma_{\goth{f}, \goth{b}, \goth{a}}^{+}(\eps^{-1}, h)$ is independent of the choice of $h$ and is the image in $\cc$ of a unit $u_{L, \goth{a}} \in \kk^{+}(\goth{f})$ under a complex embedding $\sigma'_{\cc}$ extending $\sigma_{\cc}$. 
\item any complex embedding of $\,\kk^{+}(\goth{f})$ above the real place of $\,\kk$ sends $u_{L, \goth{a}}$ to the unit circle. 
\item for $\goth{c} \in I(\goth{f})$, the reciprocity law is explicitly given by $u_{L, \goth{a}}^{\sigma_{\goth{c}}} = u_{L\goth{c}^{-1}, \goth{a}}$ where $\sigma_{\goth{c}}$ is the corresponding element in $\mathrm{Gal}(\kk^{+}(\goth{f})/\kk)$.
\end{itemize}
\end{famous}
 
Consider the following example for the field $\kk = \qq(z)$ where $z = e^{2i\pi/3}10^{1/3}$ is the root of the polynomial $x^3-10$ in the upper half-plane, $\goth{f}^3 = (5)$, $\goth{b} = (1)$ and $\goth{a}$ is the unique degree one prime above $11$ in $\kk$. The fundamental unit is $\eps = (2z^2 - z - 7)/3$. Then, for $h = -(35z^2 + 20z + 35)/3$ we compute
$$\Gamma^{-}_{\goth{f}, \goth{b}, \goth{a}}(\eps, h)^{-1} = \frac{\Gamma\left(\frac{-1}{5}, \frac{\eps^{-1}-1751}{495}, \frac{\eps-776}{495}\right)^{-11}}{\Gamma\left(\frac{-11}{5}, \frac{\eps^{-1}-1751}{45}, \frac{\eps-776}{45}\right)^{-1}} \approx -27.5333588... - i\cdot32.7146180...$$
with 1000 digits of precision and find it to be close to a root of the polynomial $x^{12} + 57x^{11} + 1956x^{10} + 4640x^9 + 35415x^8 - 109818x^7 + 150139x^6 - 109818x^5 + 35415x^4 + 4640x^3 + 1956x^2 + 57x + 1$ which defines $\kk^{+}(\goth{f})$. This example is developped in section \ref{excubictext} below.

In the cubic case, this conjecture has been tested numerically on hundreds of examples and we want to generalise this to ATR fields of degree $d \geq 4$ in some specific cases supported by numerical evidence. In section \ref{sectioncomputing}, we will describe a specific set of conditions under which we formulate our conjecture. To state it properly, we must discuss the choice of orientations, that is to say the choice of signs $\mu_{\rho}, \nu_{\rho}$. The cycle we use to evaluate our $G_r$ functions mimics the cycle used in \cite{Colmez}, \cite{Diazydiaz} in the context of explicit Shitani cones (or signed Shintani cones) in totally real fields, as well as in \cite{Sczech} in the context of the Eisenstein cocycle for totally real fields. In the context of number fields with exactly one complex place, explicit signed Shintani cones have been built in \cite{Espinoza}. Therefore, we choose orientations to mimic the orientations given in \cite{Colmez}. Remember that we have fixed a set $\eps_1, \dots, \eps_r$ of fundamental units for $\opc{\goth{f}}$. We have also fixed the ordering of the real embeddings so we may consider the sign of the regulator of these units.

\begin{general}{Definition}\label{deforientations}
The sign of the unit system $\eps_1, \dots, \eps_r$ is defined by:
$$sg(\eps_1, \dots, \eps_r) = sign(\det(\log(\sigma_{j}(\eps_k))_{1 \leq j,k \leq r}))$$
Then, we define the following orientations
$$\mu_{\rho} = \nu_{\rho} = sg(\eps_1, \dots, \eps_r) \signature(\rho)$$
where $\signature(\rho)$ is the signature of the permutation $\rho$.
\end{general}

\begin{general}{Conjecture}\label{conjecture}
Suppose that $\kk$ is an ATR field of degree $d = r+2$. Suppose $\goth{f} \neq \ok$ is an integral ideal such that $\kk^{+}(\goth{f})$ is totally complex. Suppose further that all units $\eps \in \ok^{\times}$ that are congruent to $1$ mod $\goth{f}$ have norm $+1$. Take $\goth{b}, \goth{a}$ two integral ideals as in Proposition \ref{propfinal}. Suppose that there exists a system of fundamental units $\eps_1, \dots, \eps_r$ satisfying the conditions $(\ast)$ and $(\ast\ast)$ associated to $\goth{f}, \goth{b}, \goth{a}$ described in section \ref{choiceh}. For any family of vectors $\bars{h}$ described in section \ref{choiceh}, Proposition \ref{propfinal}, and the orientations $\bars{\mu}, \bars{\nu}$ given in Definition \ref{deforientations}, the following Kronecker type limit formula holds: 
\begin{equation}\label{eqconjecture}
\norm{\goth{a}}\zeta'_{\goth{f}}([\goth{b}], 0) - \zeta'_{\goth{f}}([\goth{a}\goth{b}], 0) = \frac{1}{\cardinalshort{\mathcal{Z}_{\goth{f}}^1}}\log\left|\prod_{(n,\eps) \in \mathcal{Z}_{\goth{f}}^1}I_{r,\goth{f}, \goth{b}, \goth{a}}(\eps_1, \dots, \eps_r ; n\eps\bars{h}, \bars{\mu}, \bars{\nu})\right|^2
\end{equation}
where the finite set $\mathcal{Z}_{\goth{f}}^1$ is defined in section \ref{sectionchoiceh}, Definition \ref{defzfone}. In addition:
\begin{enumerate}
\item The product in the right-hand side is independent of the choice of vectors $\bars{h}$ in Proposition \ref{propfinal} up to a root of unity of order $m(r,N)$, where the number $m(r,N)$ is defined by $m(r,2) = 4$, $m(r, N) = N$ for $3 \leq N \leq r+3$ and $m(r,N) = 1$ for $N > r+3$. 
\item If $r \geq 1$, the product in the right-hand side, raised to the power $m(r,N)$ is the image in $\cc$ of an algebraic unit $\ula$ in the class field $\kk^{+}(\goth{f})$ under a complex embedding $\sigma'_{\cc}$ extending $\sigma_{\cc}$.
\item Any embedding of $\kk^{+}(\goth{f})$ above a real embedding of $\kk$ sends $\ula$ to the unit circle.
\end{enumerate}
\end{general}

\noindent \textbf{Remarks : } \begin{enumerate} 
\item The case $r= 0$ in the conjecture is already known as a consequence of the theory of complex multiplication and we refer to \cite{Robert} for a presentation on the subject. Note that in this case the value
$$\theta_{\goth{f}, \goth{b}, \goth{a}}^{\pm}(h)^{m(0, N)}$$
is not necessarily a unit but rather a $q$-unit where once again $q\zz = \goth{f} \cap \zz$. In the case $r = 1$, formula \refp{eqconjecture} was proven in \cite{BCG} and the rest of the conjecture was also formulated in \cite{BCG}.
\item In the conjecture, the condition that all units congruent to $1 \mod \goth{f}$ should have norm $+1$ is not restrictive. Indeed, if a unit $\eps$ existed with $\eps \equiv 1 \mod \goth{f}$ and $\norm{\eps} = -1$ then all $L$-functions attached to the extension $\kk^{+}(\goth{f})/\kk$ would have a zero at $0$ of order at least $2$. Therefore, the left-hand side in formula \refp{eqconjecture} would vanish.
\item This conjecture is formulated with restrictions, but it has been successfully tested in many other cases outside of the conditions ($\ast$), ($\ast\ast$). It is our aim to complete the formulation of this conjecture outside of these restrictions and to extend the range of ``compatible'' vectors $\bars{h}$ for which the conjecture holds. Once completed, it should be part of the conjecture that the explicit reciprocity law is given by $\ula^{\sigma_{\goth{c}}} = u_{L\goth{c}^{-1}, \goth{a}}$ for an integral ideal $\goth{c}$ in $I(\goth{f})$.
\item Our computations show a clear change of behavior between the smoothings with $N \leq r+3$ and the smoothings with $N > r+3$. In the case of ``small smoothings'', the product in the right-hand side of \refp{eqconjecture} should be an algebraic unit inside $\kk^{+}(\goth{f})(\exp(2i\pi/m(r,N))$ whereas in the case of ``big smoothings'', it should be an algebraic unit inside $\kk^{+}(\goth{f})$. We expect the exact bound $r+3$ because of the von Staudt-Clausen theorem on the denominator of Bernoulli numbers: a prime divisor $p$ of the denominator of one of the $B_1, B_2, \dots, B_{r+2}$ must satisfy $p \leq r+3$. The value of $m(r,N)$ has been tested on many examples for $0 \leq r \leq 3$. The table below presents the expected values for $m(r,N)$ in these cases.

\begin{center}
\begin{tabular}{| C | C | C | C | C |}
\hline
N & \mathrm{quartic} & \mathrm{cubic} & \mathrm{quartic} & \mathrm{quintic}\\
& m(0, N) & m(1, N) & m(2, N) & m(3, N) \\
\hline
2 & 4 & 4 & 4 & 4\\
3 & 3 & 3 & 3 & 3\\ 
5 & 1 & 1 & 5 & 5\\
7\leq & 1 & 1 & 1 & 1\\
\hline
\end{tabular}
\end{center}
\end{enumerate}

In practice, the product in the conjecture often contains redundant terms, and it is thus already a power inside $\kk^{+}(\goth{f})$. We therefore define the integer $\kappa$ as the maximal integer $k \geq 1$ dividing $\cardinalshort{\mathcal{Z}_{\goth{f}}^1}$, such that there exist subsets $(U_{\rho, j})_{\rho \in \goth{S}_r, 1 \leq j \leq k}$ of $\mathcal{Z}_{\goth{f}}^1$ satisfying
\begin{itemize}
\item[(i)] for all $\rho \in \goth{S}_r$, for all $1 \leq j \leq k$, $\cardinalshort{U_{\rho, j}} = \cardinalshort{\mathcal{Z}_{\goth{f}}^1}/k$
\item[(ii)] for all $\rho \in \goth{S}_r$, $j \neq l$, $U_{\rho, j} \cap U_{\rho, l} = \emptyset$.
\item[(iii)] for all $1 \leq j \leq k$, the product
$$\left(\prod_{\rho \in \goth{S}_r} \prod_{(n, \eps) \in U_{\rho, j}}G_{r,\goth{f}, \goth{b}, \goth{a}}^{\mu_{\rho}}(\eps_1, \dots, \eps_r ; n\eps\bars{h})^{\nu_{\rho}} \right)^{m(r,N)}$$ is independent of $j$.
\end{itemize}
Conditions (i) and (ii) express the fact that the sets $U_{\rho, j}$ form a special partition of the set $(\mathcal{Z}_{\goth{f}}^1)^{r!}$. Condition (iii) implies in particular that 
$$\prod_{(n, \eps) \in \mathcal{Z}_{\goth{f}}^1}I_{r,\goth{f}, \goth{b}, \goth{a}}(\eps_1, \dots, \eps_r ; n\eps\bars{h}, \bars{\mu}, \bars{\nu}) = \left(\prod_{\rho \in \goth{S}_r} \prod_{(n, \eps) \in U_{\rho, 1}}G_{r,\goth{f}, \goth{b}, \goth{a}}^{\mu_{\rho}}(\eps_1, \dots, \eps_r ; n\eps\bars{h})^{\nu_{\rho}}\right)^{\kappa}$$ 
is already a $\kappa$-th power of what is expected to be an algebraic unit in $\cc$. In the cubic case, it is expected that $\kappa = \cardinalshort{\mathcal{Z}_{\goth{f}}^1}$ and this is expected to be true in most simple examples for $d \geq 4$. However, the worst case $\kappa = 1$ with $\cardinalshort{\zfone} > 1$ already appears for $d = 4$ (see section \ref{lowkappaex}).

We give here a quartic example supporting the conjecture with a ``small'' smoothing: take the quartic field $\kk = \qq(z)$ where $z$ is the root of the polynomial $x^4 -6x^3-x^2-3x+1$ in the upper half-plane and fix the ordering on the real embeddings of $\kk$ such that $\sigma_1(z) < \sigma_2(z)$. Take $\goth{f}$ the unique degree one prime above $q = 2$ in $\kk$. A possible choice for the fundamental units is given by $\eps_1 = (-2z^3 + 13z^2 - z + 3)/7,~\eps_2 = (-5z^3 + 29z^2 + 15z + 18)/7 $. We choose $\goth{a}$ the degree one prime above $N = 5$ in $\kk$ and $\goth{b} = \ok$. Here, $q = 2$ therefore $\zfone = \{(1,1)\}$ and $\kappa = \cardinalshort{\zfone} = 1$. A possible choice of vectors $\bars{h}$ as in Proposition \ref{propfinal} is $h_1 = h_2 = (-18z^3 + 96z^2 + 82z + 62)/7$. We may compute the orientations $\bars{\mu} = \bars{\nu} = [-1, 1]$ as in Definition \ref{deforientations} and the two quotients
\begin{align*}
 v_1 & = \frac{G_2\left(\frac{1}{2}, \frac{5z^3 - 29z^2 - 15z-81}{70}, \frac{6z^3 - 39z^2 + 10z + 5}{70}, \frac{-2z^3 + 13z^2 - z + 24}{70}\right)^{-5}}{G_2\left(\frac{5}{2}, \frac{5z^3 - 29z^2 - 15z-81}{14}, \frac{6z^3 - 39z^2 + 10z + 5}{14}, \frac{-2z^3 + 13z^2 - z - 24}{14}\right)^{-1}}, \\
v_2 & = \frac{G_2\left(\frac{-1}{2}, \frac{-2z^3 +13z^2 - z + 24}{70}, \frac{-5z^3 + 29z^2 + 15z + 81}{70}, \frac{2z^3 - 13z^2 - 6z - 101}{70}\right)^{5}}{G_2\left(\frac{-5}{2}, \frac{-2z^3 +13z^2 - z + 24}{14}, \frac{-5z^3 + 29z^2 + 15z + 81}{14}, \frac{2z^3 - 13z^2 - 6z - 101}{14}\right)}
\end{align*}

\noindent The values obtained are respectively $v_1 \approx -2.0167576... - i\cdot5.8008598...$ and $v_2 \approx -0.4159958... + i\cdot0.0018434...$ Their product $v_1v_2 \approx 0.8496565... - i\cdot2.4094157...$ is not an algebraic integer in $\kk^{+}(\goth{f})$. Yet the \textit{fifth} power $(v_1v_2)^5 \approx 108.0070738... -i\cdot13.4979021...$ of this product coincides up to at least 1000 digits with a root of the polynomial $x^8 - 215x^7 + 11629x^6 + 11941x^5 + 3913x^4 + 11941x^3 + 11629x^2 - 215x + 1$ which defines $\kk^{+}(\goth{f})$. The constant term of this polynomial is $1$, so the roots of this polynomial are units inside $\kk^{+}(\goth{f})$. We may also check formula \refp{eqconjecture} up to 1000 digits as:
$$\norm{\goth{a}}\zeta'_{\goth{f}}([\goth{b}], 0) - \zeta'_{\goth{f}}([\goth{a}\goth{b}], 0) \approx \log|v_1v_2|^2 \approx 1.8759781...  $$ 
We will present an example in the same conditions except for the smoothing which will be a ``big smoothing'' in section \ref{quarticsimple}.

To end this section, we point out that the way our conjecture was formulated is also inspired by the rank one abelian Stark conjecture in the ATR case:

\begin{famous}{Rank one abelian Stark conjecture (ATR case, see \cite{rankoneStark})}
Write $e_{\goth{f}}$ for the number of roots of unity in $\kk^{+}(\goth{f})$. There is a unit $u_{Stark}$ in $\kk^{+}(\goth{f})$ such that 
\begin{itemize}
\item for all class $[\goth{b}]$ in $Cl^{+}(\goth{f})$ and corresponding $\sigma_{\goth{b}} \in \mathrm{Gal}(\kk^{+}(\goth{f})/\kk)$,
$$\zeta'_{\goth{f}}([\goth{b}], 0) = -\frac{1}{e_{\goth{f}}} \log|u_{Stark}^{\sigma_{\goth{b}}}|^2$$
\item every complex embedding of $\,\kk^{+}(\goth{f})$ above a real embedding of $\,\kk$ sends $u_{Stark}$ to the unit circle.
\item $\kk^{+}(\goth{f})(u_{Stark}^{1/e_{\goth{f}}})$ is an abelian extension of $\,\kk$.
\end{itemize}
\end{famous}

\noindent There are ways to make this statement a bit more precise so that the unit $u_{Stark}$ is unique and hopefully the unit $\ula$ from Conjecture \ref{conjecture} will be related to the Stark unit by the formula:
$$u_{Stark}^{m(r,N)(N-\sigma_{\goth{a}})\cardinalshort{\mathcal{Z}_{\goth{f}}^1}} = \ula^{-e_{\goth{f}}}$$
In the cases where $\zeta_{\goth{f}}'([\goth{b}], 0) \neq 0$ for some class $[\goth{b}]$, our conjecture relates to Hilbert's twelfth problem, allowing for a construction of certain class fields using specific multivariate analytic functions.

\section{Computing the arithmetic $G_r$ functions}\label{sectioncomputing}

In order to understand how the choice of vectors $(h_{\rho})_{\rho \in \goth{S}_r}$ should be made, we must carry out an explicit description of the geometric aspects of the construction of the arithmetic $G_{r,\goth{f}, \goth{b}, \goth{a}}$ functions. Thus, we will identify key quantities in the construction, namely the parameters $t$ and $l$ defined in section \ref{geomsetup} and then show how we can control these quantities. We then try to minimize the quantities that have the greatest impact on the computation time and describe an efficient algorithm to compute the arithmetic products $\spvgrcomplete$ given in Definition \ref{defIrfba}.

\subsection{Geometric setup}\label{geomsetup}

Let $\kk$ be an ATR field of degree $d = r+2$ and let $\goth{f}$, $\goth{b}$, $\goth{a}$ be given as in section \ref{arithmeticsection}. Remember that we have fixed an ordering $\sigma_1, \dots, \sigma_{r}, \sigma_{r+1} = \sigma_{\cc}, \overline{\sigma_{r+1}}$ on the embeddings of $\kk$. Put $L = \goth{f}\goth{b}^{-1}$. Fix a system of fundamental units $(\eps_1, \dots, \eps_r)$ of $\opc{\goth{f}}$. Write $(u_1, \dots u_r) = [ \eps_{\rho(1)} | \dots | \eps_{\rho(r)}]$ for any permutation $\rho \in \goth{S}_r$. For convenience, we will always write $u_0 = 1$. Choose $h \in L$ an admissible vector (see Definition \ref{defadmissible}). Put $mh' = h$ where $m$ is a positive integer and $h'$ is primitive in $L$. This means that for all non-zero integer $n \in \zz$, $h'/n \in L$ implies that $n = \pm 1$.

Recall that a basis $B = (e_0, \dots, e_{r+1})$ of the $\qq$-vector space $\kk$ is called positive if $i\cdot\det((\sigma_j(e_{k-1}))_{1\leq j,k\leq d}) >0$. We may choose a positive basis $B_L = (e_0 = h', e_1, \dots, e_{r+1})$ of $L$ such that $u_jh' = \sum_{k = 0}^j c_{jk0}e_k$ where the coefficients $c_{jk0}$ are integers and $c_{jj0} > 0$. To evaluate geometric $G_r$ functions we fix the orientation form on the lattice $L$ given by $\det_{B_L}$. In the basis $B_L$ we identify the units $u_j$ with the matrix $(c_{jkl})_{k,l} \in \slnz{d}$ such that $u_je_l = \sum_{k = 0}^{r+1}c_{jkl}e_k$. Fix the dual basis $C = (f_0, \dots, f_{r+1})$ of $\Lambda = \homlz$ such that $f_j(e_k) = \delta_{jk}$ where $\delta_{jk}$ is the Kronecker symbol. In this basis the linear form $a$ defined in section \ref{arithmeticsection} is exactly $f_{r+1}$. For $1 \leq j \leq r$, the composition of $a$ with multiplication by $u_j$ is written $au_j = \sum_{l = 1}^{r+1} c_{j(r+1)l} f_l$. Then there is a unique primitive vector $\gamma \in L$ such that $\det_C(a, au_1, \dots, au_r, \cdot) = s\gamma = \pm s h'$ with $s$ a positive integer. 

Next, we identify the family $\abar = (a = au_0, au_1, \dots, au_r)$ with the matrix obtained by concatenation of the coefficients of the linear forms $au_i$ in the basis $C$:
$$\abar = \begin{pmatrix} au_0 \\ au_1 \\ \vdots \\ au_i \\ \vdots \\ au_r \end{pmatrix} = \begin{pmatrix} 0 & 0 & \dots & 0 & \dots & 1 \\ 0 & c_{1(r+1)1} & \dots & c_{1(r+1)k} &\dots& c_{1(r+1)(r+1)} \\ \vdots & \vdots & \vdots & \vdots & \vdots & \vdots \\ 0 & c_{i(r+1)1} & \dots & c_{i(r+1)k} & \dots & c_{i(r+1)(r+1)} \\\vdots & \vdots & \vdots & \vdots \\ 0 & c_{r(r+1)1} & \dots & c_{r(r+1)k} & \dots & c_{r(r+1)(r+1)} \end{pmatrix} = \begin{pmatrix} \textbf{0} & \mathcal{A} \end{pmatrix}$$
This matrix has a first column filled with zeroes and the submatrix $\mathcal{A}$ is a square matrix of size $r+1$. The elementary divisors of $\mathcal{A}$ associated to its Smith normal form are positive integers denoted $[A_r, \dots, A_1, A_0 = 1]$ such that $A_i | A_{i+1}$. Then we define 
\begin{align*}
\lambda(u_1, \dots, u_r; h) & = \lambda = \prod_{j = 1}^r c_{jj0} \\
s(u_1, \dots, u_r; h) & = s = \prod_{j=1}^r A_j \\
t(u_1, \dots, u_r; h) & = t = A_r
\end{align*}
so that 
\begin{equation}\label{defa}
\lambda a(u_1, \dots, u_r; h) = \lambda a = \det_{B_L}(h', u_1h', \dots, u_rh', \cdot)
\end{equation}
and $\det_{C}(a, au_1, \dots, au_r, \cdot) = \pm s h'$. Recalling the construction in Proposition \ref{propdefgeom}, we must now consider positive dual families $\alphabar = (\alpha_0, \dots, \alpha_r) \subset L$ to $\abar$. Recall that this means that for all $0 \leq j \leq r$:
$$au_j(\alpha_j) > 0,~~~~ au_k(\alpha_j) = 0, ~\forall\, k \neq j$$
with once again the convention that $u_0 = 1$. Consider the matrix $\mathcal{B} = (\prod_{i = 1}^r A_i^{-1}) \mathrm{com}(\mathcal{A})^{T}$ which corresponds to the rescaled comatrix in formula \refp{comatrixchoice}. Then the product $\mathcal{A}\cdot \mathcal{B} = t I_{r+1}$. Write $\mathcal{B} = (b_{ij})_{1\leq i,j \leq r+1}$. Then a possible choice for $\alphabar$ is given by $\alpha_j = \sum_{i = 1}^{r+1} b_{ij} e_i \in L$. This family satisfies for all $0 \leq j \leq r$:
\begin{equation}\label{defalpha}
au_j(\alpha_j) = t, ~~~~ au_k(\alpha_j) = 0, ~\forall\, k \neq j
\end{equation} 
We write this as the multiplication of two matrices $\abar \cdot \alphabar = tI_{r+1}$. We say that this choice of $\alphabar$ is a uniform positive dual family to $\abar$ because the value $au_j(\alpha_j) = t$ is independent of $j$. We argue that this choice is minimal amongst the uniform positive dual families to $\abar$ in $L$, as explained in the following lemma.

\begin{general}{Lemma}\label{tmin}
Let $\mathcal{B}' \in M_{r+1}(\zz)$ be a square matrix of size $r+1$ such that $\mathcal{A} \cdot \mathcal{B}' = n I_{r+1}$. Then $t$ divides $n$. Consequently, any uniform positive dual family $\alphabar'$ to $\abar$ satisfies $\abar \cdot \alphabar' = n I_{r+1}$ where $t$ divides $n$. Furthermore, the minimal uniform positive dual family $\alphabar$ to $\abar$ described above satisfies the following property: for all prime $p$ there exists an index $0 \leq j \leq r$ such that $\forall\, n' \in \zz, \,\alpha_j - n'\gamma \not\in pL$.
\end{general} 

\begin{proof}
Consider the matrices $U,V \in \glnz{r+1}$ such that $U\mathcal{A}V$ is the diagonal matrix in Smith normal form $[A_r = t, \dots, A_1, 1]$. Then
$$ nI_{r+1} = \left(U \mathcal{A} V\right) \cdot \left(V^{-1} \mathcal{B}' U^{-1}\right) = \begin{pmatrix} t & 0 & \dots & 0 \\ 0 & A_{r-1} & \dots & 0 \\ \vdots & \vdots & \vdots & \vdots \\ 0 & \dots & \dots & 1  \end{pmatrix} V^{-1} \mathcal{B}' U^{-1} $$
This shows that $t | n$. Consider now a uniform positive dual family $\alphabar' = (\alpha'_0, \dots, \alpha'_r)$ to $\abar$ in $L$ such that $\abar \cdot \alphabar' = nI_{r+1}$. Write the coordinates of the $\alpha'_j$ in the basis $B_L$ as $\alpha'_j = \sum_{i = 0}^{r+1} b'_{ij} e_i$. Then the matrix $\mathcal{B}' = (b_{ij}')_{1 \leq i,j \leq r+1}$ satisfies $\mathcal{A} \cdot \mathcal{B}' = n I_{r+1}$ and $t$ divides $n$.

Consider now the minimal choice $\alphabar$ given above associated to the matrix $\mathcal{B} = (\prod_{i = 1}^r A_i^{-1}) \mathrm{com}(\mathcal{A})^{T}$ such that $\abar \cdot \alphabar = tI_{r+1}$. The elementary divisors of $\mathcal{B}$ are given by the positive integers $[A_r, A_r/A_1, \dots, A_r/A_{r-1}, 1]$. The last elementary divisor is the $gcd$ of the coordinates of the $\alpha_j$ in the basis $B_L$. Therefore, for each prime $p$ there exists an integer $j$ and a coefficient $b_{ij}$ not divisible by $p$. For any $n' \in \zz$, we get $\alpha_j - n'\gamma = \pm n' e_0 + \sum_{i = 1}^{r+1} b_{ij} e_i \not\in pL$.
\end{proof}

From now on we will always fix $\alphabar = \alphabar(u_1, \dots, u_r; h)$ the minimal uniform positive dual family $\alphabar$ to $\abar$ given above. Define $M = M(u_1, \dots, u_r; h) = \zz h' \oplus (\oplus_{j = 0}^r \zz \alpha_j)$ in $L$ associated to this family. This sublattice has index $\det(\mathcal{B}) = t^{r+1}/s$ in $L$. Recalling the construction in section \ref{arithmeticsection} and formula \refp{defgeom} for the geometric $G_{r, \abar}$ function we may identify $F(\alphabar)/\zz\gamma \iso L/M$ so that the geometric $G_{r, \abar}$ function is a product of $t^{r+1}/s$ ordinary elliptic $G_r$ functions. In the worst case scenario $t = s$ and $G_{r,\abar}$ is a product of $t^r$ ordinary elliptic $G_r$ functions.

Let us now analyse $\alphabar$ more closely. It follows from the definition of the linear form $a$ in formula \refp{defa} and lemma \ref{tmin} that for $j \neq k$, $\alpha_j$ must satisfy
$$ \det_{B_L}(u_0h', u_1h', \dots, u_rh', u_k\alpha_j) = 0$$ with once again the convention that $u_0 = 1$, so that $\alpha_j/h' = \sum_{l = 0}^r d_{jkl}u_k^{-1}u_l$ with rational coefficients $d_{jkl}$. This quotient depends rather loosely upon $h'$ and this leads us to investigate the case $h = 1$ (which is not an admissible vector). In what follows, we will write with a $\sim$ the counterpart of the quantities defined in the construction above for $h=1$ and $\ok$ instead of the lattice $L$. There is a positive basis $B_{\kk} = [1, \te_1, \dots, \te_{r+1}]$ of $\ok$ such that $u_j = \sum_{k = 0}^j \tc_{jk0}\te_k$ with $\tc_{jj0} > 0$. Consider 
\begin{equation}\label{defta}
\tlambda \ta = \left(\prod_{j = 1}^r \tc_{jj0}\right) \ta = \det_{B_{\kk}}(1, u_1, \dots, u_r, \cdot)
\end{equation}
and put $\tabar = (\ta, \ta u_1, \dots, \ta u_r)$. Associate to $\tabar$ a square matrix $\tilde{\mathcal{A}}$ of size $r+1$ as above for $\abar$. Write as above the elementary divisors of this matrix $\tilde{\mathcal{A}}$ associated to $\tabar$ as $[\tA_r, \dots, \tA_1, \tA_0 = 1]$ such that $\tA_i | \tA_{i+1}$. Then put $\ts = \prod_{i = 1}^r \tA_i$, $\ttt = \tA_r$ and construct as above $\talphabar$ associated to the square matrix $\tilde{\mathcal{B}} = (\prod_{i = 1}^r \tA_i^{-1})\mathrm{com}(\tilde{\mathcal{A}})^{T}$ such that $\tabar \cdot \talphabar = \ttt I_{r+1}$. Recall that this means that for all $0 \leq j \leq r$:
\begin{equation}\label{deftalpha}
\ta u_j(\talpha_j) = \ttt, ~~~~ \ta u_k(\talpha_j) = 0, ~\forall\, k \neq j
\end{equation}

In the case $r = 1$, $\tlambda$ is the content of the smallest ring $\mathcal{O}$ containing the unit $\eps$ in the sense of [\hspace{1sp}\cite{Bhargava}, Definition 14] and $\ttt$  measures how far the $\zz$-module $\zz + \zz \eps$ is from being a ring. Therefore, we will call $\tlambda$ the content of the unit system $(u_1, \dots, u_r)$ and $\ttt$ the overflow of this unit system. We insist that these quantities only depend on the unit system $(u_1, \dots, u_r)$.
We relate this general construction to the one carried out for a particular admissible vector $h$ in $L$ as follows:

\begin{general}{Lemma}\label{fundam}
The family $\epsilon\talphabar h' = (\epsilon\talpha_0 h', \dots, \epsilon\talpha_r h')$ is a uniform positive dual family to $\abar$ in $L$ for some sign $\epsilon = \pm 1$. Therefore there exists a positive integer $l$ such that $\abar \cdot (\epsilon\talphabar h') = lt I_{r+1}$. The value of $l$ is given by the following equality which holds in $\zz$
\begin{equation}\label{formulafundam}
\epsilon \lambda l t = \frac{\norm{h'}}{\norm{L}} \tlambda \ttt
\end{equation}
\end{general}

\begin{proof}
\noindent We use formulae \refp{defa} relating $\abar$ to the form $\det_{B_{\kk}}$ and \refp{defta} relating $\tabar$ to the form $\det_{B_L}$ to evaluate $\lambda a u_k \talpha_j h'$.
$$\lambda a u_k \talpha_j h' = \det_{B_L}(h', u_1 h', \dots, u_r h', u_k\talpha_j h')$$
Using standard linear algebra and the definition of the norm of the fractional ideal $L$ we may rewrite this equality using $\det_{B_{\kk}}$  instead of $\det_{B_L}$ as
\begin{align*}
\lambda a u_k \talpha_j h' & = \det_{B_L}(B_{\kk}) \det_{B_{\kk}}(h', u_1h', \dots, u_rh', u_k \talpha_jh')\\
\lambda a u_k \talpha_j h' &=\frac{1}{\norm{L}}\det_{B_{\kk}}(h', u_1h', \dots, u_rh', u_k \talpha_jh')
\end{align*}
Then, using the definition of the norm $\norm{h'}$ of $h'$ we get
\begin{align*}
\lambda a u_k \talpha_j h' & = \frac{\norm{h'}}{\norm{L}} \det_{B_{\kk}}(1, u_1, \dots, u_k \talpha_j)\\
\lambda a u_k \talpha_j h' & = \frac{\norm{h'}}{\norm{L}} \tlambda \ta u_k \talpha_j
\end{align*}
We now use formula \refp{deftalpha} to conclude that for $k \neq j$, $\lambda a u_k \talpha_j h' = 0$ and for $k =j$: 
$$\lambda a u_j \talpha_j h' = \frac{\norm{h'}}{\norm{L}} \tlambda \ttt \in \zz$$
$$a u_j \talpha_j h' = \frac{\norm{h'}}{\norm{L}} \frac{\tlambda}{\lambda} \ttt$$
The lattice $L$ is a fractional ideal so that $\talphabar h' \subset L$. The linear forms $au_k$ take integral values on $L$ so the number $n = \frac{\norm{h'}}{\norm{L}} \frac{\tlambda}{\lambda} \ttt$ is a non-zero integer. Write $\epsilon = \pm 1$ for the sign of $n$. Then $\epsilon\talphabar h'$ is a uniform positive dual family to $\abar$ in $L$ such that $\abar \cdot (\epsilon\talphabar h') = |n|I_{r+1}$. Lemma \ref{tmin} gives $n = \epsilon lt$ for a positive integer $l$ and 
$$\epsilon \lambda l t  = \lambda n = \lambda au_j \talpha_j h' =  \frac{\norm{h'}}{\norm{L}} \tlambda \ttt $$
\end{proof}

\noindent We now give a precise statement on the relation between $\alphabar$ and $\talphabar$.

\begin{general}{Proposition}\label{propl}
Let $h$ be an admissible vector in $L$ and consider a unit system $u_1, \dots, u_r$ for $\opc{\goth{f}}$. Let $\alphabar$ and $\talphabar$ be as defined in \refp{defalpha} and \refp{deftalpha}. Write $\epsilon$ for the sign in formula \refp{formulafundam}.
Then for all $0 \leq j \leq r$, there exists an integer $m_j \in \zz$ such that $\epsilon l\alpha_j = \talpha_jh' + m_jh'$.
\end{general}

\begin{proof}
Using lemma \ref{lemmachoice} for $\abar$ and the two positive dual families $\alphabar$, $\epsilon\talphabar h'$ to $\abar$ in $L$ we get $\epsilon l\alpha_j = \talpha_j h' + m_j h'/n_j$ with $m_j/n_j$ a rational number in irreducible form. Write $p_jm_j + q_jn_j = 1$. Then 
$$p_j(\epsilon l\alpha_j - \talpha_j h') + q_j h' = h'/n_j \in L$$
Because $h'$ is a primitive vector in $L$, we conclude that $n_j  = 1$.
\end{proof}

This proposition shows that the choice of $h$ only affects the value of $l$ and the values of the $m_j$ and not the main content of $\alphabar/h'$. Indeed, the special value $\spvgr$ we aim to compute is a product of $t^{r+1}/s$ ordinary elliptic $G_r$ functions with parameters $\alpha_0/h', \dots, \alpha_r/h'$, which we may write as $\epsilon(\talpha_0+m_0)/l, \dots, \epsilon(\talpha_r+m_r)/l$. For instance in the cubic case, $\Gamma^{\pm}_{\goth{f}, \goth{b}, \goth{a}}$ is a product of $t$ ordinary elliptic Gamma functions with parameters
$$\tau = \pm\frac{\eps + n_0}{l\tlambda} = \frac{\tau_0}{l}, \sigma = \pm \frac{\eps^{-1} + n_1}{l\tlambda} = \frac{\sigma_0}{l}$$ for some integers $n_0, n_1$. Therefore, the integer $l$ can be thought of as a level and will be called the level of the value $\spvgr$.

\subsection{Computing the $G_r$ functions}

\subsubsection{Computing the ordinary elliptic $G_r$ functions}\label{ordinarycomputations}

To compute the $G_r$ functions we make use of [\hspace{1sp}\cite{Nishizawa}, Proposition 3.6] which we write as
\begin{equation}\label{sumcomputations}
G_r(z, \taubar) = \exp\left(-\sum_{j \geq 1} \frac{1}{j} \frac{q_0^j\dots q_r^jx^{-j} +(-1)^r x^{j}}{\prod_{k = 0}^r (1 - q_k^j)}\right)
\end{equation}
This formula is only valid for $\taubar \in \hh^{r+1}$ and $0 < \Im(z) < \sum_{k = 0}^r \Im(\tau_k)$. We call this domain the center strip. We now only need to make sure that we can use the properties of the $G_r$ functions to reach this domain. The first step in doing so, starting from $\taubar \in (\cc - \rr)^{r+1}$ is to use the property 
$$G_r(z, \tau_0, \dots, \tau_{k-1}, -\tau_k, \tau_{k+1}, \dots, \tau_r) = G_r(z+\tau_k, \taubar)^{-1}$$
to bring back all arguments in the upper half-plane. Then, to bring back $\Im(z)$ in the center strip, use recursively the property
$$G_r(z+\tau_k, \taubar) = G_{r-1}(z, \taubar^{-}[k]) G_r(z, \taubar)$$  
Doing so will require the computation of lower degree functions, that is until we reach the case $r = 0$. Lastly, to compute $G_0 = \theta$ one may use Jacobi's triple product formula and compute directly $\theta$ using this fast-converging sum:
$$\theta(z, \tau) = \frac{q^{1/24}}{\eta(\tau)} \sum_{n \in \zz} x^n (-1)^n q^{n(n-1)/2}$$
where $\eta(\tau)$ is the Dedekind $\eta$ function. Let us say a word about the complexity of such computations. Fix $\tau_0, \dots, \tau_r$ in the upper half-plane and suppose we want to compute the level $l$ function $G_r(z, \taubar/l)$ with precision $\delta > 0$. Then the number of terms we need to compute in the sum \refp{sumcomputations} is $O(l |\log(\delta)|)$ where the constant explicitly depends on $\taubar$ and $z$.

The above algorithm will be useful to compute the value $G_{r, \goth{f}, \goth{b}, \goth{a}}^{\mu}(u_1, \dots, u_r ; h)$ when for all $0 \leq j \leq r$, $\sigma_{\cc}(\talpha_j) \in \cc - \rr$. This is always the case except possibly when $\kk$ is a field of even degree $d = 2d'$ and contains a real subfield $\mathbb{F}$ of degree $\leq d'$. In this case, because the family $1, \talpha_0, \dots, \talpha_r$ is free over $\qq$, at most $d'-1$ elements amongst the $\talpha_j$ may belong to $\mathbb{F}$. We will avoid this case for now under the ($\ast$) conditions below and discuss one example of such a case in section \ref{quarticreal}.

\subsubsection{Computing the arithmetic $G_{r, \goth{f}, \goth{b}, \goth{a}}^{\mu}(u_1, \dots, u_r ; h)$ functions}

We analyse Definition \ref{defgrarithmetic} for $G_{r, \goth{f}, \goth{b}, \goth{a}}^{\mu}(u_1, \dots, u_r ; h)$ using the geometric setup in section \ref{geomsetup} for the choice of sign $\mu = \pm 1$. Recall Definition \ref{defgrarithmetic} for the arithmetic function:
$$G_{r, \goth{f}, \goth{b}, \goth{a}}^{\mu}(u_1, \dots, u_r ; h, \sigma_{\cc}) = \frac{G_{r, (\mu a, \mu au_1, \dots, \mu au_r)}(\frac{h}{q}, \sigma_{\cc},\mathrm{Hom}_{\zz}(\goth{f}\goth{b}^{-1}, \zz))^N}{G_{r, (\mu a, \mu au_1, \dots, \mu au_r)}(\frac{h}{q}, \sigma_{\cc}, \mathrm{Hom}_{\zz}(\goth{f}(\goth{b}\goth{a})^{-1},\zz))}$$
We drop once again the $\sigma_{\cc}$ and we write using formula \refp{defgeom}:
$$G_{r, \goth{f}, \goth{b}, \goth{a}}^{\mu}(u_1, \dots, u_r ; h) = \prod_{\delta \in F(\alphabar)/\zz\gamma} \frac{G_r\left(\frac{h+q\mu \delta}{q\mu^{r+1}\gamma}, \frac{\mu\alpha_0}{\mu^{r+1}\gamma}, \dots, \frac{\mu\alpha_r}{\mu^{r+1}\gamma} \right)^N}{G_r\left(\frac{N(h+ q\mu\delta)}{q\mu^{r+1}\gamma}, \frac{N\mu\alpha_0}{\mu^{r+1}\gamma}, \dots, \frac{N\mu\alpha_r}{\mu^{r+1}\gamma} \right)}$$
where $\alphabar$ and $\gamma$ are defined in section \ref{geomsetup}. We use formula \refp{inversiontotale} as well as the fact that $F(\alphabar)/\zz\gamma \iso L/M$ where $M = M(u_1, \dots, u_r; h)$ is defined in section \ref{geomsetup} to write 
$$G_{r, \goth{f}, \goth{b}, \goth{a}}^{\mu}(u_1, \dots, u_r ; h)^{\mu^r} = \prod_{\delta \in L/M} \frac{G_r\left(\frac{\mu h+q\delta}{q\gamma}, \frac{\alpha_0}{\gamma}, \dots, \frac{\alpha_r}{\gamma} \right)^N}{G_r\left(\frac{N(\mu h+q\delta)}{q\gamma}, \frac{N\alpha_0}{\gamma}, \dots, \frac{N\alpha_r}{\gamma} \right)}$$
\begin{equation}\label{truecomputations}
G_{r, \goth{f}, \goth{b}, \goth{a}}^{\mu}(u_1, \dots, u_r ; h)^{\mu^r} = \prod_{\delta \in L/M} \frac{G_r\left(\frac{\mu h+q\delta}{q\gamma}, \epsilon\frac{\talpha_0+m_0}{l}, \dots, \epsilon\frac{\talpha_r+m_r}{l} \right)^N}{G_r\left(\frac{N(\mu h+q\delta)}{q\gamma}, \epsilon\frac{N(\talpha_0+m_0)}{l}, \dots, \epsilon\frac{N(\talpha_r+m_r)}{l} \right)}
\end{equation}
where the last expression is obtained using Proposition \ref{propl}. This special value is a product of $t^{r+1}/s = \#|L/M|$ ordinary elliptic $G_r$ functions. In the worst case scenario we get $t = s$ and $t^{r+1}/s = t^r$. This means that the computation time of this value is $O(l|\log\delta|t^r)$ with a constant that depends on $u_1, \dots, u_r$. The parameters $l$ and $t$ which depend on $h$ for fixed $u_1, \dots, u_r$ are therefore crucial in our computations. For instance, if we want to compute this special value with a precision of $1000$ digits (i.e. $|\log\delta| = 1000$) where $r = 2$, $l = 15$ and $t = 100$, then the computation requires the computation of $150.10^6$ terms appearing in the sum \refp{sumcomputations}. With a personal computer it is generally an overwhelming computation. In the next section, we describe how we can choose the vector $h$ to obtain efficient computations. Namely, because the parameter $t$ has the biggest impact on computation time, we will show that under some assumptions we may choose a vector $h$ such that $t = 1$ at the cost of a moderate elevation of the level $l$.

\subsection{The choice of vectors $(h_{\rho})_{\rho \in \goth{S}_r}$}\label{choiceh}

In view of the above estimation on the computation time for the special value $\spvgrrho$ for a permutation $\rho \in \goth{S}_r$ we must be careful in choosing the vector $h_{\rho}$ to perform efficient computations. Also, because we want to build a product of these values for $\rho \in \goth{S}_r$, we must find a way to unify our choices for each permutation $\rho$. To this end, we will fix the unit system $\eps_1, \dots, \eps_r$ and the permutation $\rho$. Write once again $(u_1, \dots, u_r) = [\eps_{\rho(1)}|\dots| \eps_{\rho(r)}] = (\eps_{\rho(1)}, \dots, \prod_{j = 1}^r\eps_{\rho(j)})$. Our aim is to reduce at all costs the value of $t = t(u_1, \dots, u_r; h)$ by choosing $h$ carefully. Formula \refp{formulafundam}
$$\pm \lambda l t = \frac{\norm{h'}}{\norm{L}} \tlambda \ttt$$
will be fundamental to understand how to minimize the value of $t$. The unit system $(u_1, \dots, u_r)$ is fixed and we therefore have to deal with the fixed constants $\tlambda, \ttt$ defined in section \ref{geomsetup}. In what follows, a goal will be to understand how the prime factors of $\tlambda$ and $\ttt$ transfer to the left-hand side of this equation. From our computations it seems that the value of $\tlambda$ is a hard limit for the value of $t$, which means that for all admissible vector $h$, $\tlambda$ divides $t$. We also believe that this value of $\tlambda$ is exactly the minimal achievable value for $t$. We are far from being able to prove this statement, but it will be proven in the case $\tlambda = 1$ under the conditions ($\ast$) and ($\ast\ast$) below. The rest of this section contains the study of the $p$-adic valuations of the quantities appearing in formula \refp{formulafundam} first in the cases where $p$ divides $\ttt$ and then in the case where $p$ does not divide $\ttt$. Under some conditions we show that in the first case the $p$-part of $h$ (as in definition \ref{ppart} below) is best chosen as specific ideals which we call \textit{target ideals} (see definition \ref{deftarget} and lemma \ref{lemmatarget} below). In the second case, it will be shown how a vector $h$ may have a non-trivial $p$-part while not affecting the value of $t$ using \textit{helper ideals} (see lemma \ref{lemmahelper} below). In the rest of this section, we will discuss a uniform choice for the vectors $h_{\rho}, \rho \in \goth{S}_r$ under the following conditions, referred to as the ($\ast$) conditions: suppose that the ideal $\goth{f} = \goth{l}^k$ is a power of a degree one prime ideal $\goth{l}$ above a prime number $\ell$ and that the choice of fundamental units $\eps_1, \dots, \eps_r$ for $\opc{\goth{f}}$ is such that for all $\rho \in \goth{S}_r$:
\begin{itemize} \label{conditions}
\item $\tilde{\lambda}_{\rho} = 1$ 
\item $\ttt_{\rho}$ and $qN\mathcal{N}(\goth{b})$ are relatively coprime \hfill ($\ast$)
\item for all $0 \leq j \leq r$, $\sigma_{\cc}(\talpha_{\rho, j}) \not\in \rr$
 \end{itemize}

\noindent where the subscript $\rho$ indicates that the quantities are related to the unit system $[\eps_{\rho(1)}| \dots | \eps_{\rho(r)}]$. Note that the value of $\ttt_{\rho}$ does not depend on $\goth{b}$ or $\goth{a}$ which means that in the second line, we only require $\ttt_{\rho}$ to be coprime to $q$ and we decide to choose $\goth{b}$ and $\goth{a}$ conveniently afterwards, so that this part is not restrictive. The last assumption is only relevant if $\kk$ is of even degree $2d' = d$ and contains a real subfield $\mathbb{F}$ of degree $\leq d'$, in which case at most $d'-1$ of the $\talpha_{\rho, j}$ may belong to $\mathbb{F}$. This particular situation is discussed in section \ref{quarticreal}.
We will review the other assumptions when we discuss ($\ast\ast$) below.

\subsubsection{Target ideals}

While choosing the vector $h = h_{\rho}$, we must be very careful of the value of the overflow $\ttt = \ttt_{\rho}$. Indeed, for most admissible vectors $h$, the constant $\ttt$ depending only on $u_1, \dots, u_r$ will divide $t = t(u_1, \dots, u_r ; h)$ and the computations will require $t^r$ evaluations of ordinary $G_r$ functions, which can become quickly overwhelming. We now introduce the notion of target ideal to understand how to choose $h$ such that $t = 1$. Explicitly, we explain that there is a tradeoff with an increase of the level $l$ as we diminish the value of $t$. In the rest of this section we will make extensive use of the following definitions:

\begin{general}{Definition}\label{ppart}
Let $x \in \kk$ and $I$ be an ideal such that $\norm{I}$ is a power of the prime number $p$. We say that the $p$-part of $x \in \kk$ is the ideal $I$ if for all prime ideal $\goth{P}$ above $p$, $v_{\goth{P}}(x) = v_{\goth{P}}(I)$. 
\end{general}

\begin{general}{Definition}
Let $I$ be a non zero integral ideal in $\ok$ such that $I \cap \zz = n\zz$. Then we define the complementary ideal $I^c$ of $I$ to be the integral ideal $n\ok/I$ so that $II^c = n\ok$. We call a non zero integral ideal primitive if $I$ is not contained in any ideal of the form $m\ok$ with $m \in \zz, m \geq 2$ or equivalently if $(I^c)^c = I$.
\end{general}

Before we define the notion of target ideal, we prove a lemma which we will need to control the $p$-adic valuation of the integer $\lambda$ in the construction of section \ref{geomsetup}.

\begin{general}{Lemma}\label{lemmalambda}
Take $h = mh'$ an admissible element in $L$. Suppose that there is a linear combination 
$$\left(u_j + \sum_{i = 0}^{j-1} n_i u_i\right)h' \in p^kL$$
where $j \geq 1$ and the coefficients $n_i$ are integers. Then $p^k$ divides the integer $c_{jj0}$ and therefore $p^k$ divides the integer $\lambda$.
\end{general}

\begin{proof}
Recall that the basis $B_L = (e_0, \dots, e_{r+1})$ defined at the beginning of section \ref{geomsetup} is such that $e_0 = h'$ and for all $1 \leq i \leq r$, $u_ih' = \sum_{k = 0}^i c_{ik0}e_k$ where the $c_{ik0}$ are integers and $c_{ii0} > 0$. This means that 
$$\left(u_j + \sum_{i = 0}^{j-1} n_i u_i\right)h' = c_{jj0}e_j + \left(\sum_{k = 0}^{j-1} c_{jk0} e_k\right) + \sum_{i = 0}^{j-1} \left(\sum_{k = 0}^{i} n_ic_{ik0} e_k\right) \in p^kL$$
so the coordinate of this element on the vector $e_j$ is exactly $c_{jj0}$ and $p^k$ divides $c_{jj0}$. Recalling that $\lambda = \prod_{j = 1}^r c_{jj0}$, we get $p^k$ divides $\lambda$.
\end{proof}

We now define the notion of target ideal to deal with the $p$-adic valuation of the overflow $\ttt$.

\begin{general}{Definition}\label{deftarget}
Assume $p$ is a prime number such that $v_p(\ttt) > 0$. Let $I$ be a non zero primitive integral ideal of $\ok$ such that $II^c = p^{v_p(\ttt)}\ok$. We say that $I$ is a target ideal associated to the unit system $u_0 = 1, u_1, \dots, u_r$ and to the prime $p$ if $\norm{I} = p^{v_p(\ts)}$ and the two following conditions are satisfied:
\begin{itemize}
\item there are distinct integers $1 \leq m_1 < \dots < m_{k_0} \leq r$ and there are linear combinations $u_{m_k} + \sum_{j= 0}^{m_k-1}n_{k,j} u_j \in J_k$ with $n_{k,j} \in \zz$ and $J_k$ a divisor of $I^c$ such that $IJ_k \subset p^{w_k}\ok$ for some valuations $w_k \in \zz_{>0}$ satisfying $\sum_{k = 1}^{k_0} w_k = v_p(\ts)$.
\item for all $0 \leq j \leq r$, $\talpha_j$ is congruent to an integer $n'_j$ modulo $I^c$
\end{itemize}
\end{general}

\noindent The importance of target ideals is expressed by the following lemma.

\begin{general}{Lemma}\label{lemmatarget}
Suppose $(\ast)$. Let $p$ be a prime number such that $v_p(\ttt) > 0$ and suppose that $I$ is a target ideal for the unit system $u_0, \dots, u_r$. Take $h = mh'$ an admissible vector such that $h'$ has $p$-part exactly $I$. Then:
$$v_p(\lambda) = v_p(\ts),~~ v_p(l) = v_p(\ttt),~~v_p(t) = 0$$
\end{general}

\begin{proof}
The ideal $I$ is a target ideal therefore there are distinct integers $1 \leq m_1 < \dots < m_{k_0} \leq r$ and linear combinations $(u_{m_k} + \sum_{j= 0}^{m_k-1}n_{k,j} u_j)h' \in IJ_kL \subset p^{w_k}L$ with $n_{k,j} \in \zz$. using lemma \ref{lemmalambda} we obtain that $p^{w_k}$ divides $c_{m_k,m_k,0}$ and therefore $p^{v_p(\ts)} = \prod_{k = 1}^{k_0} p^{w_k}$ divides $\lambda$. 

On the other hand, from lemma $\ref{tmin}$ there is an index $j$ such that $p$ does not divide $\alpha_j - nh'$ for all $n \in \zz$ and from Proposition \ref{propl} there is an integer $m_j$ such that $\epsilon l\alpha_j = \talpha_jh' + m_jh'$ for some sign $\epsilon = \pm 1$. The ideal $I$ is a target ideal, therefore there is an integer $n'_j$ such that $\talpha_j - n'_j \in I^c$. This implies that $\epsilon l\alpha_j - m_jh' - n'_jh' = (\talpha_j -n'_j)h' \in p^{v_p(\ttt)}L$. Suppose for now that $p$ does not divide $l$. Then $\exists\, n \in \zz$ such that $\epsilon ln - m_j - n'_j \in p\zz$. Therefore $\epsilon l(\alpha_j-nh') = (\talpha_j -n'_j)h' + (m_j + n'_j - \epsilon ln)h' \in pL$ but this is a contradiction. Therefore $p$ must divide the integer $l$. We then prove by induction that $p^2, \dots, p^{v_p(\ttt)}$ divide $l$. Indeed, if $p^k | l$ with $k < v_p(\ttt)$ then $\epsilon l\alpha_j - (m_j+n'_j)h' \in p^{v_p(\ttt)}L$ which means that $m_j + n'_j \in p^{k}\zz$. Suppose that $l \not\in p^{k+1}\zz$. Then there is an integer $n \in \zz$ such that $\epsilon ln - m_j - n'_j \in p^{k+1}\zz$ and then $\epsilon l(\alpha_j -nh') \in p^{k+1}L$ which yields again a contradiction. Now, using formula \refp{formulafundam} one finds
$$v_p(\lambda) + v_p(l) + v_p(t) = v_p(\norm{h'}) - v_p(\norm{L}) + v_p(\tlambda) + v_p(\ttt)$$
$$v_p(\lambda) + v_p(l) + v_p(t) = v_p(\ts) - 0 + 0 + v_p(\ttt),~~v_p(\lambda) \geq v_p(\ts), v_p(l) \geq v_p(\ttt)$$
And this leads to $v_p(\lambda) = v_p(\ts)$, $v_p(l) = v_p(\ttt)$ and $v_p(t) = 0$.
\end{proof}

This allows us to strengthen the conditions under which we might understand special values of arithmetic $G_r$ functions. We add to the conditions ($\ast$)  the following conditions ($\ast\ast$):
\begin{itemize}
\item for all $\rho \in \goth{S}_r$, for all $p | \ttt_{\rho}$, there exists a target ideal $I_{\rho, p}$  for the unit system $[\eps_{\rho(1)}|\dots| \eps_{\rho(r)}]$ and for the prime $p$ \hfill ($\ast\ast$)
\item for all $\rho \in \goth{S}_r$, there exists a positive integer $k_{\rho}$ coprime to $q$ such that the class of $k_{\rho}\prod_{p|\ttt} I_{\rho, p}$ in the wide class group $Cl(\goth{f})$ is independent of $\rho$
\end{itemize} 
We discuss briefly what we expect from the conditions ($\ast$) and ($\ast\ast$). In regards to our computations for now, the most restrictive conditions would be the content condition $\tlambda_{\rho} = 1$, as well as the conditions $\goth{f} = \goth{l}^k$ and $\ttt_{\rho}$ coprime to $q$ as together they rule out many possible ideals $\goth{f}$. We expect the target ideal condition $(\ast\ast)$ to be true generally under the conditions ($\ast$) and we expect the target ideals to be uniquely determined. Finally, the condition on the classes represented by $\prod_{p|\ttt} I_{\rho, p}$ is expected to be often true.

\subsubsection{Helper ideals}
Target ideals helped us understand how to choose the $p$-part of $h$ when $p | \ttt$. It is not clear however how to choose an admissible vector with this specific $p$-part. Especially when the class number of $\kk$ is not $1$, it is most likely that the $p$-part of an admissible vector $h$ for some $p$ not dividing $\ttt$ will be non-trivial, in which case it must be carefully chosen. To understand the best case scenario, we introduce the notion of \textit{helper ideals}.

\begin{general}{Definition}
A primitive ideal $I$ is a helper ideal if the norm $\norm{I^c}$ of its complementary ideal $I^c$ is squarefree.
\end{general}

These ideals are interesting because they do not impact the value of $t$ as explained by the following lemma.

\begin{general}{Lemma}\label{lemmahelper}
Let $p$ be a prime number not dividing $\ttt qN\norm{\goth{b}}$. Suppose that $h = mh'$ is an admissible vector such that the ideal generated by $h'$ has $p$-part a helper ideal dividing $p\ok$. This means that there is a prime ideal $\goth{P}$ of norm $p$ such that $h'\goth{P} \in pL$ and $h' \not\in pL$. Then under the assumptions  $(\ast)$, $v_p(l) = 1$ and $v_p(t) = 0$.
\end{general}

\begin{proof} 
We use once again formula \refp{formulafundam}. Under the assumptions ($\ast$), we have $\tlambda = 1$. Here, by assumption we also have $v_p(\ttt) = 0$, $v_p(\norm{h'}/\norm{L}) = d-1$. Note that $\goth{P}$ is a degree one prime so all elements in $\ok$ are congruent to some integer mod $\goth{P}$ via the isomorphism $\ok/\goth{P} \simeq \zz/p\zz$. Hence, there are integers $r_j$ such that $u_j - r_j \in \goth{P}$ for all $1 \leq j \leq r$ and therefore $(u_j - r_j)h' \in pL$. Using lemma \ref{lemmalambda} successively for all indices $1 \leq j \leq r$ with integers $n_0 = -r_j$ and $n_i = 0$ for $1 \leq i \leq j-1$ we get that $p$ divides $c_{jj0}$. Therefore, $v_p(\lambda) \geq r = d-2$ and using formula \refp{formulafundam} we get $v_p(lt) \leq 1$. We now want to prove that $v_p(l) = 1$ so that $v_p(t) = 0$. This is essentially the same proof as for target ideals. There exists $0 \leq j \leq r$ such that $p$ does not divide $\alpha_j -nh'$ for all $n \in \zz$ (in fact this is true for all $j$ here). Take $m_j$ an integer such that $\talpha_j - m_j \in \goth{P}$. There are an integer $k_j$ and a sign $\epsilon = \pm 1$ such that $\epsilon l \alpha_j = (\talpha_j - m_j)h' + k_jh'$. Suppose that $p$ does not divide $l$. Then there is an integer $n$ such that $\epsilon ln - k_j \in p\zz$ and $\epsilon l(\alpha_j -nh') = (\talpha_j - m_j)h' - (\epsilon ln - k_j)h' \in pL$. This is once again a contradiction and we must conclude that $p$ divides $l$ and $p$ does not divide $t$.
\end{proof}

\subsubsection{Using target and helper ideals to choose $(h_\rho)_{\rho \in \goth{S}_r}$}

We now use both target ideals and helper ideals to describe how we choose the vectors $h_{\rho}$ in Conjecture \ref{conjecture}. 

\begin{general}{Proposition}\label{propfinal}
Suppose that $\eps_1, \dots, \eps_r$ is a set of fundamental units for $\opc{\goth{f}}$ satisfying $(\ast)$ and $(\ast\ast)$. For all $\rho \in \goth{S}_r$, write $p_{\rho, j}$ for the prime numbers dividing $\ttt_{\rho}$ for $1 \leq j \leq j_{\rho}$ and $I_{\rho, j}$ for any corresponding target ideal in $(\ast\ast)$. Put $D = q \prod_{\rho}\prod_{j =1}^{j_{\rho}} p_{\rho, j}$. Choose for the smoothing ideal $\goth{a}$ a degree one prime ideal of norm $N$ coprime to $D$. Choose a representative integral ideal $\goth{b}$ for a class in $Cl^{+}(\goth{f})$ coprime to $DN$. Then there is an ideal $\goth{I}$ independent of $\rho$ which is either $\ok$ or a helper ideal coprime to $DN\norm{\goth{b}}$ and a positive integer $m_{\rho}$ coprime to $q$ such that for all $\rho \in \goth{S}_r$ the ideal
$$ q\frac{N}{\goth{a}\goth{b}} m_{\rho}\goth{I}\prod_{j = 1}^{j_{\rho}} I_{\rho, j}$$ 
is principal, with an admissible generator $h_{\rho}$. For all orientations $\bars{\mu}, \bars{\nu}$, the value $\spvgrcomplete$ is a product of $r!$ smoothed ordinary elliptic $G_r$ functions.
\end{general}

\begin{proof}
We want to pick an admissible generator of the ideal $q\frac{N}{\goth{a}\goth{b}} \prod_{j = 1}^{j_{\rho}} I_{\rho, j}$ but this ideal may not be principal nor possess an admissible generator. We use helper ideals to deal with those as follows. By Cebotarev's density theorem (see [\hspace{1sp}\cite{Cox}, Theorem 8.17]), there are infinitely many degree one prime ideals in each class of the wide ideal class group $Cl(\goth{f})$. For any prime integer $p$ coprime to $q$ there exists an integer $m(p)$ coprime to $q$ such that $pm(p) \equiv 1 \mod q$. Conditions $(\ast\ast)$ states that there exist integers $k_{\rho}$ such that all ideals $k_{\rho}\frac{N}{\goth{a}\goth{b}} \prod_{j = 1}^{j_{\rho}} I_{\rho, j}$ belong to the same class $\goth{c}$ in $Cl(\goth{f})$. Consider a degree one prime $\goth{P}$ of norm $p$ coprime to $DN\norm{\goth{b}}$ in the class $\goth{c}$. Then we may choose $m_{\rho} = k_{\rho}m(p)$ and $\goth{I} = p/\goth{P}$. One notices that if the class $\goth{c}$ is already the trivial class, then one may take $m_{\rho} = k_{\rho}$ and $\goth{I} = \ok$. In either case, we may find a generator $g_{\rho}$ for the ideal $\frac{N}{\goth{a}\goth{b}} m_{\rho}\goth{I}\prod_{j = 1}^{j_{\rho}} I_{\rho, j}$ which is congruent to $1 \mod \goth{f}\goth{b}^{-1}$. We then obtain an admissible vector $h_{\rho} = q g_{\rho}$ because $h_{\rho}/q \equiv 1 \mathrm{~mod~} L$ and $h_{\rho}/N \in \goth{a}^{-1}L - L$ is a generator of the cyclic group $\goth{a}^{-1}L/L$ of prime order $N$.

Using lemma \ref{lemmatarget} we get $v_{p_{\rho,j}}(\lambda_{\rho}) = v_{p_{\rho,j}}(\ts_{\rho}), v_{p_{\rho,j}}(l_{\rho}) = v_{p_{\rho,j}}(\ttt_{\rho})$ and $v_{p_{\rho,j}}(t_{\rho}) = 0$ for all $\rho, j$. Using lemma \ref{lemmahelper} for all prime divisors of $n$ such that $\norm{\goth{I}} = n^{d-1}$ we get $l_{\rho} = nl'_{\rho}$ with $gcd(n, l'_{\rho}) = 1$ and $gcd(n, t_{\rho}) = 1$. Now we have to understand what happens for the $\ell$-adic valuations in formula \refp{formulafundam} where $\norm{\goth{f}} = q = \ell^k$. In this setting, we get $v_{\ell}(\lambda_{\rho}) = rk = (d-2)k$, $v_{\ell}(l_{\rho}) = k$, and $v_{\ell}(\norm{h'_{\rho}}/\norm{L}) = kd-k$ such that using formula \refp{formulafundam} we get $v_{\ell}(t_{\rho}) = 0$. The value of $\norm{\goth{b}}$ is completely invisible in formula \refp{formulafundam} because 
$$\norm{h'_{\rho}}/\norm{L} = q^{d-1}N^{d-1}\prod_{j = 1}^{j_{\rho}} p_{\rho, j}^{v_{p_{\rho,j}}(\ts_{\rho})} n^{d-1}$$
Lastly, $N/\goth{a}$ behaves in fact as a helper ideal for which we have $v_N(l_{\rho}) = 1$ and $v_N(t_{\rho}) = 0$. We use formula \refp{formulafundam} and we obtain the values 
$$l_{\rho} = qNn\prod_{j = 1}^{j_{\rho}} p_{\rho, j}^{v_{p_{\rho, j}}(\ttt_{\rho})}$$
for the levels in the computations of the values $\spvgrrho$, as well as $t_{\rho} = 1$, which implies that the value $\spvgrcomplete$ is a product of $r!$ smoothed ordinary elliptic $G_r$ functions.
\end{proof}

\noindent \textbf{Remark:} Our choices for the vectors $h_{\rho}$ are made to satisfy two important conditions. The main focus is to produce values for $\spvgrcomplete$ with arithmetical properties. It is plainly not true that all such products of $G_r$ functions yield algebraic numbers. A key point in computing interesting $\spvgrcomplete$ is to make sure that the vectors $h_{\rho}$ share some common property as $\rho$ varies in $\goth{S}_r$. We believe that this property has to do with the lattices $M_{\rho}$ defined in section \ref{geomsetup} being somewhat ``compatible'' if not equal. The condition $t_{\rho} = 1$ for all $\rho$ guarantees that $M_{\rho} = L$ for all $\rho$, which can be seen as the strongest form of this common property. On the other hand, this condition is useful to perform quick computations.

\subsubsection{Invariance under the choice of vectors $(h_{\rho})_{\rho \in \goth{S}_r}$}\label{sectionchoiceh}

In this section we define the set $\zfone$ appearing in Conjecture \ref{conjecture} and explain why it is necessary to compute an average on this set. We first begin by two simple lemmas.

\begin{general}{Lemma}\label{lemmaintvar}
Take $h$ an admissible vector for $\goth{f}, \goth{b}, \goth{a}$. Fix a system $u_1, \dots, u_r$ of linearly independent units. Consider $n$ a positive integer congruent to $1$ mod $q$. Then
$$ G_{r, \goth{f}, \goth{b}, \goth{a}}^{\mu}(u_1, \dots, u_r; nh) = G_{r, \goth{f}, \goth{b}, \goth{a}}^{\mu}(u_1, \dots, u_r; h)$$
\end{general}

\begin{proof}
Take $\alphabar, \gamma, M$ as defined in section \ref{geomsetup} for the orientation $\mu$. Write $h = m'\gamma$ for some integer $m' \in \zz -\{0\}$. Then
$$G_{r, \goth{f}, \goth{b}, \goth{a}}^{\mu}(u_1, \dots, u_r; h) = \prod_{\delta \in L/M} \frac{G_r\left(\frac{m'}{q} + \frac{\delta}{\gamma}, \frac{\alpha_0}{\gamma}, \dots, \frac{\alpha_r}{\gamma}\right)}{G_r\left(\frac{Nm'}{q} + \frac{N\delta}{\gamma}, \frac{N\alpha_0}{\gamma}, \dots, \frac{N\alpha_r}{\gamma}\right)}$$
Because $n>0$ the same $\alphabar, \gamma, M$ are used to compute
$$G_{r, \goth{f}, \goth{b}, \goth{a}}^{\mu}(u_1, \dots, u_r; nh) = \prod_{\delta \in L/M} \frac{G_r\left(\frac{m'}{q} + \frac{m'(n-1)}{q} + \frac{\delta}{\gamma}, \frac{\alpha_0}{\gamma}, \dots, \frac{\alpha_r}{\gamma}\right)}{G_r\left(\frac{Nm'}{q} + \frac{Nm'(n-1)}{q} + \frac{N\delta}{\gamma}, \frac{N\alpha_0}{\gamma}, \dots, \frac{N\alpha_r}{\gamma}\right)}$$
Because the function $G_r$ is $1$-periodic in its first argument and $(n-1)/q \in \zz$ we have the desired equality.
\end{proof}

\begin{general}{Lemma}\label{unitinvar}
Take $h$ an admissible vector for $\goth{f}, \goth{b}, \goth{a}$. Fix a system $u_1, \dots, u_r$ of linearly independent units. Then for all orientations $\mu, \nu$ and for all unit $\eps \in \ok^{\times}$:
$$G_{r, \goth{f}, \goth{b}, \goth{a}}^{\mu}(u_1, \dots, u_r; \eps h)^{\nu} = G_{r, \goth{f}, \goth{b}, \goth{a}}^{\mu\norm{\eps}}(u_1, \dots, u_r; h)^{\nu\norm{\eps}}$$
\end{general}

\begin{proof}
Let us write $\eta = \norm{\eps} = \pm 1$. We analyse how the construction in section \ref{geomsetup} changes when the vector $h$ is replaced with $\eps h$. First, the linear forms $a_h$ and $a_{\eps h}$ are related by: 
\begin{align*}
\lambda a_{\eps h} &= \det_{B_L}(\eps h, u_1 \eps h, \dots, u_r \eps h, \cdot) \\
\lambda a_{\eps h} &= \eta \det_{B_L}(h, u_1 h, \dots, u_r h, \eps^{-1} \cdot) \\
\lambda a_{\eps h} &= \lambda \eta a_h(\eps^{-1}\cdot)
\end{align*}
This implies directly that for all $0 \leq j \leq r$ we have $\alpha_{\eps h, j} = \eta \eps \alpha_{h,j}$ and 
$$s\gamma_{\eps h} = \det_{B_{\Lambda}}(a_{\eps h}, u_1 a_{\eps h}, \dots, u_r a_{\eps_h}, \cdot) = s\eta^r\eps\gamma_h$$
This also gives $\eta\eps L/M_h = L/M_{\eps h}$. Thus, using $\gamma = \gamma_h$, $\alpha_j = \alpha_{h, j}$, $M = M_h$ we get:
$$G_{r, \goth{f}, \goth{b}, \goth{a}}^{\mu}(u_1, \dots, u_r; \eps h)^{\nu} = \prod_{\delta \in L/M} \frac{G_r\left(\frac{\eps h +q\eta\eps\delta}{q\eta^r\eps\gamma}, \frac{\eta\eps\alpha_0}{\eta^r\eps\gamma}, \dots, \frac{\eta\eps\alpha_r}{\eta^r\eps\gamma}\right)^{\nu}}{G_r\left(\frac{N(\eps h + q\eta\eps \delta)}{q\eta^r\eps\gamma}, \frac{N\eta\eps\alpha_0}{\eta^r\eps\gamma}, \dots, \frac{N\eta\eps\alpha_r}{\eta^r\eps\gamma}\right)^{\nu}}$$
It then follows from formula \refp{inversiontotale} that
$$G_{r, \goth{f}, \goth{b}, \goth{a}}^{\mu}(u_1, \dots, u_r; \eps h)^{\nu} = \prod_{\delta \in L/M} \frac{G_r\left(\frac{\eta h +q\delta}{q\gamma}, \frac{\alpha_0}{\gamma}, \dots, \frac{\alpha_r}{\gamma}\right)^{\nu\eta^{r+1}}}{G_r\left(\frac{N (\eta h+q\delta)}{q\gamma}, \frac{N\alpha_0}{\gamma}, \dots, \frac{N\alpha_r}{\gamma}\right)^{\nu\eta^{r+1}}}$$
Using formula \refp{truecomputations} we recognize the right-hand side as the expression of
$$G_{r, \goth{f}, \goth{b}, \goth{a}}^{\mu\eta}(u_1, \dots, u_r; h)^{\nu\eta}= G_{r, \goth{f}, \goth{b}, \goth{a}}^{\mu\norm{\eps}}(u_1, \dots, u_r; h)^{\nu\norm{\eps}}$$
which gives the desired result.
\end{proof} 

Consider now two possible vectors $g, h$ associated to $[\eps_{\rho(1)}|\dots|\eps_{\rho(r)}]$ in Proposition \ref{propfinal}, for the same choice of helper ideal $\goth{I}$. This means that there is a fractional ideal $J$ of $\kk$ such that $(g) = mJ$ and $(h) = m'J$ for some integers $m,m' \in \zz_{> 0}$ coprime to $q$. Then the ideal $mJ$ is generated by $mh/m'$ so there exists a unit $\eps \in \ok^{\times}$ such that $g = m\eps h/m'$. Moreover, because $g, h \equiv 1 \mod \goth{f}$ we get $m\eps/m'\equiv 1 \mod \goth{f}$. This leads us to define
\begin{equation}\label{defprezfone}
\zz_{\goth{f}}^1 = \{(n, \eps) \in \zszcross{q} \times \ok^{\times}, n\eps \equiv 1 \mod \goth{f}\}
\end{equation}
In addition, if $n\eps \equiv n \eps' \equiv 1 \mod \goth{f}$ then $\eps\eps'^{-1} \equiv 1 \mod \goth{f}$ and we work under the assumption that all units congruent to $1$ mod $\goth{f}$ have norm $+1$. Therefore, $\norm{\eps} = \norm{\eps'}$. We may then consider the following quotient space.

\begin{general}{Definition}\label{defzfone}
The elements $(n, \eps)$ and $(n', \eps')$ in $\zz_{\goth{f}}^1$ are considered equivalent if $n = n'$ (remember that $n$ and $n'$ represent classes in $\zsz{q}^{\times}$). We write $\zfone = \zz_{\goth{f}}^1/\sim$ for the corresponding equivalence classes. This set forms a group which is naturally embedded in $\zszcross{q}$. 
\end{general}

It is easy to see that $\zfone$ contains at least $(1,1)$ and $(q-1, -1)$ and that these two elements differ unless $\kk$ is of even degree and $q = 2$. It follows from lemma \ref{unitinvar} that the product
$$\prod_{(n, \eps) \in \zfone} G_{r, \goth{f}, \goth{b}, \goth{a}}^{\mu}(u_1, \dots, u_r; n\eps h)^{\nu}$$ 
is well-defined for any orientations $\mu, \nu$ where it is understood that $(n,\eps) \in \zfone$ always comes with $n > 0$. We now want to prove the following invariance proposition:

\begin{general}{Proposition}
Consider two vectors $g, h$ associated to the unit system $[\eps_{\rho(1)}|\dots, \eps_{\rho(r)}]$ in Proposition \ref{propfinal} with the same helper ideal $\goth{I}$. Fix the orientations $\mu, \nu$. Then:
\begin{enumerate}
\item There exists $(n, \eps) \in \zfone$ such that
$$G_{r, \goth{f}, \goth{b}, \goth{a}}^{\mu}(u_1, \dots, u_r; g)^{\nu} = G_{r, \goth{f}, \goth{b}, \goth{a}}^{\mu}(u_1, \dots, u_r; n\eps h)^{\nu}$$
\item The following equality holds: 
$$\prod_{(n, \eps) \in \zfone} G_{r, \goth{f}, \goth{b}, \goth{a}}^{\mu}(u_1, \dots, u_r; n\eps g)^{\nu} = \prod_{(n, \eps) \in \zfone} G_{r, \goth{f}, \goth{b}, \goth{a}}^{\mu}(u_1, \dots, u_r; n\eps h)^{\nu}$$ 
\end{enumerate}

\end{general}

\begin{proof}
1. Write once again $g = m\eps h/m'$ where $m,m'>0$ are coprime to $q$ and $m\eps \equiv m' \mod \goth{f}$. This gives:
$$G_{r, \goth{f}, \goth{b}, \goth{a}}^{\mu}(u_1, \dots, u_r; g)^{\nu} = G_{r, \goth{f}, \goth{b}, \goth{a}}^{\mu}(u_1, \dots, u_r; m\eps h/m')^{\nu}$$
Consider a positive integer $m''$ such that $m'm'' \equiv 1 \mod q$. Then $mm''\eps \equiv m'm'' \equiv 1 \mod \goth{f}$, so $(mm'', \eps) \in \zfone$ and it follows from lemma \ref{lemmaintvar} that 
$$G_{r, \goth{f}, \goth{b}, \goth{a}}^{\mu}(u_1, \dots, u_r; m\eps h/m')^{\nu} = G_{r, \goth{f}, \goth{b}, \goth{a}}^{\mu}(u_1, \dots, u_r; mm''\eps h)^{\nu}$$
2. Write once again $g = m\eps h/m'$ with $m\eps \equiv m' \mod \goth{f}$ and fix a positive integer $m''$ such that $mm'' \equiv 1 \mod q$. Then for all $(n', \eps') \in \zfone$, we have $n'\eps' g = (mn')(\eps\eps')h/m'$ and
$$G_{r, \goth{f}, \goth{b}, \goth{a}}^{\mu}(u_1, \dots, u_r; n'\eps'g)^{\nu} = G_{r, \goth{f}, \goth{b}, \goth{a}}^{\mu}(u_1, \dots, u_r; mm''n'\eps\eps' h)^{\nu}$$
Then we get for the whole product:
$$\prod_{(n', \eps') \in \zfone}G_{r, \goth{f}, \goth{b}, \goth{a}}^{\mu}(u_1, \dots, u_r; n'\eps'g)^{\nu} = \prod_{(n', \eps') \in \zfone} G_{r, \goth{f}, \goth{b}, \goth{a}}^{\mu}(u_1, \dots, u_r; mm''n'\eps\eps'h)^{\nu}$$
with $mm''n'\eps\eps' \equiv mm''\eps \equiv m'm'' \equiv 1 \mod \goth{f}$. We put $n'' = mm''n', \eps'' = \eps\eps'$ and we get $(n'', \eps'') \in \zfone$. The map $(n', \eps') \to (n'', \eps'')$ defines a bijection of $\zfone$ and we get:
$$\prod_{(n', \eps') \in \zfone}G_{r, \goth{f}, \goth{b}, \goth{a}}^{\mu}(u_1, \dots, u_r; n'\eps'g)^{\nu} = \prod_{(n'', \eps'') \in \zfone} G_{r, \goth{f}, \goth{b}, \goth{a}}^{\mu}(u_1, \dots, u_r; n''\eps''h)^{\nu}$$
\end{proof}

We have thus successfully proved that for a fixed helper ideal $\goth{I}$ the product in the right-hand side of Conjecture \ref{conjecture} is independent of the choice of vectors $h_{\rho}$ in Proposition \ref{propfinal}. We now want to give another formulation for the term
$$\prod_{(n, \eps) \in \zfone} G_{r, \goth{f}, \goth{b}, \goth{a}}^{\mu}(u_1, \dots, u_r; n\eps h)^{\nu}$$
which is useful for computations. Let us define the following character on $\zfone$:
$$ \chi_{\goth{f}} := \begin{cases} \zfone \to \{-1, 1\} \\ (n, \eps) \to \norm{\eps} \end{cases}$$
This character is well-defined because units congruent to $1$ mod $\goth{f}$ all have norm $+1$ and it only depends on $n$, thus we will often write $\chi_{\goth{f}}(n)$ for $\chi_{\goth{f}}(n, \eps)$. This specific character is closely related to norm-type characters of the class group $Cl^{+}(\goth{f})$ as follows. Consider any such character $\chi$. This means that there exists a character $\chi_{\mathrm{finite}}$ of $(\ok/\goth{f})^{\times}$ such that for all element $a \in \kk$ coprime to $\goth{f}$, $\chi((a)) = \chi_{\mathrm{finite}}(a) sg(\norm{a})$.
Using this for $(n, \eps) \in \zfone$ we have 
$$\chi((n)) = \chi((n\eps)) = \chi_{\mathrm{finite}}(n\eps)\norm{\eps} = \norm {\eps}$$
because $n\eps \equiv 1 \mod \goth{f}$ so $\chi((n)) = \chi_{\goth{f}}(n)$. Therefore $\chi_{\goth{f}}(n)$ is the common value of all norm-type characters of the class group $Cl^{+}(\goth{f})$ at $(n)$. Note that the functional equation of Hecke $L$-functions in ATR fields implies that if $L'(0, \chi) \neq 0$ for a character $\chi$ of the class group $Cl^{+}(\goth{f})$ then $\chi$ is of norm-type (see for instance \cite{Neukirch} for the functional equation). It is worth noting that $\chi_{\goth{f}}(1) = 1$ and $\chi_{\goth{f}}(q-1) = (-1)^d$ where $d$ is the degree of $\kk$. We use this character to write using lemma \ref{unitinvar}:
\begin{equation}\label{chieq}
\prod_{(n, \eps) \in \zfone} G_{r, \goth{f}, \goth{b}, \goth{a}}^{\mu}(u_1, \dots, u_r; n\eps h)^{\nu} = \prod_{(n, \eps) \in \zfone} G_{r, \goth{f}, \goth{b}, \goth{a}}^{\mu\chi_{\goth{f}}(n)}(u_1, \dots, u_r; n h)^{\chi_{\goth{f}}(n)\nu}
\end{equation}
The right-hand side of formula \refp{chieq} may be explictly written as
\begin{multline}\label{completecomputations}
\prod_{(n, \eps) \in \zfone} G_{r, \goth{f}, \goth{b}, \goth{a}}^{\mu\chi_{\goth{f}}(n)}(u_1, \dots, u_r; n h)^{\chi_{\goth{f}}(n)\nu} = \\
\prod_{(n, \eps) \in \zfone} \prod_{\delta \in L/M} ~\left(\frac{G_r\left(\frac{\mu\chi_{\goth{f}}(n)nh+q\delta}{q\gamma}, \frac{\alpha_0}{\gamma}, \dots, \frac{\alpha_r}{\gamma} \right)^N}{G_r\left(\frac{\mu\chi_{\goth{f}}(n)nNh+Nq\delta}{q\gamma}, \frac{N\alpha_0}{\gamma}, \dots, \frac{N\alpha_r}{\gamma} \right)}\right)^{\nu\mu^r\chi_{\goth{f}}(n)^{r+1}}
\end{multline}
where $\alphabar, \gamma, M$ are defined in section \ref{geomsetup}. Most of the computation of this double product can be done at once for all terms since the parameters $\alpha_0/\gamma, \dots, \alpha_r/\gamma$ are fixed. Remember that in the context of ($\ast$), ($\ast\ast$) this double product reduces to a simple product as $M = L$ and we get:
\begin{multline}\label{completecomputationsstar}
\prod_{(n, \eps) \in \zfone} G_{r, \goth{f}, \goth{b}, \goth{a}}^{\mu\chi_{\goth{f}}(n)}(u_1, \dots, u_r; n h)^{\chi_{\goth{f}}(n)\nu} = \\
\prod_{(n, \eps) \in \zfone} ~\left(\frac{G_r\left(\frac{\mu\chi_{\goth{f}}(n)nh}{q\gamma}, \frac{\alpha_0}{\gamma}, \dots, \frac{\alpha_r}{\gamma} \right)^N}{G_r\left(\frac{\mu\chi_{\goth{f}}(n)nNh}{q\gamma}, \frac{N\alpha_0}{\gamma}, \dots, \frac{N\alpha_r}{\gamma} \right)}\right)^{\nu\mu^r\chi_{\goth{f}}(n)^{r+1}}
\end{multline}

\subsection{The main algorithm}

Before giving numerical evidence to support our conjecture, we explain how proposition \ref{propfinal} is used to compute the arithmetic $G_r$ functions. 
\begin{enumerate}
\item Fix $\kk$ a degree $d \geq 2$ ATR field and $\goth{f} \neq \ok$. Compute fundamental units for $\ok^{\times}$ (as given by Pari/GP for example). Use these units to compute a set of fundamental units for $\opc{\goth{f}}$ (see \refp{defopcf}). 
\item Compute the values $\ts_{\rho}, \ttt_{\rho}$ and $\tlambda_{\rho}$ for all $\rho$ (as defined in section \ref{geomsetup}). Check conditions ($\ast$) and ($\ast\ast$) and compute target ideals (see section \ref{choiceh}). Note that when $\ttt_{\rho} = 1$ there are no target ideals for the permutation $\rho$ and the empty product $\prod_{j = 1}^{j_{\rho}} I_{\rho, j}$ is equal to the trivial ideal $\ok$.
If the conditions ($\ast$) and ($\ast\ast$) are not satisfied, try small base changes for the units.
\item Choose $\goth{a}$ and $\goth{b}$ as in Proposition \ref{propfinal}, with the smallest possible norms.
\item Use Proposition \ref{propfinal} to construct the admissible generators $h_{\rho}$. One way to find $\goth{I}$ in Proposition \ref{propfinal} is to try $\goth{I} = \ok$ first, and if unsuccessful, try successively the available helper ideals with increasing norms until successful.
\item Compute the set $\zfone$ given in Definition \ref{defzfone} and the associated character $\chi_{\goth{f}}$ (see section \ref{sectionchoiceh}). 
\item Compute $h'_\rho$, $\abar_{\rho}$, $\alphabar_{\rho}$ as defined in section \ref{geomsetup} and compute the value $\prod_{(n,\eps) \in \zfone}G_{r, \goth{f}, \goth{b}, \goth{a}}^{\mu_{\rho}}([\eps_{\rho(1)} | \dots | \eps_{\rho(r)}] ; n\eps h_\rho)^{\nu_{\rho}}$ using formula \refp{sumcomputations} as well as formula \refp{completecomputationsstar} for the desired orientations.
\item Put everything together to compute the full product 
$$\prod_{\rho \in \goth{S}_{r}} \prod_{(n,\eps) \in \zfone} G_{r, \goth{f}, \goth{b}, \goth{a}}^{\mu_{\rho}}([\eps_{\rho(1)} | \dots | \eps_{\rho(r)}] ; n\eps h_\rho)^{\nu_{\rho}}$$
for the desired orientations $\bars{\mu}$ and $\bars{\nu}$.

\end{enumerate}

\section{Numerical evidence to support the conjecture}\label{sectionnumerical}

In this section, we provide numerical examples to support our conjecture. They may be computed with high precision in a low amount of time. In what follows, we will give computation times for $1000$ digits precision on a personal computer. Computations were carried out using number fields found in the LMFDB database \cite{lmfdb} as well as the computer algebra system PARI/GP \cite{parigp}, making extensive use of algebraic number theory tools it provides. To understand the result of the computations, we use the commands \textbf{lindep} to test if the value obtained from our computions is close to the expected value for $\norm{\goth{a}}\zeta'_{\goth{f}}([\goth{b}], 0) - \zeta'_{\goth{f}}([\goth{a}\goth{b}], 0)$ and $\textbf{algdep}$ to test if the value obtained is close to some algebraic integer in an abelian extension of $\kk$. Note that if $[\kk^{+}(\goth{f}):\kk] = m$ and $[\kk:\qq] = d$ then we may compute using our method $2m$ roots of a polynomial defining a subfield of $\kk^{+}(\goth{f})$ over $\qq$ out of a maximal amount of $dm$. This is achieved by varying the class $[\goth{b}]$ and using both complex embeddings $\sigma_{\cc}$ and $\overline{\sigma_{\cc}}$. This way, we may also compute all the roots of a polynomial in $\kk[x]$ which defines a subextension of $\kk^{+}(\goth{f})/\kk$ and we can thus compute the elementary symmetric polynomials in these roots to recover the polynomial. 

In what follows, we will define our fields as $\kk = \qq(z)$ where $z$ is the complex root of some polynomial $P = (X-z)(X-\bar{z})\prod_{j = 1}^r (X-z_j) \in \qq[X]$ lying in the upper half-plane (here the $z_j$ are real numbers with $z_1 < \dots < z_r$). Thus, to define our orientations we will fix the ordering on the real embeddings of $\kk$ such that $\sigma_{j}(z) = z_j$ for $1 \leq j \leq r$. Throughout this section, we focus mainly on computations for $\goth{b} = (1)$ because most of the work on target ideals is independent of the choice for the ideal $\goth{b}$. For ease of presentation, we will name prime ideals above a prime $p$ as $\goth{P}_p, \goth{P}'_p, \goth{P}''_p, \dots$ in their order of appearance using the commands \textbf{idealprimedec} or \textbf{idealfactor} in Pari/GP version 2.15.4. We provide examples in the cubic, quartic and quintic cases with $\kappa = \cardinalshort{\zfone}$ first and then quartic examples with $\kappa < \cardinalshort{\zfone}$. In all examples, the orientations have been computed as in Definition \ref{deforientations} and all values $v$ computed will satisfy directly formula \refp{eqconjecture} or
$$\norm{\goth{a}}\zeta'_{\goth{f}}([\goth{b}], 0) - \zeta'_{\goth{f}}([\goth{a}\goth{b}], 0) = \frac{\kappa}{\cardinalshort{\mathcal{Z}_{\goth{f}}^1}}\log|v|^2$$ which is an avatar of formula \refp{eqconjecture} with $\kappa$ defined in section \ref{formulation}. For ease of presentation we will not write the partition sets $U_{\rho, j}$ (see section \ref{formulation}) in the examples below.

\subsection{Cubic examples}

Here we present six examples in the cubic case. The first example showcases the work on target ideals with $\ttt = 131$ in the spirit of ($\ast$) and ($\ast\ast$). The second example is one of the simplest ones, with $\ttt = 1$ and class number one. The third example is a case where $q = 5$. The fourth and fifth examples showcase the work on helper ideals when the class number gets larger, for $q = 3$ and $q = 7$. The sixth example shows what can be done outside ($\ast$), ($\ast\ast$) in the cubic case. We also provide a table in Appendix \ref{cubicappendix} giving the values of $\tlambda$ and $\ttt$ for the pure cubic fields defined by $X^3-m$ for $ 2 \leq m \leq 50$ up to isomorphism.

\stepcounter{example}
\subsubsection{Example \exc}
We first discuss in full detail a cubic example. Let $\zexc$ be the complex root of the polynomial $x^3-13$ lying in the upper half-plane. Then $\kk_{\exc} = \qq(\zexc)$ has class number $3$. We choose the ideal $\goth{f}$ above $q = 3$ such that $\goth{f}^3 = (3)$. The corresponding narrow ray class group is $Cl^{+}(\goth{f}) \simeq \zsz{6}$. The unit group $\opc{\goth{f}}$ of totally positive units congruent to $1 \mod \goth{f}$ is generated by $\eps = 2\zexc^2-3\zexc -4$. 
\begin{itemize}
\item We compute $\tlambda = 1$, $\ts = \ttt = 131$ as defined in section \ref{geomsetup}.
\item Then, we search for target ideals above the prime number $131$ (see Definition \ref{deftarget}) and we find the unique degree one prime $\goth{P}_{131}$ above $131$. 
\item We choose for the smoothing ideal $\goth{a}$ the unique degree one prime above $N = 5$ in $\kk_{\exc}$. The ideals $\goth{P}_2^k$ for $0 \leq k \leq 5$ represent all classes in $Cl^{+}(\goth{f})$, where $\goth{P}_2$ is the only prime ideal of norm $2$ in $\kk$. These integral ideals satisfying the conditions in Proposition \ref{propfinal}.
\item Using Proposition \ref{propfinal} for $0 \leq k \leq 5$ we obtain the following ideals and admissible generators (see Definition \ref{defadmissible}):
\begin{center}
\begin{tabular}{| C | C | C | C | C |}
\hline
k & m & \goth{I} & \text{Ideal in Prop. \ref{propfinal}} & \text{Admissible generator } h_k \\
\hline
0 & 1 & (1) & qN\goth{P}_{131}/\goth{a} & -21\zexc^2 - 42\zexc - 114\\
\hline
1 & 1 & \goth{P}'_7\goth{P}''_7 & qN\goth{P}_{131}\goth{P}'_7\goth{P}''_7/(\goth{a}\goth{P}_2) & (87\zexc^2 + 9\zexc - 477)/2\\ 
\hline
2 & 1 & \goth{P}'_{11} & qN\goth{P}_{131}\goth{P}'_{11}/(\goth{a}\goth{P}_2^2) & -(3\zexc^2 + 309\zexc - 807)/4\\
\hline
3 & 1 & (1) & qN\goth{P}_{131}/(\goth{a}\goth{P}_2^3) & (-9\zexc^2 + 27\zexc + 159)/8 \\
\hline
4 & 1 & \goth{P}'_7\goth{P}''_7 & qN\goth{P}_{131}\goth{P}'_7\goth{P}''_7/(\goth{a}\goth{P}_2^4) & (-159\zexc^2 - 603\zexc - 1431)/16\\
\hline
5 & 1 & \goth{P}'_{11} & qN\goth{P}_{131}\goth{P}'_{11}/(\goth{a}\goth{P}_2^5) & (921\zexc^2 + 2157\zexc + 4209)/32\\
\hline
\end{tabular}
\end{center}
The corresponding levels will be $l_0 = l_3 = 3\cdot 5\cdot 131 = 1965$, $l_1 = l_4 = 7l_0$, $l_2 = l_5 = 11l_0$ as defined in section \ref{geomsetup}.
\item We compute the orientations $\mu = \nu = -1$ as in Definition \ref{deforientations}.
\item We compute $\zfone = \{ (1,1), (2,-1)\}$ (see Definition \ref{defzfone}) and we check that the equality $\Gamma^{-}_{\goth{f}, \goth{b}, \goth{a}}(\eps, -2h)^{-1} = \Gamma^{-}_{\goth{f}, \goth{b}, \goth{a}}(\eps, h)^{-1}$ holds. This gives $\kappa = 2$ as defined in section \ref{formulation}.
\end{itemize}
We may compute the values
\begin{align*}
\Gamma^{-}_{\goth{f}, (1), \goth{a}}(\eps, h_0)^{-1} = & \frac{\Gamma\left(\frac{-1}{3}, \frac{\eps^{-1} + 5348}{1965}, \frac{\eps+467}{1965}\right)^{-5}}{\Gamma\left(\frac{-5}{3}, \frac{\eps^{-1}+5348}{393}, \frac{\eps+467}{393}\right)^{-1}} \approx -0.0660917... + i\cdot0.0932299... \\
\Gamma^{-}_{\goth{f}, \goth{P}_2, \goth{a}}(\eps, h_1)^{-1} = & \frac{\Gamma\left(\frac{1}{3}, \frac{\eps^{-1} + 1418}{13755}, \frac{\eps-3463}{13755}\right)^{-5}}{\Gamma\left(\frac{5}{3}, \frac{\eps^{-1}+ 1418}{2751}, \frac{\eps-3463}{2751}\right)^{-1}} \approx 0.0059953... +i\cdot0.0047179... \\
\Gamma^{-}_{\goth{f}, \goth{P}_2^2, \goth{a}}(\eps, h_2)^{-1} = & \frac{\Gamma\left(\frac{-1}{3}, \frac{\eps^{-1} - 547}{21615}, \frac{\eps + 6362}{21615}\right)^{-5}}{\Gamma\left(\frac{-5}{3}, \frac{\eps^{-1} - 547}{4323}, \frac{\eps + 6362}{4323}\right)^{-1}} \approx -289.3045814... -i\cdot127.6382732... \\
\Gamma^{-}_{\goth{f}, \goth{P}_2^3, \goth{a}}(\eps, h_3)^{-1} = & \frac{\Gamma\left(\frac{1}{3}, \frac{\eps^{-1} + 5348}{1965}, \frac{\eps+467}{1965}\right)^{-5}}{\Gamma\left(\frac{5}{3}, \frac{\eps^{-1}+5348}{393}, \frac{\eps+467}{393}\right)^{-1}} \approx -5.0606452... - i\cdot7.1386178... \\
\Gamma^{-}_{\goth{f}, \goth{P}_2^4, \goth{a}}(\eps, h_4)^{-1} = & \frac{\Gamma\left(\frac{-1}{3}, \frac{\eps^{-1} + 1418}{13755}, \frac{\eps-3463}{13755}\right)^{-5}}{\Gamma\left(\frac{-5}{3}, \frac{\eps^{-1}+1418}{2751}, \frac{\eps-3463}{2751}\right)^{-1}} \approx 103.0063956... - i\cdot81.0605592... \\
\Gamma^{-}_{\goth{f}, \goth{P}_2^5, \goth{a}}(\eps, h_5)^{-1} = & \frac{\Gamma\left(\frac{1}{3}, \frac{\eps^{-1} -547}{21615}, \frac{\eps+6362}{21615}\right)^{-5}}{\Gamma\left(\frac{5}{3}, \frac{\eps^{-1}-547}{4323}, \frac{\eps+6362}{4323}\right)^{-1}} \approx -0.0028933... + i\cdot0.0012765... 
\end{align*}
Using all 6 values $I_k = \Gamma^{-}_{\goth{f}, \goth{P}_2^k, \goth{a}}(\eps, h_k)^{-1}$ we may compute and identify the relative polynomial:
\begin{align*}
\prod_{k = 0}^5 (X - I_k) =~& (X^6+1) + (-34\zexc^2 + 26\zexc + 128)(X^5+X)\\ &+  (1127\zexc^2 +8879\zexc -27106)(X^4+X^2) \\ &+ (40740\zexc^2-16965\zexc -185350)X^3
\end{align*}
Alternatively, all values $I_k$ are close to roots of the absolute palindromic polynomial $x^{18} + 384x^{17} + 2310x^{16} - 10646490x^{15} + 1596241353x^{14} + 18608357181x^{13} \linebreak + 156933809421x^{12} + 215098256580x^{11} + 381407365338x^{10} + 338205493469x^9 + \dots + 384x + 1$ 
which defines the degree $6$ extension $\kk_{\exc}^{+}(\goth{f})/\kk$. The computation time for $1000$ digits is 15 seconds for $I_0 = I_3^{-1}$, 65 seconds for $I_1 = I_4^{-1}$ and 105 seconds for $I_2 = I_5^{-1}$. We may also check formula \refp{eqconjecture} as
$$ \norm{\goth{a}}\zeta'_{\goth{f}}([(1)], 0) - \zeta'_{\goth{f}}([\goth{a}.(1)], 0) = \frac12\log\left|\prod_{(n, \eps') \in \zfone}\Gamma^{-}_{\goth{f}, \goth{b}, \goth{a}}(\eps, n\eps'h)^{-1}\right|^2 \approx -4.3382052...$$
and similarly for the other classes. Consider now another admissible vector $h = 3\zexc^2+6\zexc + 12$ for $\goth{b} = (1)$. The valuation of $h$ at any prime ideal dividing $131$ is $0$. This is one of the simplest admissible vectors we could find. Using this vector $h$ we have to compute $t = 131$ ordinary elliptic Gamma functions and the computation time is 30 seconds for $1000$ digits. In the cubic case, because the computation time for the precision $\delta$ is $O(l|\log(\delta)|t)$, reducing the value of $t$ by raising the level $l$ may not result in a massive time gain. In the cases $r \geq 2$, the gain will be more substantial.

\stepcounter{example}
\subsubsection{Example \exc}\label{cubicsimple}

We now discuss one of the simplest cubic examples which was presented in the introduction. Let $\zexc$ be the complex root of the polynomial $x^3-2$ in the upper half-plane. Then $\kk_{\exc} = \qq(\zexc)$ has class number $1$. We choose the ideal $\goth{f}$ above $q = 3$ such that $\goth{f}^3 = (3)$. The corresponding narrow ray class group is $Cl^{+}(\goth{f}) \simeq \zsz{2}$. The unit $\eps = \zexc -1$ is a generator for $\opc{\goth{f}}$ such that $\tlambda = \ttt = 1$. We choose $\goth{b} = (1)$ and $\goth{a}$ the degree one prime above $N = 5$ in $\kk_{\exc}$. The ideal $qN/\goth{a}$ has an admissible generator $h = -3\zexc^2 +6\zexc + 3$. The corresponding level is $l = 15$. Here we have $\zfone = \{(1,1), (2, -1)\}$ and $\kappa = 2$. The value
$$\Gamma^{-}_{\goth{f}, \goth{b}, \goth{a}}(\eps, h)^{-1} = \frac{\Gamma\left(\frac{-1}{3}, \frac{\eps^{-1}+2}{15}, \frac{\eps-7}{15}\right)^{-5}}{\Gamma\left(\frac{-5}{3}, \frac{\eps^{-1}+2}{3}, \frac{\eps-7}{3}\right)^{-1}} \approx -1.2937005... + i\cdot1.4743341...$$
is close to a root of the polynomial $x^6 + 3x^5 + 6x^4 + 5x^3 + 6x^2 + 3x + 1$ which defines $\kk_{\exc}^{+}(\goth{f})$. The computation time for $1000$ digits is 1 second.

\stepcounter{example}
\subsubsection{Example \exc}\label{excubictext}
Here we explain how we have obtained the example in section \ref{formulation} with $q = 5$. Let $\zexc =  e^{2i\pi/3}10^{1/3}$ be the root of the polynomial $x^3-10$ in the upper half-plane. Then $\kk_{\exc} = \qq(\zexc)$ has class number $1$. We choose the ideal $\goth{f}$ above $q = 5$ such that $\goth{f}^3 = (5)$. The corresponding narrow ray class group is $Cl^{+}(\goth{f}) \simeq \zsz{4}$. The unit $\eps = (2\zexc^2 - \zexc - 7)/3$ is a generator for $\opc{\goth{f}}$ such that $\tlambda = 1$ and $\ttt = 9$. We may choose $\goth{b} = (1)$ and $\goth{a}$ the unique degree one prime above $N = 11$ in $\kk_{\exc}$. We find the target ideal $\goth{P}_3^2$ above $3$ where $\goth{P}_3\goth{P}_3^{\prime2} = (3)$. Then, $h = -(35z^2 + 20z + 35)/3$ is an admissible generator for the ideal $qN\goth{P}_3^2/\goth{a}$. Here we have $\zfone = \{(1,1), (4, -1)\}$ and $\kappa = 2$. The value
$$\Gamma^{-}_{\goth{f}, \goth{b}, \goth{a}}(\eps, h)^{-1} = \frac{\Gamma\left(\frac{-1}{5}, \frac{\eps^{-1}-1751}{495}, \frac{\eps-776}{495}\right)^{-11}}{\Gamma\left(\frac{-11}{5}, \frac{\eps^{-1}-1751}{45}, \frac{\eps-776}{45}\right)^{-1}} \approx -27.5333588... - i\cdot32.7146180...$$
is close to a root of the polynomial $x^{12} + 57x^{11} + 1956x^{10} + 4640x^9 + 35415x^8 - 109818x^7 + 150139x^6 - 109818x^5 + 35415x^4 + 4640x^3 + 1956x^2 + 57x + 1$ which defines $\kk_{\exc}^{+}(\goth{f})$. The computation time for 1000 digits is 7 seconds.

\stepcounter{example}
\subsubsection{Example \exc}

Let $\zexc$ be the complex root of the polynomial $x^3-65$ in the upper half-plane. Then $\kk_{\exc} = \qq(\zexc)$ has class number $18$. This means that most ideals won't be principal ideals and we will need to use helper ideals to build an admissible vector $h$. We choose the ideal $\goth{f}$ above $q = 3$ such that $\goth{f}^3 = (3)$. The corresponding narrow ray class group is $Cl^{+}(\goth{f}) \simeq \zsz{6}\times \zsz{6}$. The unit $\eps = \zexc -4$ is a generator for $\opc{\goth{f}}$ such that $\tlambda = \ttt = 1$. We choose $\goth{b} = (1), \goth{b} = \goth{P}_{59}, \dots$ representatives for the $36$ classes in $Cl^{+}(\goth{f})$, and $\goth{a}$ the unique degree one prime above $N = 5$ in $\kk_{\exc}$. The ideals $qN/(\goth{a}\goth{b})$ are unfortunately not principal in general so we use helper ideals. In the case $\goth{b} = (1)$ we find the helper ideal $\goth{P}_{13}^2$ and in the case $\goth{b} = \goth{P}_{59}$ we find the helper ideal $\goth{P}'_{41}$. The corresponding admissible generators we found for $qN\goth{P}_{13}^2/\goth{a}$ and $qN\goth{P}'_{41}/(\goth{a}\goth{P}_{59})$ are $h_{\goth{b} = (1)} = 3\zexc^2$ and $h_{\goth{b} = \goth{P}_{59}} = -(72\zexc^2 + 75\zexc + 1590)/59$ respectively. The corresponding levels will be $l_{(1)} = 3\cdot5\cdot13$, $l_{\goth{P}_{59}} = 3\cdot5\cdot41$. Here we have $\zfone = \{(1,1), (2, -1)\}$ and $\kappa = 2$. We may compute
\begin{align*}
\Gamma^{-}_{\goth{f}, (1), \goth{a}}(\eps, h_{\goth{b} = (1)})^{-1} =
 \frac{\Gamma\left(\frac{-1}{3}, \frac{\eps^{-1}-211}{195}, \frac{\eps-61}{195}\right)^{-5}}{\Gamma\left(\frac{-5}{3}, \frac{\eps^{-1}-211}{39}, \frac{\eps-61}{39}\right)^{-1}} &\approx -1.6691052... +i\cdot5.7493283...\hfill\\
\Gamma^{-}_{\goth{f}, \goth{P}_{59}, \goth{a}}(\eps, h_{\goth{b} = \goth{P}_{59}})^{-1} =
 \frac{\Gamma\left(\frac{1}{3}, \frac{\eps^{-1}+1034}{615}, \frac{\eps-91}{615}\right)^{-5}}{\Gamma\left(\frac{5}{3}, \frac{\eps^{-1}+1034}{123}, \frac{\eps-91}{123}\right)^{-1}} &\approx 0.0344135... - i\cdot0.0123218...\hfill
\end{align*}
 and the remaining 34 out of 36 values $I_{\goth{b}} = \Gamma_{\goth{f}, \goth{b}, \goth{a}}^{-}(\eps, h_{\goth{b}})^{-1}$ attached to the 36 classes in $Cl^{+}(\goth{f})$. We may then compute the polynomial 
 $$\prod_{\goth{b} \in Cl^{+}(\goth{f})} (X - I_{\goth{b}}) = X^{36}-(u +v\zexc +w\zexc^2)X^{35} + \dots + 1 \in \mathcal{O}_{\kk_{\exc}}[X]$$
where $u = 93377174024326, v = 769211619985, w = -5967373310133$. This palindromic polynomial 
defines a relative equation of the class field $\kk_{\exc}^{+}(\goth{f})$ above $\kk_{\exc}$ and we identify the rest of its coefficients in $\mathcal{O}_{\kk_{\exc}}$. The computation time for $1000$ digits and for all of the 36 computations is 8 minutes, which gives 13 seconds per individual computation on average.

It is interesting to note that this field as well as the field defined by $x^3-2$ belong to a special family parametrized by $x^3 - 2^{3k} - (-1)^k$ for $k \geq 0$ which behaves nicely compared to other pure cubic fields with respect to the values $\tlambda, \ttt$. Generally speaking, the values of these parameters vary significantly with the discriminant, as shown in the table in Appendix \ref{cubicappendix}. Putting $\bb{L}_k = \qq(z_k)$ where $z_k$ is the complex root of the polynomial $x^3 - 2^{3k} - (-1)^k$ lying in the upper half-plane, it seems that the positive fundamental unit in $\bb{L}_k$ is given by $\eps_k = (-1)^k(z_k-2^k)$ and its inverse by $\eps_k^{-1} = z_k^2 + 2^kz_k + 2^{2k}$. When $D_k = 2^{3k} + (-1)^k$ is cube-free (at least for $0 \leq k \leq 100$ except for $k = 49, 50$), because $D_k \not\equiv \pm 1 \mod 9$, the field $\bb{L}_k$ is a so-called pure cubic field of the first kind and an integral basis of $\mathcal{O}_{\bb{L}_k}$ is given by $(1, z_k, z_k^2/s_k)$ where $D_k = r_k s_k^2$ with $r_k, s_k$ square-free and relatively coprime. In that case, we may compute $\tlambda = 1$ and $\ttt = s_k$. When $D_k$ is also square-free (at least for $0 \leq k \leq 100$ except for $k = 7, 10, 21, 26, 30, 35, 63, 68, 70, 77, 78, 90, 91$) we get $\tlambda_k = \ttt_k = 1$. The class group grows very rapidly in this family (ex: $h(\bb{L}_{11}) = 2^2\cdot3^5\cdot5\cdot 3191$ and
$h(\bb{L}_{20}) = 2^6\cdot3^4\cdot5\cdot11\cdot19\cdot79\cdot 863\cdot18047$). It would be interesting to understand if for $k \geq 100$ there are infinitely many fields in this family for which the conditions $\tlambda = \ttt = 1$ are satisfied for the fundamental unit. Our construction would then conjecturally allow to construct algebraic units of very high degree.

\stepcounter{example}
\subsubsection{Example \exc}
Here we give an example with $q = 7$. Let $\zexc$ be the complex root of the polynomial $x^3-14$ in the upper half-plane. Then $\kk_{\exc} = \qq(\zexc)$ has class number $3$. We choose the ideal $\goth{f}$ above $q = 7$ such that $\goth{f}^3 = (7)$. The corresponding narrow ray class group is $Cl^{+}(\goth{f}) \simeq \zsz{18}$. The unit $\eps = -\zexc^2+2\zexc +1$ is a generator for $\opc{\goth{f}}$ such that $\tlambda = 1$, $\ttt = 2\cdot11$. We choose $\goth{b} = (1)$ and $\goth{a}$ the unique degree one prime above $5$ in $\kk_{\exc}$. Computing target ideals gives $\goth{P}_2$ and $\goth{P}_{11}$ the degree one primes above $2$ and $11$. The ideal $qN\goth{P}_2\goth{P}_{11}/(\goth{a}\goth{b})$ is unfortunately not principal so we need to look for a helper ideal. In Proposition \ref{propfinal} we may take $m_1 = 3$ and $\goth{I} = \goth{P}_3^2$ so that the ideal $3qN\goth{P}_2\goth{P}_{11}\goth{P}_3^2/(\goth{a}\goth{b})$ is principal with an admissible generator $h = -21\zexc^2 + 21\zexc - 336$. The corresponding level will be $l = 2\cdot3\cdot5\cdot7\cdot11$. Here we have $\zfone = \{(1,1), (6, -1)\}$ and $\kappa = 2$. The value
$$\Gamma^{-}_{\goth{f}, \goth{b}, \goth{a}}(\eps, h)^{-1} =
 \frac{\Gamma\left(\frac{3}{7}, \frac{\eps^{-1}-3067}{2310}, \frac{\eps+1007}{2310}\right)^{-5}}{\Gamma\left(\frac{15}{7}, \frac{\eps^{-1}-3067}{462}, \frac{\eps+1007}{462}\right)^{-1}} \approx -0.1700923... +i\cdot3.8609499...$$
is close to a root of the palindromic polynomial $x^{54} - 4167x^{53} + 7931535x^{52} - 259219286x^{51}... -4167x +1$ which defines $\kk_{\exc}^{+}(\goth{f})$ over $\qq$. This polynomial has very large coefficients, and we could alternatively compute the remaining 17 out of 18 values associated to the 18 classes in $Cl^{+}(\goth{f})$ to identify a relative polynomial in $\mathcal{O}_{\kk_{\exc}}[x]$ instead, as we did in the previous example. The computation time for $1000$ digits and for all $18$ values is 18 minutes, which gives 1 minute per individual computation on average.

\stepcounter{example}
\subsubsection{Example \exc}
This example falls outside of the conditions ($\ast$) and ($\ast\ast$) because we will have $\tlambda = 2$. Yet, it shows that these conditions are a set of $p$-adic conditions which can be satisfied for most $p$. In other words, the unique prime number for which ($\ast$) and ($\ast\ast$) are not true in the present case is $2$. Because of this, we can try to reduce the number of computations outside of the prime number $2$.

Let $\zexc$ be the complex root of the polynomial $x^3-5$ in the upper half-plane. Then $\kk_{\exc} = \qq(\zexc)$ has class number $1$. We choose the ideal $\goth{f}$ above $q= 3$ such that $\goth{f}^3 = (3)$. The corresponding narrow ray class group is $Cl^{+}(\goth{f}) \simeq \zsz{2}$. The unit $\eps = 2\zexc^2-4\zexc+1$ is a generator for $\opc{\goth{f}}$ such that $\tlambda = 2$ and $\ttt = 2\times 13$. We choose $\goth{b} = 1$ and $\goth{a}$ the unique degree one prime above $N = 5$ in $\kk_{\exc}$. We find $\goth{P}_{13}$ a target ideal above $13$ and we find an admissible generator $h = 6\zexc^2+15$ for the ideal $qN\goth{P}_{13}/\goth{a}$. The corresponding level will be $l = 3\cdot 5\cdot 13$. Here we have $\zfone = \{(1,1), (2, -1)\}$ and $\kappa = 2$. Then $\Gamma^{-}_{\goth{f}, \goth{b}, \goth{a}}(\eps, h)^{-1}$ is a product of two ordinary smoothed elliptic Gamma functions with parameters
$$\tau = \frac{\eps^{-1}+119}{2\times195},~\sigma = \frac{\eps + 59}{2 \times 195}$$
More precisely, put $\delta = -2\zexc^2-\zexc-5$. Then the set $L/M$ in formula \refp{truecomputations} is exactly $\{0, \delta\}$ with $\delta/\gamma = (-4\zexc^2 - 5\zexc - 55)/195$. The product
$$ \frac{\Gamma\left(\frac{1}{3}, \frac{\eps^{-1}+119}{390}, \frac{\eps+59}{390}\right)^{-5}}{\Gamma\left(\frac{5}{3}, \frac{\eps^{-1}+119}{78}, \frac{\eps+59}{78}\right)^{-1}}\cdot\frac{\Gamma\left(\frac{1}{3}+\frac{-4\zexc^2 - 5\zexc - 55}{195}, \frac{\eps^{-1}+119}{390}, \frac{\eps+59}{390}\right)^{-5}}{\Gamma\left(\frac{5}{3}+\frac{-4\zexc^2 - 5\zexc - 55}{39}, \frac{\eps^{-1}+119}{78}, \frac{\eps+59}{78}\right)^{-1}}$$
is close to the root $\approx 9.8439696... + i\cdot5.1060682...$ of the palindromic polynomial $x^6 - 21x^5 + 150x^4 - 185x^3 + 150x^2 - 21x + 1$ which defines $\kk_{\exc}^{+}(\goth{f})$. The computation time for $1000$ digits is 11 seconds.

\subsection{Quartic examples}\label{quarticexamples}

We now present four examples in the quartic case. The first example is one of the simplest cases, for which $\ttt_1 = \ttt_2 = 1$. The second example showcases the work on target ideals with $\ttt_1 = 7, \ttt_2 = 7129$. The third example shows that we can sometimes hope to find better units to gain time in the computations. The fourth example is a case where the quartic field $\kk$ contains a real quadratic field and the set of units falls just outside ($\ast$), ($\ast\ast$). In all examples, we use the ordering $\{ Id, (21) \}$ of $\goth{S}_2$ and when we say that we ``choose'' fundamental units for $\opc{\goth{f}}$, we let Pari/GP compute a set of fundamental units for $\ok^{\times}$, which we use to compute a set of fundamental units for $\opc{\goth{f}}$. Most of the time Pari/GP does a great job in finding small fundamental units, and the few operations required to compute $\opc{\goth{f}}$ keep these units relatively small. 

\stepcounter{example}
\subsubsection{Example \exc}\label{quarticsimple}

We first discuss one of the simplest quartic examples, mentioned in section \ref{formulation}, this time with a ``big smoothing''. Let $\zexc$ be the complex root of the polynomial $x^4 -6x^3-x^2-3x+1$ lying in the upper half-plane. Then $\kk_{\exc} = \qq(\zexc)$ has class number $1$. We choose $\goth{f} = \goth{P}_2$ the degree one prime above $q = 2$. The corresponding narrow ray class group is $Cl^{+}(\goth{f}) \simeq \zsz{2}$. We choose the fundamental units
$$ \eps_1 = \frac{-2\zexc^3 + 13\zexc^2 - \zexc + 3}{7},~\eps_2 = \frac{-5\zexc^3 + 29\zexc^2 + 15\zexc + 18}{7} $$
for the set $\opc{\goth{f}}$ of totally positive units congruent to $1 \mod \goth{f}$. We compute the contents $\tlambda_1 = \tlambda_2 = 1$, the overflows $\ttt_1 = \ttt_2 = 1$ as defined in section \ref{geomsetup}. We may choose $\goth{b} = (1)$ and $\goth{a}$ the unique degree one prime above $N = 13$ in $\kk_{\exc}$. Since $\ttt_1 = \ttt_2 = 1$, there are no target ideals to be found here. The ideal $qN/\goth{a}$ has an admissible generator $h_1 = h_2 = (44\zexc^3 - 258\zexc^2 - 104\zexc - 80)/7$ (see Definition \ref{defadmissible}). The corresponding levels will be $l_1 = l_2 = 2\cdot13$. We compute the orientations $\bars{\mu} = \bars{\nu} = [-1, 1]$ as given in Definition \ref{deforientations}. Here, we have $\zfone = \{ (1,1) \}$ as in Definition \ref{defzfone} and $\kappa = 1$. Let us write the parameters

\begin{align*}
&\tau = \eps_2 -15, && \tau' = -6 +\frac{1}{\eps_1} + \frac{1}{\eps_1\eps_2} \\
&\sigma = -7+\frac{1}{\eps_2}, && \sigma' = -\eps_2+15 \\
&\rho = -\eps_1-3, && \rho' = 4\eps_1+19-\frac{1}{\eps_2}
\end{align*}

\noindent Then 
$$\frac{G_{2, \goth{f}, \goth{b}, \goth{a}}^{+}(\eps_2, \eps_1\eps_2; h_2)}{G_{2, \goth{f}, \goth{b}, \goth{a}}^{-}(\eps_1, \eps_1\eps_2; h_1)} = \frac{G_2\left(\frac{-1}{2}, \frac{\tau}{26}, \frac{\sigma}{26}, \frac{\rho}{26}\right)^{-13}}{G_2\left(\frac{-13}{2}, \frac{\tau}{2}, \frac{\sigma}{2}, \frac{\rho}{2}\right)^{-1}} \times \frac{G_2\left(\frac{1}{2}, \frac{\tau'}{26}, \frac{\sigma'}{26}, \frac{\rho'}{26}\right)^{13}}{G_2\left(\frac{13}{2}, \frac{\tau'}{2}, \frac{\sigma'}{2}, \frac{\rho'}{2}\right)}
$$
is close to the root $ \approx 4.1210208... - i\cdot5.0617720...$ of the polynomial $x^8 - 7x^7 + 33x^6 + 49x^5 + 17x^4 + 49x^3 + 33x^2 - 7x + 1$ which defines $\kk_{\exc}^{+}(\goth{f})$. The computation time for $1000$ digits is 6 seconds. We may also check formula \refp{eqconjecture} up to 1000 digits as:
$$\norm{\goth{a}}\zeta'_{\goth{f}}([\goth{b}], 0) - \zeta'_{\goth{f}}([\goth{a}\goth{b}], 0) \approx \log\left|\frac{G_{2, \goth{f}, \goth{b}, \goth{a}}^{+}(\eps_2, \eps_1\eps_2; h_2)}{G_{2, \goth{f}, \goth{b}, \goth{a}}^{-}(\eps_1, \eps_1\eps_2; h_1)}\right|^2 \approx 3.7519563...$$
We compare the above computation to a simple choice for the admissible vectors $h_1, h_2$ given by Pari/GP. Namely for $h_1 = h_2 = (86\zexc^3 - 496\zexc^2 - 286\zexc - 206)/7$ we compute $t_1 = t_2 = 31$ and $l = 26$. The computation time for this pair and for 1000 digits is 2 minutes, which is roughly $20$ times longer than for the previous pair.

\stepcounter{example}
\subsubsection{Example \exc}
We now discuss in full detail a quartic example with specific work on target ideals. Let $\zexc$ be the complex root of the polynomial $x^4 - 19x^3 + 18x^2 + 8x + 1$ lying in the upper half-plane. Then $\kk_{\exc} = \qq(\zexc)$ has class number $1$. We choose $\goth{f} = \goth{P}_3$ a prime ideal above $q = 3$ such that $(3) = \goth{P}_3^2\goth{P}_3^{\prime 2}$. The corresponding narrow ray class group is $Cl^{+}(\goth{f}) \simeq \zsz{2}\times\zsz{2}$. We choose the fundamental units
$$\eps_1 = \frac{19\zexc^3 - 366\zexc^2 + 438\zexc + 44}{9},~\eps_2 = -\zexc^3 + 19\zexc^2 - 18\zexc - 8 $$
for $\opc{\goth{f}}$. We compute $\tlambda_1 = \tlambda_2 = 1$, $\ttt_1 = 7$, $\ttt_2 = 7129$ where $7129$ is a prime number. We may choose $\goth{b} = (1)$ and $\goth{a}$ the unique degree one prime above $N = 13$ in $\kk_{\exc}$. We search for target ideals and find $\goth{P}_7$ in the first case and $\goth{P}_{7129}$ in the second case. Both ideals $qN\goth{P}_7/\goth{a}$ and $qN\goth{P}_{7129}/\goth{a}$ are principal with admissible 
generators $h_1 = (-32\zexc^3 + 606\zexc^2 - 543\zexc- 112)/3$ and $h_2 = (124\zexc^3 - 2319\zexc^2 + 1563\zexc + 200)/3$. The corresponding levels will be $l_1 = 3\cdot7\cdot13$ and $l_2 = 3\cdot13\cdot7129$. Here we have $\zfone = \{(1,1), (2, 1)\}$ and $\kappa = 2$. Let us write the parameters

\begin{align*}
&\tau = \frac{2\eps_1\eps_2 -92\eps_1-3}{\eps_1}, && \tau' = \frac{57488\eps_1\eps_2 + 3\eps_2+16}{\eps_1\eps_2} \\
&\sigma = \frac{925\eps_1\eps_2+367\eps_1-47}{\eps_1 \eps_2}, && \sigma' = \frac{47\eps_1\eps_2 -328\eps_2-348694}{1} \\
&\rho = \frac{-\eps_1\eps_2-3\eps_1+40}{1}, && \rho' = \frac{-57\eps_1\eps_2-90004\eps_2-2}{\eps_2}
\end{align*}

\noindent Then 
$$\frac{G_{2, \goth{f}, \goth{b}, \goth{a}}^{+}(\eps_1, \eps_1\eps_2; h_1)}{G_{2, \goth{f}, \goth{b}, \goth{a}}^{-}(\eps_2, \eps_1\eps_2; h_2)} = \frac{G_2\left(\frac{1}{3}, \frac{\tau}{273}, \frac{\sigma}{273}, \frac{\rho}{273}\right)^{13}}{G_2\left(\frac{13}{3}, \frac{\tau}{21}, \frac{\sigma}{21}, \frac{\rho}{21}\right)} \times \frac{G_2\left(\frac{-2}{3}, \frac{\tau'}{278031}, \frac{\sigma'}{278031}, \frac{\rho'}{278031}\right)^{-13}}{G_2\left(\frac{-26}{3}, \frac{\tau'}{21387}, \frac{\sigma'}{21387}, \frac{\rho'}{21387}\right)^{-1}}
$$
is close to the root $\approx 10.6409709... - i\cdot 5.9332732...$ of the polynomial $x^8 - 18x^7 + 83x^6 + 396x^5 + 597x^4 + 396x^3 + 83x^2 - 18x + 1$ which defines a subextension of $\kk_{\exc}^{+}(\goth{f})/\kk_{\exc}$. The computation time for $1000$ digits is 3 minutes and 20 seconds.

\stepcounter{example}
\subsubsection{Example \exc}
\addtocounter{example}{-1}

Keep $\kk_{\exc}$, $\goth{f}$, $\goth{b}$, $\goth{a}$ as in example \exc\,  and let us change our choice of fundamental units. We fix another set of fundamental units for $\opc{\goth{f}}$:

$$ \eps_1 = \frac{\zexc^3 - 21\zexc^2 + 54\zexc + 11}{9},~\eps_2 = -\zexc^3 + 19\zexc^2 - 18\zexc - 8 $$
Then, we compute $\tlambda_1 = \tlambda_2 = 1$, $\ttt_1 = 1$, $\ttt_2 = 25$. In the second computation, a target ideal is $\goth{P}_5^2$ the square of the unique degree one prime above $5$ and the ideals $qN/\goth{a}$ and $qN\goth{P}_5^2/\goth{a}$ are generated by the admissible vectors $h_1 =(-44\zexc^3 + 843\zexc^2 - 927\zexc - 232)/3$, $h_2 = (76\zexc^3 - 1449\zexc^2 + 1470\zexc + 344)/3$. The corresponding levels will be $l_1=3\cdot13$ and $l_2 = 3\cdot5^2\cdot13$. Let us write the parameters

\begin{align*}
&\tau = \frac{2\eps_1\eps_2 +109\eps_1-3}{\eps_1}, && \tau' = \frac{842\eps_1\eps_2 +\eps_2+3}{\eps_1\eps_2}  \\
&\sigma = \frac{-321\eps_1\eps_2-16\eps_1+7}{\eps_1 \eps_2}, && \sigma' = \frac{-7\eps_1\eps_2 +4\eps_2-63}{1} \\
&\rho = \frac{-\eps_1-47}{1}, && \rho' = \frac{-3\eps_1\eps_2+566\eps_2-2}{\eps_2}
\end{align*}

\noindent Then 
$$\frac{G_{2, \goth{f}, \goth{b}, \goth{a}}^{+}(\eps_2, \eps_1\eps_2; h_2)}{G_{2, \goth{f}, \goth{b}, \goth{a}}^{-}(\eps_1, \eps_1\eps_2; h_1)} = \frac{G_2\left(\frac{-1}{3}, \frac{\tau}{39}, \frac{\sigma}{39}, \frac{\rho}{39}\right)^{-13}}{G_2\left(\frac{-13}{3}, \frac{\tau}{3}, \frac{\sigma}{3}, \frac{\rho}{3}\right)^{-1}} \times \frac{G_2\left(\frac{2}{3}, \frac{\tau'}{975}, \frac{\sigma'}{975}, \frac{\rho'}{975}\right)^{13}}{G_2\left(\frac{26}{3}, \frac{\tau'}{75}, \frac{\sigma'}{75}, \frac{\rho'}{75}\right)}
$$

\noindent is close to the same root $\approx 10.6409709... - i\cdot 5.9332732...$ of the polynomial $x^8 - 18x^7 + 83x^6 + 396x^5 + 597x^4 + 396x^3 + 83x^2 - 18x + 1$. The computation time for $1000$ digits is 11 seconds. This value is the same as for the previous example in which the choice of fundamental units was different, which supports the idea that the product in Conjecture \ref{conjecture} is independent of said choice, even though computations may be longer when a poor choice of fundamental units is made.

\addtocounter{example}{2}
\subsubsection{Example \exc}\label{quarticreal}

This example falls outside of ($\ast$) and ($\ast\ast$) because $gcd(\ttt_1, \ttt_2, q = 2) = 2$. Still, we can find a pair of vectors $h_1, h_2$ for which the computations give interesting results.  Let $\zexc$ be the complex root of the polynomial $x^4-12$ lying in the upper half-plane. Then $\kk_{\exc} = \qq(\zexc)$ has class number $1$ and \textit{contains the real quadratic field $\qq(\sqrt{3})$}. Remember that this means that we must carefully choose the fundamental units in this case to avoid totally real elements (see ($\ast$) conditions in section \ref{choiceh}). We choose the ideal $\goth{f}$ above $q = 2$ such that $\goth{f}^4 = (2)$. The corresponding narrow ray class group is $Cl^{+}(\goth{f}) \simeq \zsz{2}$. We choose the fundamental units
$$ \eps_1 = \frac{\zexc^2 + 2\zexc + 2}{4},~\eps_2 = \frac{-\zexc^3 + \zexc^2 + 4\zexc - 2}{4} $$
for $\opc{\goth{f}}$. We compute $\tlambda_1 = \tlambda_2 = 1$, $\ttt_1 = 2\cdot3$, $\ttt_2 = 2\cdot 11$. We may choose $\goth{b} = (1)$ and $\goth{a}=\goth{P}_{23}$ a degree one prime above $N = 23$ in $\kk_{\exc}$. We search for target ideals and find $\goth{P}_3$ the degree one prime ideal above $3$ in the first case and $\goth{P}_{11}$ one of the degree one primes above $11$ in the second case. Then $qN\goth{P}_3/\goth{a}$ and $qN\goth{P}_{11}/\goth{a}$ are generated by the admissible vectors $h_1 = (-4\zexc^3 + 11\zexc^2 + 10\zexc + 30)/2$ and $h_2 = (9\zexc^3 + 4\zexc^2 - 34\zexc - 56)/2$. The corresponding levels will be $l_1 = 2\cdot3\cdot23$ and $l_2 = 2\cdot11\cdot23$. Here we have $\zfone = \{(1,1)\}$ and $\kappa = 1$. We compute the orientations $\bars{\mu} = \bars{\nu} =  [1, -1]$ and the value $I_{2, \goth{f}, \goth{b}, \goth{a}}(\eps_1, \eps_2, \bars{h}, \bars{\mu}, \bars{\nu}) \approx 13.9102308... - i\cdot24.0932265...$ with $t_1 = t_2 = 2$, i.e. we compute a product of 8 ordinary elliptic smoothed $G_2$ functions. The result is close to a root of the polynomial $x^8 - 28x^7 + 778x^6 - 112x^5 - 749x^4 - 112x^3 + 778x^2 - 28x + 1$ which defines $\kk_{\exc}^{+}(\goth{f})$. The computation time for $1000$ digits is 8 seconds. 

Now we would like to discuss the case where some of the $\talpha_j$ belong to $\rr$. Choose the fundamental units:
$$ \eps_1 = \frac{\zexc^2 + 2\zexc + 2}{4},~\eps_2 = \frac{-\zexc^2 + 4}{2} = 2+ \sqrt{3} $$
Then we may compute $\tlambda_1 = \tlambda_2 = 1$, $\ttt_1 = \ttt_2 = 1$. Keeping the same ideals $\goth{b} = (1)$, $\goth{a} = \goth{P}_{23}$, we find the admissible vectors $h_1 = h_2 = (-3\zexc^3 - 9\zexc^2 - 4\zexc + 34)/2$ for $qN/\goth{P}_{23}$. We write the obtained parameters as:

\begin{align*}
&\tau = \eps_2 +175 = 177 + \sqrt{3}, && \tau' = 39+\frac{1}{\eps_1} \\
&\sigma = -342+\frac{3}{\eps_2}-\frac{1}{\eps_1 \eps_2}, && \sigma' = \eps_2\eps_1+\eps_2-352 \\
&\rho = -\eps_1+217, && \rho' = 179-\frac{1}{\eps_2} = 177 + \sqrt{3}
\end{align*}
We get in particular $\tau, \rho' \in \rr$, $\Im(\tau'), \Im(\rho) < 0$ and $\Im(\sigma), \Im(\sigma') > 0$. In this setting, as suggested by my advisor P. Charollois, it is possible to make sense of formula \refp{exponentialformula}. Indeed, using Liouville's theorem on the bad approximation of real irrational algebraic numbers by rational numbers, we get that 
$$\frac{1}{\sin \pi j \tau} = O(j^2),~\frac{1}{\sin \pi j \rho'} = O(j^2)$$
since $\tau, \rho'$ are algebraic numbers of degree $2$. This guarantees that  formula \refp{exponentialformula} converges in this setting for $z \in \qq$. The computation of the product
$$\frac{G_{2, \goth{f}, \goth{b}, \goth{a}}^{+}(\eps_1, \eps_1\eps_2; h_2)}{G_{2, \goth{f}, \goth{b}, \goth{a}}^{-}(\eps_2, \eps_1\eps_2; h_1)} = \frac{G_2\left(\frac{1}{2}, \frac{\tau}{46}, \frac{\sigma}{46}, \frac{\rho}{46}\right)^{23}}{G_2\left(\frac{23}{2}, \frac{\tau}{2}, \frac{\sigma}{2}, \frac{\rho}{2}\right)} \times \frac{G_2\left(\frac{-1}{2}, \frac{\tau'}{46}, \frac{\sigma'}{46}, \frac{\rho'}{46}\right)^{-23}}{G_2\left(\frac{-23}{2}, \frac{\tau'}{2}, \frac{\sigma'}{2}, \frac{\rho'}{2}\right)}
$$
using formula \refp{exponentialformula} yields a number close to the same root $\approx 13.9102308... -i\cdot24.0932265...$ of the polynomial $x^8 - 28x^7 + 778x^6 - 112x^5 - 749x^4 - 112x^3 + 778x^2 - 28x + 1$. It is expected to be true generally that formula \refp{exponentialformula} may always be used to compute the right-hand side of formula \refp{eqconjecture} when some of the $\talpha_j$ are real.

\subsection{One Quintic example}

\stepcounter{example}

We have already successfully produced a dozen of quintic examples and we give here a detailed presentation of our simplest one. Let $\zexc$ be the complex root of the polynomial $x^5 - x^4 - x^3 - 2x^2 + x + 1$ lying in the upper half-plane. Then $\kk_{\exc} = \qq(\zexc)$ has class number $1$. We choose $\goth{f} = \goth{P}_3$ the degree one prime above $q = 3$. The corresponding narrow ray class group is $Cl^{+}(\goth{f}) \simeq \zsz{2}$. We choose the fundamental units (once again take the fundamental units provided by Pari/GP version 2.15.4 and then compute the simplest totally positive ones):
$$ \eps_1 = \zexc^4 - 2\zexc^3 - \zexc + 3,~\eps_2 = 2\zexc^4 - 2\zexc^3 - 3\zexc + 3,~\eps_3 = 2\zexc^4 - 3\zexc^3 - 4\zexc + 4 $$
for $\opc{\goth{f}}$. We fix the order $\{ Id, (32), (21), (231), (312), (31)\}$ of $\goth{S}_3$ given by Pari/GP. We compute $\tlambda_1 = \tlambda_2 = \tlambda_3 = \tlambda_4 = \tlambda_5 = \tlambda_6 = 1$, $\ttt_1 = 7\cdot 37 \cdot 137$, $\ttt_2 = 31\cdot53$, $\ttt_3 = 491$, $\ttt_4 = 107$, $\ttt_5 = 1$ and $\ttt_6 = 145637$ where we have written the prime number decomposition of all the overflows $\ttt$. We may choose $\goth{b} = (1)$ and $\goth{a}$ the unique degree one prime above $N = 11$ in $\kk_{\exc}$. We look for target ideals above the prime numbers listed above and find the ideals $\goth{P}_7, \goth{P}'_{37}, \goth{P}'_{137}, \goth{P}_{31}, \goth{P}_{53}, \goth{P}'_{491}, \goth{P}'_{107}, \goth{P}_{145637}$ in the order provided by Pari/GP version 2.15.4. The ideals $qN\goth{P}_7\goth{P}'_{37}\goth{P}'_{137}/\goth{a}$, $qN\goth{P}_{31} \goth{P}_{53}/\goth{a}$, $qN\goth{P}'_{491}/\goth{a}$, $qN\goth{P}'_{107}/\goth{a}$, $qN/\goth{a}$, $qN\goth{P}_{145637}/\goth{a}$ admit the following admissible generators:
\begin{align*}
h_1 & = 147\zexc^4 - 135\zexc^3 - 90\zexc^2 - 234\zexc + 72 \\
h_2 & = -15\zexc^4 - 30\zexc^3 + 90\zexc^2 + 36\zexc + 60  \\
h_ 3 & = 42\zexc^4 - 15\zexc^3 - 87\zexc^2 - 48\zexc + 30  \\
h_4 & = -21\zexc^4 + 57\zexc^3 - 6\zexc^2 - 9\zexc - 48  \\
h_5 & = -3\zexc^4 - 6\zexc^3 + 18\zexc^2 - 6\zexc + 12  \\
h_6 & = 108\zexc^4 - 246\zexc^3 + 111\zexc^2 - 81\zexc + 195 
\end{align*}

\noindent The corresponding levels are $l_1 = 3\cdot7\cdot11\cdot37\cdot137$, $l_2 = 3\cdot 11\cdot 491$, $l_3 = 3\cdot11\cdot107$, $l_4 = 3\cdot11\cdot31\cdot53$, $l_5 = 3\cdot11$, $l_6 = 3\cdot11\cdot145637$. Here we have $\zfone = \{(1,1), (2,-1)\}$ and $\kappa = 2$. We compute the six quotients: 
\allowdisplaybreaks
\begin{align*}
&\frac{G_3\left(\frac{-1}{3}, \frac{\tau_1}{l_1}, \frac{\sigma_1}{l_1}, \frac{\rho_1}{l_1}, \frac{\varpi_1}{l_1}\right)^{11}}{G_3\left(\frac{-11}{3}, \frac{11\tau_1}{l_1}, \frac{11\sigma_1}{l_1}, \frac{11\rho_1}{l_1}, \frac{11\varpi_1}{l_1}\right)} && \approx -6.9846353...10^{-18} + i\cdot6.2764063... 10^{-17} \\
&\frac{G_3\left(\frac{-1}{3}, \frac{\tau_2}{l_2}, \frac{\sigma_2}{l_2}, \frac{\rho_2}{l_2}, \frac{\varpi_2}{l_2}\right)^{-11}}{G_3\left(\frac{-11}{3}, \frac{11\tau_2}{l_2}, \frac{11\sigma_2}{l_2}, \frac{11\rho_2}{l_2}, \frac{11\varpi_1}{l_2}\right)^{-1}} && \approx
1.7434761...10^{11} -i\cdot1.7914686...10^{11} \\
&\frac{G_3\left(\frac{1}{3}, \frac{\tau_3}{l_3}, \frac{\sigma_3}{l_3}, \frac{\rho_3}{l_3}, \frac{\varpi_3}{l_3}\right)^{-11}}{G_3\left(\frac{11}{3}, \frac{11\tau_3}{l_3}, \frac{11\sigma_3}{l_3}, \frac{11\rho_3}{l_3}, \frac{11\varpi_3}{l_3}\right)^{-1}} && \approx
0.2815024... - i\cdot0.0665163... \\
&\frac{G_3\left(\frac{1}{3}, \frac{\tau_4}{l_4}, \frac{\sigma_4}{l_4}, \frac{\rho_4}{l_4}, \frac{\varpi_4}{l_4}\right)^{11}}{G_3\left(\frac{11}{3}, \frac{11\tau_4}{l_4}, \frac{11\sigma_4}{l_4}, \frac{11\rho_4}{l_4}, \frac{11\varpi_4}{l_4}\right)} && \approx 9.8457700... - i\cdot2.3506772...\\
&\frac{G_3\left(\frac{-1}{3}, \frac{\tau_5}{l_5}, \frac{\sigma_5}{l_5}, \frac{\rho_5}{l_5}, \frac{\varpi_5}{l_5}\right)^{11}}{G_3\left(\frac{-11}{3}, \frac{11\tau_5}{l_5}, \frac{11\sigma_5}{l_5}, \frac{11\rho_5}{l_5}, \frac{11\varpi_5}{l_5}\right)} && \approx -0.2518907... -i\cdot 0.2274847...\\
&\frac{G_3\left(\frac{-1}{3}, \frac{\tau_6}{l_6}, \frac{\sigma_6}{l_6}, \frac{\rho_6}{l_6}, \frac{\varpi_6}{l_6}\right)^{-11}}{G_3\left(\frac{-11}{3}, \frac{11\tau_6}{l_6}, \frac{11\sigma_6}{l_6}, \frac{11\rho_6}{l_6}, \frac{11\varpi_6}{l_6}\right)^{-1}} && \approx 100864.0193260... - i\cdot767246.5816458...
\end{align*}

\noindent where the parameters are given in the table below. The product of these six quotients is close to the root $\approx -11.6360077... + i\cdot3.4634701...$ of the polynomial $x^{10} + 24x^9 + 164x^8 + 99x^7 - 62x^6 - 89x^5 - 62x^4 + 99x^3 + 164x^2 + 24x + 1$ which defines $\kk_{\exc}^{+}(\goth{f})$. The computation time for $1000$ digits is 1 minute and 35 seconds, but the computation time for each of the individual computations is not uniform at all. The fifth computation requires 1 second whereas the sixth computation requires 58 seconds because of the level difference $l_5 = 33$ versus $l_6 = 4806021$.
\bigskip

\noindent \begin{tabular}{| C C C C C |}
\hline
\eps_1\tau_1 \hfill= &229\eps_1\eps_2\eps_3 &-1727\eps_1\eps_2 &-796386\eps_1&+3385 \\
\eps_1\eps_2\sigma_1 \hfill= &-114\eps_1\eps_2\eps_3 &-987551\eps_1\eps_2 &+718\eps_1 &+949 \\
\eps_1\eps_2\eps_3\rho_1 \hfill = &-595138\eps_1\eps_2\eps_3 +&12 \eps_1\eps_2 &+1445\eps_1 &-907 \\
\varpi_1 \hfill = &-52 \eps_1\eps_2\eps_3 &+857\eps_1\eps_2 &-295\eps_1 &+2067162\\
\hline

\eps_1\tau_2 \hfill = &99\eps_1\eps_2\eps_3 &-356\eps_1\eps_3 &+19304\eps_1 &+57 \\
\eps_1\eps_3\sigma_2 \hfill = &-17\eps_1\eps_2\eps_3 &-70870\eps_1\eps_3 &-145\eps_1 &+40 \\
\eps_1\eps_2\eps_3\rho_2 \hfill = &324170\eps_1\eps_2\eps_3 &-85 \eps_1\eps_3 &+197\eps_1 &-111 \\
\varpi_2 \hfill = &91 \eps_1\eps_2\eps_3 &-460\eps_1\eps_3 &+3\eps_1 &-28725\\
\hline

\eps_1\eps_2\tau_4 \hfill = &-62\eps_1\eps_2\eps_3 &-12636\eps_1\eps_2 &+20\eps_2 &+3 \\
\eps_1\eps_2\eps_3\sigma_4 \hfill = &327915\eps_1\eps_2\eps_3 &-37\eps_1\eps_2 &-1101\eps_2 &+154 \\
\rho_4 \hfill = &-44\eps_1\eps_2\eps_3 &+314 \eps_1\eps_2 &-65\eps_2 &-52915 \\
\eps_2\varpi_4 \hfill = &-262 \eps_1\eps_2\eps_3 &-5\eps_1\eps_2 &-768869\eps_2 &-19\\
\hline

\eps_2\eps_3\tau_3 \hfill = &-27\eps_1\eps_2\eps_3 &-84\eps_2\eps_3 &+11\eps_2 &+1 \\
\eps_1\eps_2\eps_3\sigma_3 \hfill = &2555\eps_1\eps_2\eps_3 &+19\eps_2\eps_3 &-391\eps_2 &+52 \\
\rho_3 \hfill = &-11\eps_1\eps_2\eps_3 &+49\eps_2\eps_3 &-51\eps_2 &-491 \\
\eps_2\varpi_3 \hfill = &-80 \eps_1\eps_2\eps_3 &-23\eps_2\eps_3 &-1402\eps_2 &-1\\
\hline

\eps_1\eps_2\eps_3\tau_5 \hfill = &-456\eps_1\eps_2\eps_3 &-7\eps_1\eps_3 &-2\eps_3 &+3 \\
\sigma_5 \hfill = &-4\eps_1\eps_2\eps_3 &+17\eps_1\eps_3 &+12\eps_3 &+566 \\
\eps_3\rho_5 \hfill = &9\eps_1\eps_2\eps_3 &-39 \eps_1\eps_3 &+1265\eps_3 &+13 \\
\eps_1\eps_3\varpi_5 \hfill = &\eps_1\eps_2\eps_3 &-252\eps_1\eps_3 &-3\eps_3 &+2\\
\hline

\eps_1\eps_2\eps_3\tau_6 \hfill = &9785210\eps_1\eps_2\eps_3 &+2609\eps_2\eps_3 &-523\eps_3 &-729 \\
\sigma_6 \hfill = &479\eps_1\eps_2\eps_3 &+987\eps_2\eps_3 &-13874\eps_3 &-19886676 \\
\eps_3\rho_6 \hfill = &1713\eps_1\eps_2\eps_3 &-8328 \eps_2\eps_3 &-13204073\eps_3 &+764 \\
\eps_2\eps_3\varpi_6 \hfill = &4002 \eps_1\eps_2\eps_3 &+7773526\eps_2\eps_3 &+1749\eps_3 &-2296\\
\hline
\end{tabular} \medskip







\subsection{Quartic examples with $\kappa < \cardinalshort{\zfone}$}\label{lowkappaex}

We now present three quartic examples with $\kappa < \cardinalshort{\zfone}$. The first example is a case where $q= 5$, $\kappa = 2$ and $\zfone \simeq \zsz{5}^{\times}$. The second example is a case where $q=5$, $\kappa = 1$ and $\zfone \simeq \zsz{5}^{\times}$. The third example is a case where $q = 7$, $\kappa = 1$ and $\zfone = \{(1,1), (6, -1)\}$.
 
\stepcounter{example}
\subsubsection{Example \exc}\label{kappa2}

We first discuss a quartic example with $q= 5$, $\kappa = 2$ and $\zfone \simeq \zsz{5}^{\times}$, where the character $\chi_{\goth{f}}$ is $(\frac{\cdot}{5})$. Let $\zexc$ be the complex root of the polynomial $x^4 - 3x + 1$ lying in the upper half-plane. Then $\kk_{\exc} = \qq(\zexc)$ has class number $1$. We choose $\goth{f} = \goth{P}_5$ the degree one prime above $q = 5$. The corresponding narrow ray class group is $Cl^{+}(\goth{f}) \simeq \zsz{2}$. We choose the fundamental units
$$ \eps_1 = -2\zexc^3 + \zexc^2 +3\zexc -1,~\eps_2 = \zexc^2 -2\zexc +1 $$
for $\opc{\goth{f}}$. We compute $\tlambda_1 = \tlambda_2 = 1$, $\ttt_1 = 25561$, $\ttt_2 = 59$ where $25561$ is a prime number. We may choose $\goth{b} = (1)$ and $\goth{a}$ the unique degree one prime above $N = 7$ in $\kk_{\exc}$. We search for target ideals and find $\goth{P}_{25561}$ in the first case and $\goth{P}_{59}$ in the second case. Both ideals $qN\goth{P}_{25561}/\goth{a}$ and $qN\goth{P}_{59}/\goth{a}$ are principal with admissible generators $h_1 = 130\zexc^3 - 40\zexc^2 - 125\zexc- 225$ and $h_2 = 65\zexc^3 + 50\zexc^2 + 60\zexc - 95$. The corresponding levels will be $l_1 = 5\cdot7\cdot25561$ and $l_2 = 5\cdot7\cdot59$. Let us write the parameters

\begin{align*}
&\tau = \frac{9\eps_1\eps_2+694144\eps_1+70}{\eps_1}, && \tau' = \frac{-19730\eps_2\eps_1+11\eps_2-4}{\eps_1\eps_2}  \\
&\sigma = \frac{-3353996\eps_1\eps_2+504\eps_1-215}{\eps_1\eps_2}, && \sigma' = \frac{-215\eps_2\eps_1-2436\eps_2+64874}{1}  \\
&\rho = \frac{32\eps_1\eps_2+11\eps_1+130894}{1}, && \rho' = \frac{-5\eps_2\eps_1+5654\eps_2+9}{\eps_2}
\end{align*}

\noindent Then the two products 
\begin{align*}
v_1 = &\prod_{k = 1, 4}\frac{G_2\left(\frac{k}{5}, \frac{\tau}{894635}, \frac{\sigma}{894635}, \frac{\rho}{894635}\right)^{7}}{G_2\left(\frac{7k}{5}, \frac{\tau}{127805}, \frac{\sigma}{127805}, \frac{\rho}{127805}\right)}\times\frac{G_2\left(\frac{k}{5}, \frac{\tau'}{2065}, \frac{\sigma'}{2065}, \frac{\rho'}{2065}\right)^{-7}}{G_2\left(\frac{7k}{5}, \frac{\tau'}{295}, \frac{\sigma'}{295}, \frac{\rho'}{295}\right)^{-1}}\\
v_2 = &\prod_{k = 2, 3}\frac{G_2\left(\frac{-k}{5}, \frac{\tau}{894635}, \frac{\sigma}{894635}, \frac{\rho}{894635}\right)^{-7}}{G_2\left(\frac{-7k}{5}, \frac{\tau}{127805}, \frac{\sigma}{127805}, \frac{\rho}{127805}\right)^{-1}}\times\frac{G_2\left(\frac{-k}{5}, \frac{\tau'}{2065}, \frac{\sigma'}{2065}, \frac{\rho'}{2065}\right)^{7}}{G_2\left(\frac{-7k}{5}, \frac{\tau'}{295}, \frac{\sigma'}{295}, \frac{\rho'}{295}\right)^{1}}
\end{align*}
coincide up to 1000 digits and are close to the root $\approx -33.4405123... -  i\cdot62.3162955...$ of the polynomial $x^8 + 68x^7 + 5078x^6 + 5703x^5 + 7945x^4 + 5703x^3 + 5078x^2 + 68x + 1$ which defines $\kk_{\exc}^{+}(\goth{f})$. The computation time for $1000$ digits is 48 minutes. The fact that the two products seem to be equal suggests that the integer $\kappa$ defined in section \ref{formulation} is equal to $2$ in this case. We then check formula \refp{eqconjecture} as:
$$\norm{\goth{a}}\zeta'_{\goth{f}}([\goth{b}], 0) - \zeta'_{\goth{f}}([\goth{a}\goth{b}], 0) \approx \frac{1}{4}\log\left|v_1v_2\right|^2 \approx 4.2587554...$$

\stepcounter{example}
\subsubsection{Example \exc}
We now discuss a case where $q = 5$, $\kappa = 1$ and $\zfone \simeq \zsz{5}^{\times}$ with $\chi_{\goth{f}} = 1$. In this case all $L$-functions attached to $\kk^{+}(\goth{f})$ have a zero of order at least $2$ at $0$ and $\kk^{+}(\goth{f})$ is the narrow Hilbert class field of $\kk$, therefore the Stark unit we aim to compute is just $1$. We consider $\zexc$ the complex root of the polynomial $x^4 - x^3 - 2x^2 - 3x - 2$ lying in the upper half-plane. Then $\kk_{\exc} = \qq(\zexc)$ has class number $1$. We choose $\goth{f} = \goth{P}_5$ the degree one prime above $q = 5$. The corresponding narrow ray class group is $Cl^{+}(\goth{f}) \simeq \zsz{2}$. We choose the fundamental units
$$\eps_1 = -\zexc^2 + 3\zexc + 3 ,~\eps_2 = \zexc^3 + 2\zexc^2 - 8\zexc - 7 $$
for $\opc{\goth{f}}$. We compute $\tlambda_1 = \tlambda_2 = 1$, $\ttt_1 = 2\cdot12907$, $\ttt_2 = 61\cdot73$ where $12907$ is a prime number. We may choose $\goth{b} = (1)$ and $\goth{a}$ the unique degree one prime above $N = 7$ in $\kk_{\exc}$. We search for target ideals and find $\goth{P}_2$ a degree one prime above $2$ and $\goth{P}'_{12907}$ a degree one prime above $12907$ in the first case. In the second case, we find $\goth{P}_{61}$ a degree one prime ideal above $61$ and $\goth{P}_{73}$ a degree one prime above $73$. Both ideals $qN\goth{P}_2\goth{P}'_{12907}/\goth{a}$ and $qN\goth{P}_{61}\goth{P}_{73}/\goth{a}$ are principal with admissible generators $h_1 = 40\zexc^3 - 210\zexc^2 + 95\zexc + 80$ and $h_2 = 360\zexc^3 - 525\zexc^2 - 160\zexc - 785$. The corresponding levels will be $l_1 = 2\cdot5\cdot7\cdot12907$ and $l_2 = 5\cdot7\cdot61\cdot73$. Let us write the parameters

\begin{align*}
&\tau = \frac{-59\eps_1\eps_2-194580149\eps_1+3255}{\eps_1}, && \tau' = \frac{-159315\eps_2\eps_1+1189\eps_2-76}{\eps_1\eps_2} \\
&\sigma = \frac{44607037\eps_1\eps_2+135\eps_1-49}{\eps_1 \eps_2}, && \sigma' = \frac{-49\eps_2\eps_1-108\eps_2-106940}{1} \\
&\rho = \frac{95\eps_1\eps_2+1189\eps_1+3636283}{1}, && \rho' = \frac{-108\eps_2\eps_1+36840\eps_2-59}{\eps_2}
\end{align*}

\noindent Then 
$$ \prod_{k = 1}^4 \frac{G_2\left(\frac{-k}{5}, \frac{\tau}{903490}, \frac{\sigma}{903490}, \frac{\rho}{903490}\right)^{7}}{G_2\left(\frac{-7k}{5}, \frac{\tau}{129070}, \frac{\sigma}{129070}, \frac{\rho}{129070}\right)}  \times\frac{G_2\left(\frac{-k}{5}, \frac{\tau'}{155855}, \frac{\sigma'}{155855}, \frac{\rho'}{155855}\right)^{-7}}{G_2\left(\frac{-7k}{5}, \frac{\tau'}{22265}, \frac{\sigma'}{22265}, \frac{\rho'}{22265}\right)^{-1}}
$$
is $\approx 1$ which is a correct Stark unit when all $L$ functions have a zero of order $2$ at $0$. The computation time for 1000 digits is 70 minutes. Computations show that the above product is not a power of one of its subproducts, therefore the integer $\kappa$ defined in section \ref{formulation} is equal to $1$. Amongst our examples, this one exhibits one of the largest computed value for $\cardinalshort{\zfone}/\kappa = 4$.

\stepcounter{example}
\subsubsection{Example \exc}
We now discuss a case where $q = 7$, $\kappa = 1$ and $\zfone = \{(1,1), (6, -1)\}$. Consider $\zexc$ the complex root of the polynomial $x^4 - x^3 - x^2 - 5x + 1$ lying in the upper half-plane. Then $\kk_{\exc} = \qq(\zexc)$ has class number $1$. We choose $\goth{f} = \goth{P}_7$ the degree one prime above $q = 7$ in $\kk_{\exc}$. The corresponding narrow ray class group is $Cl^{+}(\goth{f}) \simeq \zsz{6}$. We choose the fundamental units
$$\eps_1 = 2\zexc^3 - 5\zexc^2 + \zexc ,~\eps_2 = -\zexc^3 + \zexc^2 + 3\zexc$$
for $\opc{\goth{f}}$. We compute $\tlambda_1 = \tlambda_2 = 1$, $\ttt_1 = 13\cdot1069$, $\ttt_2 = 227257$ where $1069, 227257$ are prime numbers. We may choose $\goth{b} = (1)$ and $\goth{a} = \goth{P}_{11}$ the degree one prime above $N = 11$ in $\kk_{\exc}$. We search for target ideals and find $\goth{P}_{13}$ the degree one prime above $13$ and $\goth{P}'_{1069}$ a degree one prime above $1069$ in the first case. In the second case, we find $\goth{P}_{227257}$ a degree one prime above $227257$ in the second case. In Proposition \ref{propfinal} we may choose $m_1 = 1$, $m_2 = 2$ and $\goth{I} = \ok$ such that both ideals $qN\goth{P}_{13}\goth{P}'_{1069}/\goth{a}$ and $2qN\goth{P}_{227257}/\goth{a}$ are principal with admissible generators $h_1 = 182\zexc^3 + 35\zexc^2 - 84\zexc - 714$ and $h_2 = 532\zexc^3 - 490\zexc^2 + 406\zexc - 1554$. The corresponding levels will be $l_1 = 7\cdot11\cdot13\cdot1069$ and $l_2 = 7\cdot11\cdot227257$. Let us write the parameters

\begin{align*}
&\tau = \frac{5\eps_1\eps_2-5863655\eps_1+51}{\eps_1}, && \tau' = \frac{11916518\eps_2\eps_1+276\eps_2-191}{\eps_1\eps_2} \\
&\sigma = \frac{18641035\eps_1\eps_2+692\eps_1-197}{\eps_1 \eps_2}, && \sigma' = \frac{197\eps_2\eps_1+4832\eps_2+7904541}{1} \\
&\rho = \frac{-49\eps_1\eps_2-276\eps_1+2381446}{1}, && \rho' = \frac{216\eps_2\eps_1-1644334\eps_2-5}{\eps_2}
\end{align*}

\noindent For $k = 1, 2, 3, 4, 5, 6$, write $I_k$ for the value
$$ I_k = \frac{G_2\left(\frac{k}{7}, \frac{\tau}{1070069}, \frac{\sigma}{1070069}, \frac{\rho}{1070069}\right)^{11}}{G_2\left(\frac{11k}{7}, \frac{\tau}{97279}, \frac{\sigma}{97279}, \frac{\rho}{97279}\right)} \times \frac{G_2\left(\frac{-2k}{7}, \frac{\tau'}{17498789}, \frac{\sigma'}{17498789}, \frac{\rho'}{17498789}\right)^{-11}}{G_2\left(\frac{-22k}{7}, \frac{\tau'}{1590799}, \frac{\sigma'}{1590799}, \frac{\rho'}{1590799}\right)^{-1}}
$$
On the one hand, $I_1I_6 \approx -1197160922532.7652149... - i\cdot 455052155255.6575547... $ is the product in the right-hand side of formula \refp{eqconjecture}, which we check as:
$$\norm{\goth{a}}\zeta'_{\goth{f}}([\goth{b}], 0) - \zeta'_{\goth{f}}([\goth{a}\goth{b}], 0) \approx \frac{1}{2}\log\left|I_1I_6\right|^2 \approx 27.8784505...$$
On the other hand, the values $I_1I_6$, $I_2I_5 \approx 10^{-4}(6.0177934...-i\cdot6.7709930...)$ and $I_3I_4 \approx 10^{-10}(-3.0632199...-i\cdot 8.0567091...)$ seem to be algebraically dependent and we may identify the polynomial
\begin{align*}
\prod_{k = 1, 2, 3} (X - I_k I_{7-k})&(X-I_k^{-1} I_{7-k}^{-1}) = X^6 + 1  \\ &+(-153504836673\zexc^3 + 136635341749\zexc^2 \\&~~~+ 544161498400\zexc - 108155523221) (X^5+X) \\ &+(45202964386965263585\zexc^3 - 350943983942803910564\zexc^2 \\ &~~~+613535010394689372486\zexc - 104952772799918427079) (X^4 + X^2) \\ &+(138218090908063915784215\zexc^3 + 7152666826995541968563\zexc^2 \\ &~~~- 814169082025891017196278\zexc + 154702240787650756485938) X^3
\end{align*}


\noindent as a relative polynomial defining the extension $\kk_{\exc}^{+}(\goth{f})/\kk_{\exc}$. Alternatively, there are roots $u_k^{\pm 1}$ of the palindromic polynomial
\begin{align*}
&x^{24} +1+ 1292(x^{23} + x) + 610850(x^{22} + x^2) - 324070132(x^{21} +x^3) \\
&+ 43597528236(x^{20} +x^4) - 641906806311(x^{19} + x^5) \\
&+ 3661996071192(x^{18} + x^6) - 9521279678257(x^{17} +x^7) \\
&+ 11598437454383(x^{16} + x^8) - 3951777294409(x^{15} +x^9) \\
&- 4643997885968(x^{14} + x^{10}) + 3477156436009(x^{13} +x^{11}) - 43787701643x^{12}
\end{align*} defining $\kk^{+}_{\exc}(\goth{f})$ such that for $k = 1, 2, 3$, $(I_kI_{7-k})^{\pm 1} \approx u_k^{\pm 4}$ and it seems that $I_kI_{7-k}$ is already a $4$-th power inside $\kk^{+}_{\exc}(\goth{f})$.

\appendix

\section{Appendix: Pure cubic table for $2 \leq m \leq 50$}\label{cubicappendix}

The following table gives one of the positive fundamental units $\eps$ of the pure cubic field defined by $X^3-m$ in terms of the root $z = e^{2i\pi/3}m^{1/3}$ for cube-free integers $2 \leq m \leq 50$ up to isomorphism. It also gives the corresponding contents $\tlambda$ and the overflows $\ttt$ computed regardless of whether or not there is a modulus $\goth{f}$ satisfying the conditions in Conjecture \ref{conjecture} with $\opcf$ generated by $\eps$. In particular, it is clear that the parameters $\tlambda, \ttt$ vary significantly with the value of $m$. In general the value of $\tlambda$ stays rather small while the value of $\ttt$ may become very large. Our work shows that the computation time for the right-hand side of formula \refp{eqconjecture} is at least $O(qN\tlambda \ttt |\log \delta|)$ for the precision $\delta$ so computations can quickly become overwhelming, in particular for the cases $m = 23, 29, 33, 41, 46, 47$.

\begin{center}
\begin{tabular}{| C | C | C | C |}
\hline
X^3-m & \eps & \tlambda & \ttt\\
\hline
X^3-2 & z - 1 & 1 & 1\\
\hline
X^3 - 3& z^2 - 2& 1& 3\\
\hline
X^3 - 5& 2z^2 - 4z + 1& 2& 2\cdot13\\
\hline
X^3 - 6& -3z^2 + 6z - 1& 3& 2\cdot3\cdot7\\
\hline
X^3 - 7& -z + 2& 1& 1\\
\hline
X^3 - 10& 2/3z^2 - 1/3z - 7/3& 1& 3^2\\
\hline
X^3 - 11& -2z^2 + 4z + 1& 2& 2\cdot19 \\
\hline
X^3 - 12& -3/2z^2 + 3z + 1& 3& 3\cdot5 \\
\hline
X^3 - 13& 2z^2 - 3z - 4& 1& 131\\
\hline
X^3 - 14& -z^2 + 2z + 1& 1& 2\cdot11\\
\hline
X^3 - 15& -12z^2 + 30z - 1& 6& 2\cdot3\cdot5\cdot7^2\\
\hline
X^3 - 17& -7z + 18& 7& 3\cdot7 \\
\hline
X^3 - 19& 1/3z^2 - 2/3z - 2/3& 1& 3\\
\hline
X^3 - 20& -1/2z^2 + z + 1& 1& 7\\
\hline
X^3 - 21& -4z^2 - 6z + 47& 2& 2\cdot3\cdot47\\
\hline
X^3 - 22& 4z^2 - 3z - 23& 1& 5\cdot7\cdot41 \\
\hline
X^3 - 23& 6230z^2 - 3160z - 41399& 10& 2\cdot5\cdot31\cdot2659\cdot67853 \\
\hline
X^3 - 26& z - 3& 1& 3\\
\hline
X^3 - 28& 1/6z^2 - 1/3z - 1/3& 1& 1 \\
\hline
X^3 - 29& 370284z^2 + 103819462z & 2& 2\cdot7\cdot1953569\cdot16988921\\
 & - 322461439& &\cdot602078573\\
\hline
X^3 - 30& -3z^2 + 9z + 1& 3& 3^2\cdot19 \\
\hline
X^3 - 31& 20z^2 + 54z - 367& 2& 2\cdot11317\\
\hline
X^3 - 33& -394098z^2 + 97392z & 6& 2\cdot3\cdot21313\\
 & + 3742201& &\cdot224363\cdot1956481\\
\hline
X^3 - 34& 51z^2 + 24z - 613& 3& 2\cdot3^2\cdot5\cdot7\cdot13\cdot61\\
\hline
X^3 - 35& -1/3z^2 + 10/3z - 22/3& 1& 5\cdot23 \\
\hline
X^3 - 37& -3z + 10& 3& 3^2 \\
\hline
X^3 - 38& 3z^2 - 55z + 151& 1& 13\cdot79\cdot163\\ 
\hline
X^3 - 39& 2z^2 - 23& 2& 2\cdot3\cdot13 \\
\hline
X^3 - 41& -70761183382z^2 & 2& 2\cdot7\cdot17\cdot967\cdot17467\\
 & + 305478475184z & & \cdot810791\cdot126525109\\
 &-211991370839& & \cdot26087689931\\
\hline
X^3 - 42& -12z^2 + 42z - 1& 6& 2\cdot3\cdot7\cdot97\\
\hline
X^3 - 43& 2z - 7& 2& 2 \\
\hline
X^3 - 44& 17/6z^2 + 2/3z - 113/3& 1& 3^2\cdot23\cdot29\\
\hline
X^3 - 45& -104z^2 + 66z + 1081& 6& 2\cdot3\cdot37\cdot97\cdot197\\
\hline
X^3 - 46& -309z^2 - 48z + 4139& 9& 2\cdot3^2\cdot19\cdot281\cdot523\\
\hline
X^3 - 47& -64786z^2 + 69704z & 2& 2\cdot17\cdot37\cdot71\cdot2311\\
 &+ 592199 & &\cdot2791\cdot5693 \\
\hline
\end{tabular}
\end{center}

\bibliographystyle{plain}
\bibliography{bibliographie}{}

\end{document}